\documentclass[3p,preprint]{jcp}

\title{A Family of Independent Variable Eddington Factor Methods with \resp{Efficient Preconditioned Iterative Solvers}}
\author[1]{Samuel Olivier\corref{cor}}
\ead{solivier@berkeley.edu}
\author[2]{Will Pazner}
\ead{pazner1@llnl.gov}
\author[2]{Terry S. Haut}
\ead{haut3@llnl.gov}
\author[2]{Ben C. Yee}
\ead{yee26@llnl.gov}
\address[1]{Applied Science \& Technology, University of California, Berkeley, Berkeley, CA 94708, United States of America}
\address[2]{Lawrence Livermore National Laboratory, 7000 East Avenue, Livermore, CA 94550, United States of America}
\cortext[cor]{Corresponding author}

\begin{document}
\begin{abstract}
	We present a family of discretizations for the Variable Eddington Factor (VEF) equations that have high-order accuracy on curved meshes and efficient preconditioned iterative solvers. The VEF discretizations are combined with the Discontinuous Galerkin transport discretization from \cite{graph_sweeps} to form an effective high-order, linear transport method. The VEF discretizations are derived by extending the unified analysis of Discontinuous Galerkin methods for elliptic problems presented by \citet{Arnold2002} to the VEF equations. This framework is used to define analogs of the interior penalty, second method of Bassi and Rebay, minimal dissipation local Discontinuous Galerkin, and continuous finite element methods. The analysis of subspace correction preconditioners \cite{Pazner2021}, which use a continuous operator to iteratively precondition the discontinuous discretization, is extended to the case of the non-symmetric VEF system. Numerical results demonstrate that the VEF discretizations have arbitrary-order accuracy on curved meshes, preserve the thick diffusion limit, and are effective on a proxy problem from thermal radiative transfer in both outer transport iterations and inner preconditioned linear solver iterations. 
	\resp[red]{We demonstrate that the VEF solution converges to the \Sn transport solution as the mesh is refined on both problems with smooth and non-smooth behavior in angle. Parallel performance studies show that the interior penalty VEF discretization's linear solve weak scales out to 1024 processors and strong scales well on a single node. Particular attention is paid to the parallel performance of the VEF algorithm when used in combination with a parallel block Jacobi transport sweep.} 
\end{abstract}

\maketitle 

\section{Introduction}
\ifreview\linenumbers\fi
The Variable Eddington Factor (VEF) method \cite{mihalas}, also known as Quasidiffusion (QD) \cite{goldin}, is a rapidly converging, nonlinear scheme for solving the Boltzmann transport equation, a crucial component of high energy density physics (HEDP) simulations, nuclear reactor analysis, and medical physics. VEF has been applied to a wide range of transport and multiphysics problems including (but not limited to) nuclear reactor eigenvalue problems \cite{airstova_eigenvalue}, nuclear reactor kinetics \cite{doi:10.13182/NSE13-42}, and thermal radiative transfer (TRT) \cite{anistratov1996nonlinear}. It performs well in problems having both optically thick and thin regions and treats anisotropic scattering equally well \cite{anistratov_fvm,ARISTOVA2000139}. Robust convergence is achieved by iteratively coupling the transport equation to the VEF equations, a moment-based equivalent reformulation of transport. The exact closures used to form the VEF equations are weak functions of the solution meaning even simple iterative schemes, such as fixed-point iteration, can often converge in a small number of iterations that is independent of the mean free path. 

VEF offers significant algorithmic flexibility in that any valid discretization of the VEF equations will yield a rapidly converging algorithm. This is in stark contrast to Diffusion Synthetic Acceleration (DSA) which places severe restrictions on the discretization of the moment equations in order to guarantee stability \cite{A}. In the case where the VEF and transport discretizations are not algebraically consistent, referred to as a VEF method with an ``independent'' discretization \cite{doi:10.1080/00411459308203810,two-level-independent-warsa}, the discrete solutions of the transport and VEF equations will differ on the order of the discretization error and will be equivalent only in the limit as the mesh is refined. However, even in an under resolved problem, VEF still produces a ``transport solution'' in that the solution of the VEF method is a discrete solution of an equivalent reformulation of the transport equation. Furthermore, VEF methods generally preserve the thick diffusion limit \cite{diflim} and have conservation even if the transport discretization in isolation does not. These properties are particularly useful in multiphysics calculations since the lower-dimensional VEF equations can be directly coupled to the other physics components in place of the high-dimensional transport equation. In addition, discretizations for the transport and VEF equations can be designed independently so that they are in some sense optimal for their intended uses. 

This flexibility has been exploited to improve efficiency in relation to all seven dimensions of the transport equation. 
\citet{GHASSEMI2020109315} showed that different order temporal discretizations can be applied to the transport and VEF equations. Ongoing work suggests that time-stepping stability and accuracy can be maintained when just one transport inversion is performed per time step \cite{yee_mc21}. \citet{anistratov2021implicit} used data compression techniques to reduce storage costs in time-dependent calculations. In astrophysics, VEF is used to simplify the implementation of coupling TRT to hydrodynamics and to avoid the memory cost of solving the time-dependent transport equation \cite{Jiang_2012,GNEDIN2001437,GEHMEYR1994320}. 
\citet{Davis_2012} used a short characteristics discretization of the transport equation. \citet{me} and \citet{LOU2019258} designed a spatial discretization of the VEF equations to increase multiphysics compatibility. \citet{YEE2020109696} showed that robust convergence is maintained even when positivity-preserving methods are used inside the iteration. 
\citet{ANISTRATOV2019186} solved the multigroup TRT equations by using a VEF method with multiple levels in frequency. It is also well-known that the multigroup eigenvalue problem can be solved with only the need for eigenvalue iterations on the one-group VEF equations \cite{AL}. 

The above techniques rely on the efficient solution of the discretized VEF equations. VEF methods reduce the overall cost of the simulation by trading inversions of the high-dimensional transport equation for inversions of the lower-dimensional VEF equations. In all of VEF's applications, the inversion of the discretized VEF equations is buried under multiple nested loops corresponding to time integration, Newton iterations, eigenvalue iterations, multi-group iterations, and/or fixed-point iterations. The efficient iterative inversion of the VEF equations is then crucial to the efficiency of the overall algorithm and is a prerequisite for the practicality of any VEF method. 

The unusual structure of the VEF equations and their lack of self-adjointness make the development of discretizations and their corresponding preconditioned iterative solvers difficult.
While considerable effort has been placed into discretizing the VEF equations, to our knowledge, existing methods either rely on expensive and unscalable preconditioners such as block incomplete LU (BILU) factorization, cannot be solved with iteration counts independent of the mesh size, or do not mention solvers entirely. Previous work on discretizing the VEF equations includes finite volume \cite{anistratov_fvm,doi:10.1080/00411459308203810,QDBC,Jiang_2012,Jones2019TheQM}, finite difference \cite{WIESELQUIST2014343}, mixed finite element \cite{vallette,me,olivier_mandc,LOU2019258}, continuous finite element \cite{wieselquist,two-level-independent-warsa}, and discontinuous finite element \cite{dima_dfem} techniques. Most VEF methods are designed to be algebraically consistent with their application's discretized transport equation which typically requires discretizing the first-order form of the VEF equations. Such discretizations solve for both the zeroth and first moment of the solution and thus have significantly more unknowns than discretizations of the second-order form. In addition, block preconditioners \cite{benzi_golub_liesen_2005} are required to efficiently solve discretizations of the first-order form. Such solvers can require nested iteration for robustness (see \cite{warsa_mfem} for a radiation diffusion example). 

\citet{two-level-independent-warsa} showed that VEF methods with and without algebraic consistency converge equivalently as long as the transport data is properly represented. In particular, computing the Eddington tensor and boundary factor using finite element interpolation and Discrete Ordinates (\Sn) angular quadrature enables rapid convergence for any valid discretization of the VEF equations. An independent discretization of the second-order form of the VEF equations then has the potential to provide the rapid convergence of a consistent VEF method while solving for fewer unknowns and avoiding the need for block preconditioners. Such a method also has the flexibility to discretize the VEF equations in a manner that can leverage existing linear solver technology.

Our motivation for this research is in the context of HEDP experiments where the tightly coupled simulation of hydrodynamics and TRT is required, the latter of which typically includes the \Sn transport equation. For hydrodynamics, it has been shown that, compared to low-order methods, high-order methods on curved meshes have improved robustness, symmetry preservation, and strong scaling on emerging high performance computer architectures \cite{blast,blast2,blast3}. Transport methods compatible with this multiphysics framework are desired. \citet{graph_sweeps} showed that adequately approximating realistic meshes generated from a high-order hydrodynamics code as straight-edged required a significant number of mesh refinements leading to an impractical increase in transport unknowns. It is also possible that high-order accurate transport methods could be beneficial in terms of multiphysics compatibility with high-order hydrodynamics. High-order transport methods compatible with curved meshes have been developed recently in \cite{woods_thesis,graph_sweeps} with corresponding consistent DSA discretizations in \cite{ldrd_dsa,doi:10.1080/00295639.2020.1799603}. However, high-order discretizations of the VEF equations compatible with curved meshes have not yet been developed. 

In this paper, we design a family of independent VEF discretizations for the linear, steady-state transport problem that can be efficiently and scalably solved with high-order accuracy, in multiple dimensions, and on curved meshes. Our approach is to begin with discretization techniques known to have effective preconditioners on the simpler case of radiation diffusion (i.e.~the model Poisson problem) and adapt them to the VEF equations. 
By using the Eddington tensor and boundary factor interpolation procedure established in \cite{two-level-independent-warsa}, these methods achieve both rapid convergence in outer fixed-point iterations and in inner linear solver iterations when paired with a high-order Discontinuous Galerkin (DG) discretization of \Sn transport. 

In particular, we extend the unified analysis of DG methods developed for elliptic problems presented by \citet{Arnold2002} to the VEF equations to derive analogs of the interior penalty (IP), second method of Bassi and Rebay (BR2), minimal dissipation local Discontinuous Galerkin (MDLDG), and continuous finite element (CG) techniques.
We show that the IP and BR2 VEF methods are effectively preconditioned by the subspace correction method from \citet{Pazner2021} and that Algebraic Multigrid (AMG) is effective for the CG and MDLDG discretizations.
\citet{dima_dfem} also applied DG techniques to the VEF equations but they discretize the first-order form of the VEF equations and only consider lowest-order elements in one dimension.
We note that our CG operator is equivalent to extending the continuous finite element VEF discretization in \cite{two-level-independent-warsa} to multiple dimensions, arbitrary-order, and curved meshes.

The paper proceeds as follows.
First, we describe the VEF method analytically and discuss iterative schemes to solve the coupled transport-VEF system.
Then, we provide background on representing high-order meshes and finite element solutions and present the mathematical notation that will be used in the remainder of the paper.
We derive the extension of the unified framework for DG to the VEF equations.
The IP, BR2, MDLDG, and CG VEF discretizations are derived from this framework.
\S \ref{sec:subspace} discusses the design and analysis of subspace correction preconditioners and extends their analysis to the case of non-symmetric linear systems.

We next give computational results.
We show that all the methods presented achieve high-order accuracy on a curved mesh through the method of manufactured solutions, preserve the thick diffusion limit both on an orthogonal and a severely distorted curved mesh, and are effective on the linearized, steady-state crooked pipe problem, a challenging proxy problem from TRT, in both outer fixed-point iterations and inner linear solver iterations.
\resp[red]{Next, the IP VEF solution is shown to converge in space to the DG \Sn transport solution computed using the DSA preconditioner of \citet{ldrd_dsa} on both a problem with smooth and non-smooth solution in angle. 
We then present a parallel weak scaling study for the IP discretization which demonstrates the scalability of the algorithm out to 1024 processors and 40 million VEF scalar flux unknowns. This is followed by a strong scaling study showing the performance of the IP VEF method on a single node. The parallel scaling studies include an investigation of the performance consequences associated with using a parallel block Jacobi transport sweep.}
Finally, we give conclusions and recommendations for future work.

\section{The VEF Method}
The steady-state, mono-energetic, fixed-source transport problem with isotropic scattering and inflow boundary conditions is: 
	\begin{subequations} 
	\begin{equation} \label{eq:transport}
		\Omegahat\cdot\nabla\psi + \sigma_t \psi = \frac{\sigma_s}{4\pi} \int \psi \ud \Omega' + q \,, \quad \x \in \D \,,
	\end{equation}
	\begin{equation} \label{eq:transport_inflow}
		\psi(\x,\Omegahat) = f(\x,\Omegahat) \,, \quad \x \in \partial \D \ \mathrm{and} \ \Omegahat\cdot\hat{n} < 0 \,,
	\end{equation}
	\end{subequations}
where $\psi(\x,\Omegahat)$ is the angular flux, $\Omegahat\in \mathbb{S}^2$ the direction of particle flow, $\D$ the spatial domain of the problem with $\partial\D$ its boundary, $\sigma_t(\x)$ and $\sigma_s(\x)$ the total and scattering macroscopic cross sections, respectively, $q(\x,\Omegahat)$ the fixed-source, and $f(\x,\Omegahat)$ the inflow boundary function. The VEF equations are given by 
	\begin{subequations} \label{eq:vef_mixed}
	\begin{equation} \label{eq:vef0}
		\nabla\cdot\vec{J} + \sigma_a \varphi = Q_0 \,, 
	\end{equation}
	\begin{equation} \label{eq:vef1}
		\nabla\cdot\paren{\E\varphi} + \sigma_t\vec{J} = \Qone \,,
	\end{equation}
	\end{subequations}
where $\sigma_a(\x) = \sigma_t(\x) - \sigma_s(\x)$ is the absorption macroscopic cross section, $\varphi(\x)$ and $\vec{J}(\x)$ the zeroth and first angular moments of the angular flux, and 
	\begin{equation} \label{eq:eddington_tensor}
		\E(\x) = \frac{\int \Omegahat\otimes\Omegahat\, \psi \ud \Omega}{\int \psi \ud \Omega} 
	\end{equation}
is the Eddington tensor. We refer to $\varphi(\x)$ as the scalar flux and $\vec{J}(\x)$ as the current. In addition, $Q_i = \int \Omegahat^i q \ud \Omega$ are the angular moments of the fixed-source, $q$. The VEF equations are derived by taking the zeroth and first angular moments of the transport equation and closing the second moment of the angular flux, ${\P = \int \Omegahat\otimes\Omegahat\, \psi \ud \Omega}$, with 
	\begin{equation}
		\P = \E\varphi \,. 
	\end{equation}
By eliminating the current, the VEF equations can be cast as a drift-diffusion equation: 
	\begin{equation} \label{eq:vef_2nd_order}
		-\nabla\cdot \frac{1}{\sigma_t}\nabla\cdot\paren{\E\varphi} + \sigma_a \varphi = Q_0 - \nabla\cdot \frac{\Qone}{\sigma_t} \,. 
	\end{equation}
In both the first-order form (Eq.~\ref{eq:vef_mixed}) and second-order form (Eq.~\ref{eq:vef_2nd_order}), the presence of the Eddington tensor inside the divergence leads to diffusion, advection, and reaction-like terms that make applying existing discretization techniques difficult. 

The Miften-Larsen transport-consistent boundary conditions \cite{QDBC} are 
	\begin{equation} \label{eq:mlbc}
		\vec{J}\cdot\hat{n} = 2g + E_b \varphi \,, \quad \x \in \partial \D 
	\end{equation}
where 
	\begin{equation}
		g(\x) = \int_{\Omegahat\cdot\hat{n}<0} \Omegahat\cdot\hat{n} \, f(\x,\Omegahat) \ud \Omega 
	\end{equation}
is the incoming partial current computed from the transport boundary inflow function and 
	\begin{equation} \label{eq:eddington_bdr_factor}
		E_b = \frac{\int |\Omegahat\cdot\hat{n}|\, \psi \ud \Omega}{\int \psi \ud \Omega}
	\end{equation}
is the Eddington boundary factor. This boundary condition is derived by manipulating partial currents and using an analogous nonlinear closure. In equations, with the partial currents defined as $J_n^\pm = \int_{\Omegahatn\gtrless0} \Omegahatn\,\psi \ud \Omega$, 
	\begin{equation}
	\begin{aligned}
		\vec{J}\cdot\hat{n} &= J_n^+ + J_n^- \\
		&= \paren{J_n^+ - J_n^-} + 2J_n^- \\
		&= \int |\Omegahatn|\, \psi \ud \Omega + 2J_n^- \\
		&\rightarrow E_b \varphi + 2J_n^- \,, 
	\end{aligned}
	\end{equation}
where $g$ in Eq.~\ref{eq:mlbc} plays the role of $J_n^-$ using the transport equation's inflow boundary condition. 

If the Eddington tensor and boundary factor are known, the VEF equations define the zeroth and first moments of the angular flux. In other words, the VEF equations with Miften-Larsen boundary conditions are an equivalent reformulation of the transport equation. However, this is a trivial closure in that the solution to the transport equation must already be known to define the VEF data. VEF methods rely on the fact that the VEF data are weak functions of the angular flux and thus simple iterative schemes can converge rapidly. 

Note that when an independent discretization is used for the VEF equations, the discretized VEF scalar flux and VEF current will not be equivalent to the zeroth and first angular moments of the discrete angular flux; the two solutions will differ on the order of the spatial discretization error. To notationally separate the two scalar flux solutions, we use $\varphi$ (varphi) to denote the VEF scalar flux and $\phi = \int \psi \ud \Omega$ (phi) as the zeroth moment of the angular flux.

VEF methods seek the solution of the nonlinearly coupled system of equations: 
	\begin{subequations}
	\begin{equation}
		\Omegahat\cdot\nabla\psi + \sigma_t \psi = \frac{\sigma_s}{4\pi}\varphi + q \,, 
	\end{equation}
	\begin{equation}
		-\nabla\cdot \frac{1}{\sigma_t}\nabla\cdot\paren{\E\varphi} + \sigma_a \varphi = Q_0 - \nabla\cdot \frac{\Qone}{\sigma_t}\,,
	\end{equation}
	\end{subequations}
where the drift-diffusion form of VEF is used for brevity. Boundary conditions are specified by Eqs.~\ref{eq:transport_inflow} and \ref{eq:mlbc} for the transport and VEF drift-diffusion equation, respectively. Here, the transport equation's scattering source is now coupled to the VEF drift-diffusion equation and the data for the VEF drift-diffusion equation are nonlinearly coupled to the transport equation. We have increased the complexity of the problem by adding the VEF scalar flux as an additional unknown and by casting the linear transport problem as nonlinear. However, properties of the VEF data allow this nonlinear, coupled system to be solved more efficiently than algorithms based on the transport equation alone.  

Let 
	\begin{equation}
		\mat{L}\psi = \Omegahat\cdot\nabla\psi + \sigma_t \psi\,,
	\end{equation}
	\begin{equation}
		\mat{R}(\psi)\varphi = -\nabla\cdot \frac{1}{\sigma_t}\nabla\cdot\paren{\E(\psi)\varphi} + \sigma_a \varphi \,,
	\end{equation}
be the streaming and collision operator and VEF drift-diffusion operator, respectively, where $\paren{\cdot}$ indicates a nonlinear dependence on the argument. By linearly eliminating the angular flux, the transport-VEF system is equivalent to 
	\begin{equation} \label{eq:lin_elim}
		\mat{R}\!\paren{\mat{L}^{-1}\!\paren{\frac{\sigma_s}{4\pi}\varphi + q}}\!\varphi = Q_0 - \nabla\cdot \frac{\Qone}{\sigma_t} \,. 
	\end{equation}
Applying the inverse of the drift-diffusion operator, we see that the solution of the coupled transport-VEF system is the fixed-point: 
	\begin{equation} \label{eq:fixed_point}
		\varphi = G(\varphi) 
	\end{equation}
where 
	\begin{equation} \label{eq:fp_op}
		G(\varphi) = \mat{R}\!\paren{\mat{L}^{-1}\!\paren{\frac{\sigma_s}{4\pi}\varphi + q}}^{-1}\!\paren{Q_0 - \nabla\cdot \frac{\Qone}{\sigma_t}} \,. 
	\end{equation}
The fixed-point operator $G$ is applied in two stages: 1) solve the transport equation using a scattering source formed from the VEF  scalar flux and 2) solve the VEF drift-diffusion equation using the VEF data computed from the angular flux from stage 1). 

The simplest algorithm to solve Eq.~\ref{eq:fixed_point} is fixed-point iteration: 
	\begin{equation}
		\varphi^{k+1} = G(\varphi^k) 
	\end{equation}
where $k$ denotes the iteration index and $\varphi^0$ is an initial guess for the solution. This process is repeated until the difference between successive iterates is small enough. Since the Eddington tensor and boundary factor are weak functions of the angular flux even this simple iteration strategy often converges rapidly. 

Iterative efficiency can be improved with the use of Anderson acceleration. Anderson acceleration defines the next iterate as the linear combination of the previous $m$ iterates that minimizes the residual $\varphi - G(\varphi)$. For the storage cost of $m$ previous iterates, Anderson acceleration increases the convergence rate and improves robustness. While it is not practical to store multiple copies of the \emph{angular} flux, it is reasonable to expect that a small set of \emph{scalar} flux-sized vectors can be stored. The process of linearly eliminating the transport equation, codified in Eq.~\ref{eq:lin_elim}, allows the Anderson space to be built from the much smaller scalar flux-sized vectors only. In the case where a subset of the angular flux unknowns are not eliminated, such as when a parallel block Jacobi sweep is used to avoid communication costs or when mesh cycles or reentrant faces are present, the solution vector can be augmented with these un-eliminated unknowns so that they are included in the Anderson space. This is the nonlinear analog to the ideas used for Krylov-accelerated source iteration \cite{doi:10.13182/NSE02-14}. 

In addition, defining the nonlinear residual as 
	\begin{equation}
		F(\varphi) = \varphi - G(\varphi) = 0 \,,
	\end{equation}
root-finding methods such as Jacobian-free Newton Krylov (JFNK) can be used. JFNK builds a new Krylov space to approximate the gradient of $F$ at each iteration meaning information across iterations is not kept. JFNK typically required significantly more evaluations of $G$ than Anderson-accelerated fixed-point iteration. Thus, we present results using fixed-point iteration and Anderson-accelerated fixed-point iteration only. 

The following sections present the discretizations and solvers needed to efficiently evaluate $G$ numerically. 

\section{Mesh and Finite Element Preliminaries}
\subsection{Description of the Mesh}
\begin{figure}
	\centering
	\begin{subfigure}{.49\textwidth}
		\centering
		\includegraphics[width=.6\textwidth]{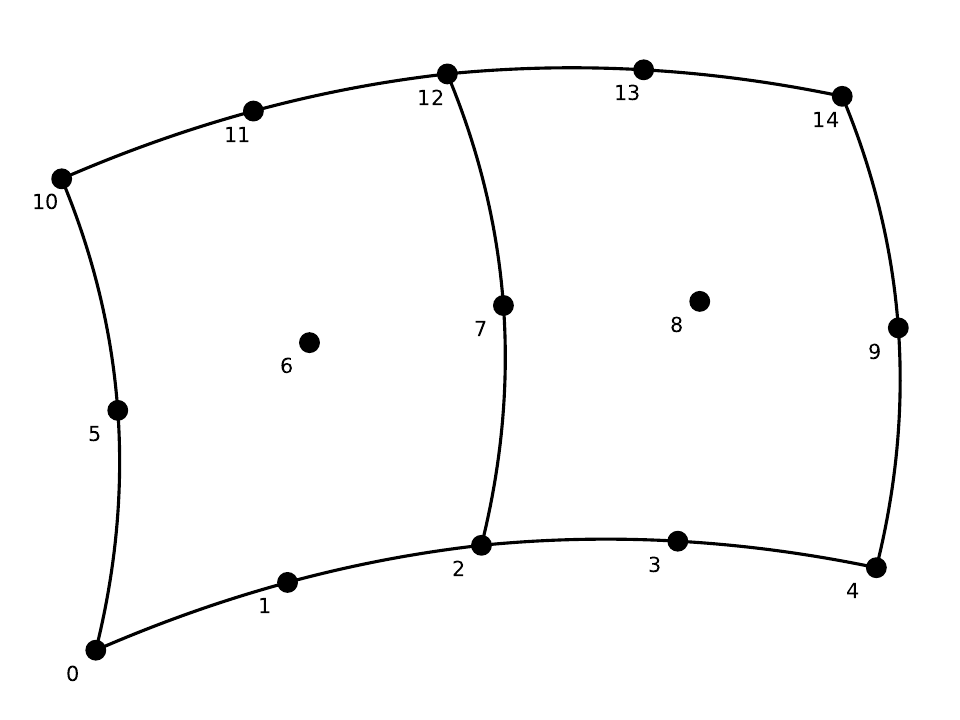}
		\caption{}
		\label{fig:curved_mesh}
	\end{subfigure}
	\begin{subfigure}{.49\textwidth}
		\centering 
		\includegraphics[width=\textwidth]{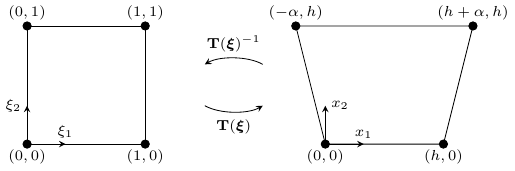}
		\caption{}
		\label{fig:trap_trans}
	\end{subfigure}
	\caption{Depictions of (a) the continuity of an interior face in a high-order curved mesh and (b) the reference transformation for a non-affine, linear, quadrilateral element.}
	\label{fig:mesh_fig}
\end{figure}

Let $\D \subset \R^{\dim}$ with $\dim=2,3$ be the domain of the problem. Consider the tessellation $$\D = \bigcup_{K_e\in\meshT} K_e$$ with $K_e$ the $e^{th}$ element in the mesh $\meshT$. Each coordinate of the mesh is represented by a piecewise continuous polynomial. In other words, the mesh itself is a member of an $\bracket{H^1(\D)}^{\dim}$ finite element space. This allows representation of  curved surfaces and enforces continuity of the mesh coordinates along the interfaces between elements. Figure \ref{fig:curved_mesh} depicts a mesh of two quadratic, quadrilateral elements where the mesh control points labeled 2, 7, and 12 are shared between the two elements to enforce continuity of the shared interior interface between them. 

The mesh element $K_e$ is given as the image of the reference element $\hat{K}$ under an invertible, polynomial mapping $\T_e: \hat{K} \rightarrow K_e$ where $\T_e \in [\mathcal{P}_m(\hat{K})]^{\dim}$ for simplicial elements (triangles and tetrahedra) or $\T_e \in [\mathcal{Q}_m(\hat{K})]^{\dim}$ for tensor product elements (quadrilaterals and hexahedra). Here, $\mathcal{P}_m(\hat{K})$ is the space of polynomials of total degree at most $m$ in \emph{all} variables and $\mathcal{Q}_m(\hat{K})$ the space of polynomials of degree at most $m$ in \emph{each} variable. For example, in two dimensions,
	\begin{equation}
		\mathcal{P}_1(\hat{K}) = \{ 1, \xi_1, \xi_2 \} 
	\end{equation}
while 
	\begin{equation}
		\mathcal{Q}_1(\hat{K}) = \{ 1, \xi_1, \xi_2, \xi_1\xi_2 \} \,. 
	\end{equation}
We do not consider the use of $\mathcal{P}_m(\hat{K})$ on tensor-product elements for either the mesh or the solution. 

The reference element is the unit $\dim$-simplex for simplicial elements (i.e.~a triangle with coordinates (0,0), (1,0), and (0,1)) or the unit $\dim$-cube $\hat{K} = [0,1]^{\dim}$ for tensor product elements. Figure \ref{fig:trap_trans} depicts a mesh transformation for a non-affine, linear, quadrilateral element. In the remainder of this document, we assume the use of tensor product elements however the derivations apply analogously to simplicial elements. 

Let $\vec{\xi} \in \hat{K}$ denote the reference coordinate. The Jacobian matrix of the mapping is 
	\begin{equation} \label{eq:jacobian}
		\mat{F}_e = \pderiv{\T_e}{\vec{\xi}} \in \R^{\dim\times\dim}\,. 
	\end{equation}
Furthermore, we define $J_e = |\mat{F}_e|$ as the determinant of the Jacobian matrix. 
As an example, the transformation, Jacobian matrix, and determinant for the transformation depicted in Fig.~\ref{fig:trap_trans} are
	\begin{equation}
		\T = \begin{bmatrix} 
			h\xi_1 + \alpha \xi_2(2\xi_1 - 1) \\ h\xi_2
		\end{bmatrix} \,, \quad 
		\mat{F} = \begin{bmatrix} 
			2\alpha \xi_2 + h & \alpha(2\xi_1 - 1) \\ 0 & h
		\end{bmatrix} \,, \quad J = 2 \alpha \xi_2 h + h^2 \,. 
	\end{equation}
The mesh transformations are used to perform integration in reference space using: 
	\begin{equation} \label{eq:ref_int}
		\int \paren{\cdot} \ud \x = \sum_{K_e\in\meshT}\int_{K_e} \paren{\cdot} \ud \x = \sum_{K_e \in \meshT} \int_{\hat{K}} \paren{\cdot}\,J_e\!\ud \vec{\xi} \,. 
	\end{equation}
For integrands involving gradients, the chain rule implies that 
	\begin{equation} \label{eq:grad_trans}
		\nabla_{\x} = \mat{F}^{-T}\nabla_{\vec{\xi}} \,. 
	\end{equation}
Integration over surfaces is performed over the $\dim-1$ dimensional reference element using the transformed element of surface area.
In this document, integration over the domain $\D$ is implicitly performed using numerical quadrature and the relations in Eqs.~\ref{eq:ref_int} and \ref{eq:grad_trans}. Finally, the characteristic mesh length, $h$, is computed with 
	\begin{equation}
		h_e = \paren{\int_{\hat{K}} J_e \ud \vec{\xi}}^{1/\dim} \,,  
	\end{equation}
with $h = \max_{K_e\in\meshT} h_e$. 

\subsection{Finite Element Spaces}
On each element, we will seek solutions to the transport and VEF drift-diffusion equations in the space of polynomials mapped from the reference element $\hat{K}$ defined by 
	\begin{equation}
		\mathbb{Q}_p(K_e) = \{ u = \hat{u} \circ \T_e^{-1} : \hat{u} \in \mathcal{Q}_p(\hat{K}) \} \,,
	\end{equation}
where $\hat{u}$ indicates a function defined on the reference element. The delineation between $\mathcal{Q}$ and $\mathbb{Q}$ is required when non-affine\footnote{Examples of non-affine transformations include mapping the reference square to a trapezoid or any high-order, curved element.} mesh transformations are used. In such a case, $u = \hat{u}\circ \T_e^{-1} \notin \mathcal{Q}_p(K_e)$ even if $\hat{u}\in\mathcal{Q}_p(\hat{K})$. That is, the solution can be non-polynomial due to the composition with the inverse of the element transformation. For example, the inverse of the transformation in Fig.~\ref{fig:trap_trans} is 
	\begin{equation}
		\T^{-1} = \begin{bmatrix} 
			\frac{hx_1 + \alpha x_2}{h^2 + 2\alpha x_2} \\ x_2/h
		\end{bmatrix}
	\end{equation}
which is non-polynomial in the first coordinate. 

The degree-$p$ DG finite element space is: 
	\begin{equation}
		Y_p = \{ u \in L^2(\D) : u|_{K_e} \in \mathbb{Q}_p(K_e) \,, \quad \forall K_e \in \meshT \}
	\end{equation}
so that each function $u \in Y_p$ is a piecewise polynomial mapped from the reference element with no continuity enforced between elements. Its vector-valued analog is 
	\begin{equation}
		W_p = \{ \vec{v} \in [L^2(\D)]^{\dim} : \vec{v} \in [\mathbb{Q}_p(K_e)]^{\dim}\,, \quad \forall K_e \in \meshT \}\,,
	\end{equation}
which simply uses the scalar DG space for each component of the vector. 
We will also need the discrete $H^1(\D)$, or continuous finite element space, defined as: 
	\begin{equation}
		V_p = \{ u \in H^1(\D) : u|_{K_e} \in \mathbb{Q}_p(K_e) \,, \quad \forall K_e \in \meshT \} \,. 
	\end{equation}
Here, $u\in V_p$ is a piecewise continuous mapped polynomial. 

\begin{response}[orange][subparam]
For each of the above spaces, we allow the polynomial degree $p$ to be chosen independently of $m$, the polynomial degree used to describe the mesh, and thus support sub-, iso-, and super-parametric approximations. For our target application of Lagrangian hydrodynamics, the polynomial degree of the mesh is defined by the finite element representation of the fluid velocity. Typically, thermodynamic variables are approximated with polynomials one degree lower than the polynomials used for the fluid velocity. That is, for degree-$p$ transport, the mesh will be degree $m = p+1$, leading to a sub-parametric approximation for the radiation component of the multiphysics simulation.
\end{response}

A nodal basis for the element-local polynomial space is used. For a degree-$p$ element, let $\xi_i$ denote the $(p+1)$ Gauss-Lobatto or Gauss-Legendre points in the interval $[0,1]$. The $(p+1)^{\dim}$ points $\vec{\xi}_i$ on the unit cube $[0,1]^{\dim}$ are given by the $\dim$-fold Cartesian product of the one-dimensional points. Let $\ell_i$ denote the Lagrange interpolating polynomial satisfying $\ell_i(\vec{\xi}_j) = \delta_{ij}$ where $\delta_{ij}$ is the Kronecker delta. The set of functions $\{\ell_i\}$ form a basis for the space $\mathcal{Q}_p(\hat{K})$. The DG and $H^1(\D)$ finite element spaces are built element-by-element from this local basis. 

Note that the Gauss-Lobatto points include the interval end points while the Gauss-Legendre points do not. Thus, using Gauss-Lobatto points yields both points on the interior and the boundary of the element while using Gauss-Legendre leads to points on the interior of the element only. These are referred to as closed and open bases, respectively. In the case of DG, no continuity between elements is enforced so it is acceptable to use either an open or closed basis. Both Gauss-Lobatto and Gauss-Legendre have the required properties to be accurate in the limit $p\rightarrow\infty$ so the choice of Gauss-Lobatto versus Gauss-Legendre is typically dictated by other aspects of the overall algorithm such as preconditioners. The basis formed from the Gauss-Legendre points has the beneficial property of diagonal mass matrices on affine meshes, while the basis formed from Gauss-Lobatto points typically leads to sparser global systems since closed bases couple fewer unknowns on interior faces. A closed basis is required for $H^1(\D)$ finite element spaces to enable the strong enforcement of continuity between elements. 

\subsection{Mathematical Notation}
It is helpful to define the ``broken'' gradient, denoted $\nablah$, obtained by applying the gradient locally on each element. That is, 
	\begin{equation} \label{eq:nablah}
		\paren{\nablah u}|_{K_e} = \nabla\!\paren{u|_{K_e}} \,, \quad \forall K_e \in \meshT \,. 
	\end{equation}
This distinction is important since for $u\in Y_p$, $\nabla u$ is not well-defined since $u$ may be discontinuous across element interfaces.
However, $\nablah u$ is well-defined since $u$ is locally differentiable on each element. 

\begin{figure}
	\centering
	\includegraphics[width=.5\textwidth]{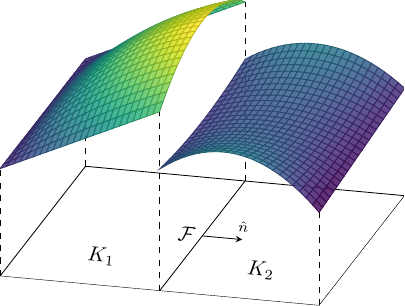}
	\caption{A depiction of a discontinuous, piecewise quadratic solution across two quadrilateral elements. The normal vector, $\hat{n}$, is defined as pointing from $K_1$ to $K_2$ along the face between $K_1$ and $K_2$.}
	\label{fig:jump_avg}
\end{figure}

We will use the following notation to describe the jump and average of a discontinuous function along an interior mesh face. Let $\Gamma$ be the set of all unique faces in the mesh and $\Gamma_0 = \Gamma \setminus \partial\D$ the set of unique interior faces. Additionally, define $\Gamma_b = \Gamma \cap \partial\D$ as the set of faces on the boundary so that $\Gamma = \Gamma_0 \cup \Gamma_b$. We define $\hat{n}_K$ as the outward unit normal to element $K$. On an interior face $\mathcal{F} \in \Gamma_0$ between elements $K_1$ and $K_2$, we use the convention that $\hat{n}$ is the the unit vector perpendicular to the shared face $K_1 \cap K_2$ pointing from $K_1$ to $K_2$ (see Fig.~\ref{fig:jump_avg}). On such an interior face, the jump, $\jump{\cdot}$, and average, $\avg{\cdot}$, are defined as
	\begin{equation} \label{eq:jumpavg_def}
		\jump{u} = u_1 - u_2 \,, \quad \avg{u} = \frac{1}{2}\paren{u_1 + u_2} \,, \quad \text{on} \ \mathcal{F} \in \Gamma_0 \,, 
	\end{equation}
where $u_i = u|_{\partial K_i}$ with analogous definitions for vectors.

Note that in contrast to the notation of \cite{Arnold2002}, our jump operator does not change the rank of its argument: the jump of a scalar is a scalar and the jump of a vector is a vector. 
Consequently, our notation is not invariant under element renumbering, since flipping the ordering of the elements negates the value of the jump.
However, the bilinear and linear forms presented in this paper always pair the jump with another normal-dependent term.
The negation of the jump induced by swapping the element ordering is then balanced by flipping the orientation of the normal vector, and so the discretizations under consideration are in fact invariant with respect to the element ordering.

On the boundary of the domain, we set the jump and average to 
	\begin{equation}
		\jump{u} = u\,, \quad \avg{u} = u \,, \quad \text{on} \ \mathcal{F} \in \Gamma_b \,,
	\end{equation}
and likewise for vector-valued functions on the boundary.
A straightforward computation shows that 
	\begin{equation} \label{ref:jump_avg_id}
		\sum_{K\in\meshT} \int_{\partial K} u\,\vec{v}\cdot\hat{n}_K \ud s = \int_{\Gamma} \jump{u\,\vec{v}\cdot\hat{n}}\ud s = \int_\Gamma \jump{u} \avg{\vec{v}\cdot\hat{n}} \ud s + \int_{\Gamma_0} \avg{u}\jump{\vec{v}\cdot\hat{n}} \ud s \,. 
	\end{equation}
We refer to this as the ``jumps and averages identity''. 
The restriction of the integration to interior faces for the second term in the last equality is consistent with the notation of \cite{Arnold2002} and is used so that only one term contributes on the boundary of the domain.

Finally, we refer to a function as ``single-valued'' on an interior face if its values obtained from approaching from each side of the face are identical so that 
	\begin{equation}
		\jump{u} = 0 \,, \quad \avg{u} = u \,. 
	\end{equation}
Note in particular that the jump and average operators are single-valued. 

\section{Transport Discretizations}
In this work, we assume the transport equation is discretized with the Discrete Ordinates (\Sn) angular model and an arbitrary-order Discontinuous Galerkin (DG) spatial discretization compatible with curved meshes (e.g. \cite{woods_thesis,graph_sweeps}). In \Sn, the transport equation is collocated at discrete angles, $\Omegahat_d$, and integration is numerically approximated using a suitable angular quadrature rule $\{\Omegahat_d, w_d\}_{d=1}^{N_\Omega}$ on the unit sphere. The VEF data are then 
	\begin{subequations} \label{eq:vef_disc}
	\begin{equation}
		\E(\x) = \frac{\sum_{d=1}^{N_\Omega}w_d\, \Omegahat_d\otimes\Omegahat_d\, \psi_d(\x)}{\sum_{d=1}^{N_\Omega} w_d \psi_d(\x)} \,,
	\end{equation} 
	\begin{equation}
		E_b(\x) = \frac{\sum_{d=1}^{N_\Omega} w_d\,|\Omegahat_d\cdot\hat{n}| \,\psi_d(\x)}{\sum_{d=1}^{N_\Omega} w_d \psi_d(\x)} \,, 
	\end{equation}
	\end{subequations}
where $\psi_d(\x) = \psi(\x,\Omegahat_d)$ is the discrete angular flux in direction $\Omegahat_d$. With degree-$p$ DG in space, the angular flux in each discrete direction $\Omegahat_d$ is a member of $Y_p$. Through the standard finite element interpolation procedure, the Eddington tensor and boundary factor in Eq.~\ref{eq:vef_disc} can be evaluated at any location in the mesh. Note that it is important to interpolate the numerator and denominator of the VEF data \emph{independently}. That is, the boundary factor and each component of the Eddington tensor are represented as degree-$p$ improper rational polynomials on each element. \resp[blue][ratint]{Improper rational polynomials cannot be integrated exactly with numerical quadrature. Thus, bilinear and linear forms involving VEF data will possess integration error. In practice, we have seen that the optimal order of convergence is maintained despite this inexact numerical integration. This observation is corroborated by \citet[\S4.1]{ciarlet2002finite} which presents an analysis of the stability and accuracy of the general finite element method when inexact numerical integration is used.}

Defining
	\begin{equation}
		\P(\x) = \sum_d w_d \, \Omegahat_d\otimes \Omegahat_d \, \psi_d(\x) 
	\end{equation}
as the discrete second moment of the angular flux and using the quotient rule, the local divergence of the Eddington tensor 
	\begin{equation} \label{eq:Ediv}
		\nablah\cdot\E = \frac{\paren{\nablah\cdot\P}\!\phi - \P\cdot\nablah \phi}{\phi^2} 
	\end{equation}
is well-defined assuming $\phi>0$. Here, the divergence of a second-order tensor is the vector formed by taking the divergence of each of the columns of the tensor. 

We restrict our attention to problems where $\psi \geq \delta > 0$ inside the domain, for some $\delta$.
This assumption is reasonable for our applications but may be violated in shielding or deep penetration problems.
Application of a positivity-preserving negative flux fixup then ensures that $\phi$ is bounded away from zero, so that $\E$, $E_b$, and $\nablah\cdot\E$ are all bounded.
Thus, through \Sn angular quadrature and finite element interpolation, the Eddington tensor, boundary factor, and the local divergence of the Eddington tensor can be evaluated at any point in any element of the mesh.
This completes the definition of the connection between the discrete transport equation and the VEF drift-diffusion equation.
Note that since the angular flux is generally discontinuous across interior mesh interfaces, the Eddington tensor and its divergence also will be.
Thus, we will carefully design the discretization of the VEF drift-diffusion equation to accommodate discontinuous data. 

The VEF scalar flux connects with the transport equation in the scattering source. To support generality, we assume that the finite element space for the VEF scalar flux and the finite element space for the angular flux are different. The scattering source is then constructed using a mixed-space mass matrix that has test functions in the space for the angular flux and trial functions in the space for the VEF scalar flux. 

\section{Derivation of DG VEF} \label{sec:dg-vef}

In this section, we adapt the derivation of the unified framework for DG methods designed for the Poisson equation in \cite{Arnold2002} to the VEF equations. This enables the use of any of the DG methods described there. 
\citet{Arnold2002} derive a family of DG methods for:
	\begin{subequations}
	\begin{equation}
		\vec{q} = \nabla u \,, 
	\end{equation}
	\begin{equation}
		-\nabla\cdot\vec{q} = f\,,
	\end{equation}
	\end{subequations}
with Dirichlet boundary conditions. The present goal is to adapt their derivation to the VEF equations: 
	\begin{subequations}
	\begin{equation}
		\nabla\cdot\paren{\E\varphi} + \sigma_t \vec{J} = \vec{Q}_1 \,,
	\end{equation}
	\begin{equation}
		\nabla\cdot\vec{J} + \sigma_a \varphi = Q_0 \,, 
	\end{equation}
	\end{subequations}
with the Robin style boundary conditions given in Eq.~\ref{eq:mlbc}. We will see significant differences in the final equation since the Eddington tensor is inside the divergence.
Additionally, the presence of a right-hand side in the first moment equation as well as non-unit coefficients introduce further complications.
We will then derive analogs of the interior penalty (IP), second method of Bassi and Rebay (BR2), and minimal dissipation local Discontinuous Galerkin (MDLDG) variants.
Finally, we will show how to extract a continuous finite element method from this framework. 

\subsection{Adaption of the Unified Framework to VEF}
We seek the VEF scalar flux in the degree-$p$ DG finite element space $Y_p$ and the current in the degree-$p$, vector-valued DG finite element space $W_p$. The weak form is then: find $(\varphi,\vec{J}) \in Y_p \times W_p$ such that for each $K \in \meshT$: 
	\begin{subequations}
	\begin{equation}
		\int_{\partial K} \vec{v}\cdot\widehat{\E\varphi}\hat{n}_K \ud s - \int_K \nabla\vec{v}|_K : \E \varphi \ud \x + \int_K \sigma_t \, \vec{v}\cdot\vec{J} \ud \x = \int_K \vec{v}\cdot\vec{Q}_1 \ud \x \,, \quad \forall \vec{v} \in [\mathbb{Q}_p(K)]^{\dim} \,, 
	\end{equation}
	\begin{equation}
		\int_{\partial K} u\,\widehat{\vec{J}}\cdot\hat{n}_K \ud s - \int_K \nabla u|_K \cdot\vec{J} \ud \x + \int_K \sigma_a\, u \varphi \ud \x = \int_K u\, Q_0 \ud \x \,, \quad \forall u \in \mathbb{Q}_p(K) \,, 
	\end{equation}
	\end{subequations}
where the \emph{numerical fluxes} $\widehat{\E\varphi}$ and $\widehat{\vec{J}}$ are approximations of $\E\varphi$ and $\vec{J}$ on the boundaries of the elements in the mesh. In the above, integration by parts was applied on each element so that only local differentiation on each element is required for functions in $Y_p$ and $W_p$. We have grouped the product $\E\varphi$ as the numerical flux to mimic the integration by parts of the product of a tensor and vector. Here, the gradient of a vector is 
	\begin{equation}
		\paren{\nabla\vec{v}}_{ij} = \paren{\pderiv{\vec{v}_i}{\x_j}} \in \R^{\dim\times \dim}
	\end{equation}
and 
	\begin{equation}
		\mat{A} : \mat{B} = \sum_{i=1}^{\dim} \sum_{j=1}^{\dim} \mat{A}_{ij} \mat{B}_{ij} \,, \quad \mat{A}, \mat{B} \in \R^{\dim\times \dim} 
	\end{equation}
is the scalar contraction of two tensors. Note that if $\E = \frac{1}{3}\I$ then 
	\begin{equation}
		\nabla\vec{v} : \E = \frac{1}{3}\nabla\cdot\vec{v} 
	\end{equation}
and the symmetric weak form for radiation diffusion can be recovered. 

Summing the zeroth moment over all elements: 
	\begin{equation}
		\int_\Gamma \jump{u} \avg{\widehat{\vec{J}}\cdot\hat{n}} \ud s + \int_{\Gamma_0} \avg{u}\jump{\widehat{\vec{J}}\cdot\hat{n}} \ud s - \int \nablah u \cdot \vec{J} \ud \x + \int \sigma_a \, u\varphi \ud \x = \int u\, Q_0 \ud \x \,,
	\end{equation}
where the jumps and averages identity (Eq.~\ref{ref:jump_avg_id}) was used along with the definition of the broken gradient from Eq.~\ref{eq:nablah}. We will now use the discrete first moment to determine a functional form for $\vec{J}$. Integrating by parts locally over element $K$, we have that 
	\begin{equation} \label{eq:identity}
		\int_K \nabla\vec{v}|_K : \E \varphi \ud \x = \int_{\partial K} \vec{v}\cdot\E\varphi\hat{n}_K \ud s - \int_K \vec{v}\cdot\nabla\cdot\paren{\E\varphi}\!|_K \ud \x \,. 
	\end{equation}
Here, a numerical flux is not required since the integration by parts is performed on the gradient restricted to to each element $K$. 
The first moment's weak form on each element becomes: 
	\begin{equation}
		\int_{\partial K} \vec{v}\cdot\paren{\widehat{\E\varphi}\hat{n}_K - \E\varphi\hat{n}_K} \ud s + \int_K \vec{v}\cdot\nabla\cdot\paren{\E\varphi}\!|_K \ud \x + \int_K \sigma_t \, \vec{v}\cdot\vec{J} \ud \x = \int_K \vec{v} \cdot\vec{Q}_1 \ud \x \,, \quad \forall \vec{v} \in [\mathbb{Q}_p(K)]^{\dim} \,. 
	\end{equation}
Summing over all elements and using the jumps and averages identity, the weak form for the first moment is:  
	\begin{multline} \label{eq:weak_strong}
		\int_\Gamma \avg{\vec{v}}\cdot\jump{\widehat{\E\varphi}\hat{n} - \E\varphi\hat{n}} \ud s + \int_{\Gamma_0} \jump{\vec{v}} \cdot \avg{\widehat{\E\varphi}\hat{n} - \E\varphi\hat{n}} \ud s \\+ \int \vec{v}\cdot\nablah \cdot\paren{\E\varphi} \ud \x + \int \sigma_t \, \vec{v}\cdot\vec{J} \ud \x = \int \vec{v}\cdot\vec{Q}_1 \ud \x \,, \quad \forall \vec{v} \in W_p \,, 
	\end{multline}
where $\nablah\cdot\paren{\E\varphi}$ is evaluated as $\nablah\cdot\paren{\E\varphi} = \E\nablah \varphi + \paren{\nablah\cdot\E}\!\varphi$, and the term $\nablah\cdot\E$ is computed using Eq.~\ref{eq:Ediv}. 

We now wish to write all terms as volumetric integrals so that a functional form for the current can be found. To that end, define \emph{lifting operators} $\vec{r}(\vec{\tau}) \in W_p$ and $\vec{\ell}(\vec{\chi}) \in W_p$ such that 
	\begin{subequations}
	\begin{equation} \label{eq:lift-r}
		\int \sigma_t\, \vec{v}\cdot\vec{r}(\vec{\tau}) \ud \x = - \int_{\Gamma} \avg{\vec{v}} \cdot \vec{\tau} \ud s \,, \quad \forall \vec{v} \in W_p \,, 
	\end{equation}
	\begin{equation} \label{eq:lift-l}
		\int \sigma_t \, \vec{v}\cdot\vec{\ell}(\vec{\chi}) \ud \x = -\int_{\Gamma_0} \jump{\vec{v}} \cdot \vec{\chi} \ud s \,, \quad \forall \vec{v} \in W_p \,, 
	\end{equation}
	\end{subequations}
where $\vec{\tau}$ and $\vec{\chi}$ are vector functions that are singled-valued on $\Gamma_0$. 
Note that the lifting operators are finite element grid functions just as the current is and that the left hand sides are simply the $W_p$ total interaction mass matrix.
Since $W_p$ is piecewise discontinuous, the $W_p$ mass matrix is block-diagonal by element and thus the systems of equations corresponding to Eqs.~\ref{eq:lift-r} and \ref{eq:lift-l} are amenable to efficient direct factorization (see \ref{sec:lifting}). 

Setting $\vec{\tau} = \jump{\widehat{\E\varphi}\hat{n} - \E\varphi\hat{n}}$ and $\vec{\chi} = \avg{\widehat{\E\varphi}\hat{n} - \E\varphi\hat{n}}$ and using the definitions of the lifting operators, Eq.~\ref{eq:weak_strong} can be written entirely in terms of volumetric integrals as: 
	\begin{equation}
		\int \sigma_t\, \vec{v}\cdot\vec{J} \ud \x = \int \sigma_t\, \vec{v}\cdot\bracket{\frac{1}{\sigma_t} \paren{\vec{Q}_1 - \nablah\cdot\paren{\E\varphi}} + \vec{r}\!\paren{\jump{\widehat{\E\varphi}\hat{n} - \E\varphi\hat{n}}} + \vec{\ell}\!\paren{\avg{\widehat{\E\varphi}\hat{n} - \E\varphi\hat{n}}}} \ud \x 
	\end{equation}
for all $\vec{v}\in W_p$. Subtracting the right hand side and setting the integrand to zero implies that 
	\begin{equation} \label{eq:current_form}
		\vec{J} = \frac{1}{\sigma_t}\paren{\vec{Q}_1 - \nablah\cdot\paren{\E\varphi}} + \vec{r}\!\paren{\jump{\widehat{\E\varphi}\hat{n} - \E\varphi\hat{n}}} + \vec{\ell}\!\paren{\avg{\widehat{\E\varphi}\hat{n} - \E\varphi\hat{n}}}\,. 
	\end{equation}
Observe that the above represents the element-local strong form of the current, $\frac{1}{\sigma_t}\paren{\vec{Q}_1 - \nablah\cdot\paren{\E\varphi}}$ found by analytically eliminating the current, with additional terms that capture the effect of the numerical fluxes. In other words, we have derived the \emph{discrete} elimination of the current. 

Using this discrete form for the current and the definitions of the lifting operators to convert from volumetric integrals back to surface integrals, the zeroth moment becomes: 
	\begin{multline} \label{eq:familynobc}
		\int_\Gamma \jump{u} \avg{\widehat{\vec{J}}\cdot\hat{n}} \ud s + \int_{\Gamma_0} \avg{u}\!\jump{\widehat{\vec{J}}\cdot\hat{n}} \ud s + \int_\Gamma \avg{\frac{\nablah u}{\sigma_t}} \cdot \jump{\widehat{\E\varphi}\hat{n} - \E\varphi\hat{n}} \ud s \\ 
		+ \int_{\Gamma_0} \jump{\frac{\nablah u}{\sigma_t}} \cdot \avg{\widehat{\E\varphi}\hat{n} - \E\varphi\hat{n}} \ud s + \int \nablah u \cdot \frac{1}{\sigma_t}\nablah\cdot\paren{\E\varphi} \ud \x + \int \sigma_a\, u \varphi \ud \x \\ 
		= \int u\, Q_0 \ud \x + \int \nablah u \cdot \frac{\vec{Q}_1}{\sigma_t} \ud \x \,, \quad \forall u \in Y_p \,. 
	\end{multline}
On boundary faces, we apply the Miften-Larsen boundary conditions by setting 
	\begin{equation} \label{eq:ml_bdr_flux}
		\widehat{\vec{J}}\cdot\hat{n} = 2g + E_b \varphi \,, \quad \widehat{\E\varphi}\hat{n} = \E\varphi\hat{n} \,, \quad \text{on} \ \mathcal{F} \in \Gamma_b \,. 
	\end{equation}
All the methods we consider use so-called conservative numerical fluxes such that 
	\begin{subequations}
	\begin{equation}
		\jump{\widehat{\vec{J}}\cdot\hat{n}} = 0 \,, \quad \avg{\widehat{\vec{J}}\cdot\hat{n}} = \widehat{\vec{J}}\cdot\hat{n} \,, \quad \text{on} \ \mathcal{F} \in \Gamma_0 \,, 
	\end{equation}
	\begin{equation}
		\jump{\widehat{\E\varphi}\hat{n}} = 0 \,, \quad \avg{\widehat{\E\varphi}\hat{n}} = \widehat{\E\varphi}\hat{n} \,, \quad \text{on} \ \mathcal{F} \in \Gamma_0 \,. 
	\end{equation}
	\end{subequations}
Using the boundary conditions and the assumption of conservative numerical fluxes, Eq.~\ref{eq:familynobc} becomes: 
	\begin{multline} \label{eq:family}
		\int_{\Gamma_b} E_b\, u \varphi \ud s + \int_{\Gamma_0} \jump{u} \widehat{\vec{J}}\cdot\hat{n} \ud s - \int_{\Gamma_0} \avg{\frac{\nablah u}{\sigma_t}} \cdot \jump{\E\varphi\hat{n}} \ud s \\+ \int_{\Gamma_0} \jump{\frac{\nablah u}{\sigma_t}} \cdot \avg{\widehat{\E\varphi}\hat{n} - \E\varphi\hat{n}} \ud s + \int \nablah u \cdot \frac{1}{\sigma_t}\nablah\cdot\paren{\E\varphi} \ud \x + \int \sigma_a\, u \varphi \ud \x \\ 
		= \int u\, Q_0 \ud \x + \int \nablah u \cdot \frac{\vec{Q}_1}{\sigma_t} \ud \x - 2\int_{\Gamma_b} u\, g \ud s \,, \quad \forall u \in Y_p \,. 
	\end{multline}
Equation \ref{eq:family} defines a \emph{family} of DG methods. That is, through the specification of the numerical fluxes on interior faces, analogs of all the methods listed in \cite{Arnold2002} can be derived.

\subsection{Specification of Numerical Fluxes}
All the methods we consider use numerical fluxes of the form 
	\begin{subequations}
	\begin{equation}
		\widehat{\vec{J}}\cdot\hat{n} = \avg{\frac{1}{\sigma_t}\paren{\vec{Q}_1 - \nablah\cdot\paren{\E\varphi}}\cdot\hat{n}} + \alpha(\varphi) \,, \quad \text{on} \ \Gamma_0 \,, 
	\end{equation}
	\begin{equation}
		\widehat{\E\varphi}\hat{n} = \avg{\E\varphi\hat{n}} + \vec{\theta}(\varphi) \,,\quad \text{on} \ \Gamma_0 \,, 
	\end{equation}
	\end{subequations}
where $\alpha(\varphi)$ and $\vec{\theta}(\varphi)$ are single-valued functions whose purpose are to ensure a stable discretization.
The IP, BR2, and LDG methods differ only in the choice of $\alpha(\varphi)$ and $\vec{\theta}(\varphi)$. With these numerical fluxes, Eq.~\ref{eq:family} becomes: 
	\begin{multline} \label{eq:family_alpha}
		\int_{\Gamma_b} E_b\, u \varphi \ud s + \int_{\Gamma_0} \jump{u} \alpha(\varphi) \ud s - \int_{\Gamma_0} \jump{u} \avg{\frac{1}{\sigma_t}\nablah\cdot\paren{\E\varphi}\cdot\hat{n}} \ud s - \int_{\Gamma_0} \avg{\frac{\nablah u}{\sigma_t}} \cdot \jump{\E\varphi\hat{n}} \ud s \\ + \int_{\Gamma_0} \jump{\frac{\nablah u}{\sigma_t}} \cdot \vec{\theta}(\varphi) \ud s
		+ \int \nablah u \cdot \frac{1}{\sigma_t}\nablah\cdot\paren{\E\varphi} \ud \x + \int \sigma_a\, u \varphi \ud \x \\ 
		= \int u\, Q_0 \ud \x + \int \nablah u \cdot \frac{\vec{Q}_1}{\sigma_t} \ud \x - \int_{\Gamma_0} \jump{u} \avg{\frac{\vec{Q}_1\cdot\hat{n}}{\sigma_t}} \ud s - 2\int_{\Gamma_b} u\, g \ud s \,, \quad \forall u \in Y_p \,. 
	\end{multline}
Recall that this form has already applied boundary conditions according to Eq.~\ref{eq:ml_bdr_flux}. In other words, the above corresponds to a DG scheme with the following numerical fluxes:
		\begin{subequations}
		\begin{equation}
			\widehat{\vec{J}}\cdot\hat{n} = \begin{cases}
				\avg{\frac{1}{\sigma_t}\paren{\vec{Q}_1 - \nablah\cdot\paren{\E\varphi}}\cdot\hat{n}} + \alpha(\varphi) \,, & \text{on} \ \Gamma_0 \\ 
				2g + E_b \varphi \,, & \text{on} \ \Gamma_b 
			\end{cases} \,, 
		\end{equation}
		\begin{equation}
			\widehat{\E\varphi}\hat{n} = \begin{cases}
				\avg{\E\varphi\hat{n}} + \vec{\theta}(\varphi) \,, & \text{on} \ \Gamma_0 \\ 
				\E\varphi\hat{n} \,, & \text{on} \ \Gamma_b 
			\end{cases} \,.
		\end{equation}
		\end{subequations}

\subsubsection{Interior Penalty}
An interior penalty (IP)-like method uses 
	\begin{equation}
		\alpha(\varphi) = \kappa \jump{\varphi} \,, \quad \vec{\theta}(\varphi) = 0 \,, 
	\end{equation}
where $\kappa$ is the penalty parameter. The full IP weak form is then: find $\varphi \in Y_p$ such that 
	\begin{multline} \label{eq:ip}
		\int_{\Gamma_b} E_b\, u \varphi \ud s + \int_{\Gamma_0} \kappa \jump{u} \jump{\varphi} \ud s - \int_{\Gamma_0} \jump{u} \avg{\frac{1}{\sigma_t}\nablah\cdot\paren{\E\varphi}\cdot\hat{n}} \ud s - \int_{\Gamma_0} \avg{\frac{\nablah u}{\sigma_t}} \cdot \jump{\E\varphi\hat{n}} \ud s \\
		+ \int \nablah u \cdot \frac{1}{\sigma_t}\nablah\cdot\paren{\E\varphi} \ud \x + \int \sigma_a\, u \varphi \ud \x \\ 
		= \int u\, Q_0 \ud \x + \int \nablah u \cdot \frac{\vec{Q}_1}{\sigma_t} \ud \x - \int_{\Gamma_0} \jump{u} \avg{\frac{\vec{Q}_1\cdot\hat{n}}{\sigma_t}} \ud s - 2\int_{\Gamma_b} u\, g \ud s \,, \quad \forall u \in Y_p \,. 
	\end{multline}
\begin{response}[orange][ipexp]
IP methods require that $\kappa \propto \sigma_t^{-1} p^2/h$ in order to guarantee stability as the mesh is refined. The constant of proportionality is a user-defined parameter that is often problem dependent. For example, we will see that severely distorted meshes require the penalty parameter to be increased in order for the IP VEF method to be stable. 
We note that the penalty bilinear form, given by
	\begin{equation} \label{eq:penalty_bilin}
		\int_{\Gamma_0} \kappa \jump{u} \jump{\varphi} \ud s \,,
	\end{equation}
is symmetric positive definite and has a nullspace corresponding to functions that are continuous on the interior of the domain. 
A large enough penalty parameter causes the penalty bilinear form to dominate the negative definite bilinear forms in the discretization making the overall system positive definite. 
However, a large penalty parameter also increases the relative dominance of the penalty bilinear form's nullspace. This has the effect of 1) regularizing the solution towards continuous functions such that the limit $\kappa \rightarrow \infty$ would yield a continuous solution and 2) increasing the linear system's condition number causing the effectiveness of standard preconditioners (e.g.~AMG) to degrade as the mesh is refined. 
This sub-optimal performance was the motivation for the development of the uniform subspace correction preconditioner \cite{Pazner2021} which achieves iterative convergence independent of the mesh size, polynomial order, and penalty parameter. In \S\ref{sec:subspace}, the analysis of this preconditioner is extended to the non-symmetric case of the VEF equations.
\end{response}

\subsubsection{BR2}
The second method of Bassi and Rebay (BR2) uses an alternative penalty term. Let $\vec{\rho}_f(\omega)\in W_p$ be a face-local lifting operator defined by 
	\begin{equation} \label{eq:rhof}
		\int \vec{v}\cdot\vec{\rho}_f(\omega) \ud \x = -\int_f \avg{\vec{v} \cdot\hat{n}} \omega \ud s \,, \quad \forall \vec{v} \in W_p\,, \quad \text{on} \ f\in \Gamma_0 \,. 
	\end{equation}
Here, $\omega$ is a scalar function that is single-valued on the interior face $f$. Note that the integration on the left hand side is over the entire domain while the right hand side is localized to a single interior face. This means the right hand side, and thus $\vec{\rho}_f(\omega)$, will be non-zero only for DOFs in elements that share the face $f$. 

A BR2-like discretization sets 
	\begin{equation}
		\alpha(\varphi) = -\eta\avg{\vec{\rho}_f(\jump{\varphi})\cdot\hat{n}} \,, \quad \text{on} \ f \in \Gamma_0 \,, \quad \vec{\theta}(\varphi) = 0 \,,
	\end{equation}
so that the relevant term is 
	\begin{equation} \label{eq:br2_stab}
	\begin{aligned}
		\int_{\Gamma_0} \jump{u} \alpha(\varphi) \ud s &= -\sum_{f\in\Gamma_0}\int_{f} \eta\,\jump{u}\!\avg{\vec{\rho}_f(\jump{u})\cdot\hat{n}} \ud s \\
		&= \sum_{f\in\Gamma_0}\int \eta\,\vec{\rho}_f(\jump{u})\cdot\vec{\rho}_f(\jump{\varphi}) \ud \x \,.
	\end{aligned}
	\end{equation}

The BR2 DG VEF discretization is then: find $\varphi \in Y_p$ such that 
	\begin{multline} \label{eq:br2}
		\int_{\Gamma_b} E_b\, u \varphi \ud s - \int_{\Gamma_0} \jump{u} \avg{\frac{1}{\sigma_t}\nablah\cdot\paren{\E\varphi}\cdot\hat{n}} \ud s - \int_{\Gamma_0} \avg{\frac{\nablah u}{\sigma_t}} \cdot \jump{\E\varphi\hat{n}} \ud s \\
		+ \sum_{f\in\Gamma_0} \int \eta\, \vec{\rho}_f(\jump{u}) \cdot \vec{\rho}_f(\jump{\varphi}) \ud \x + \int \nablah u \cdot \frac{1}{\sigma_t}\nablah\cdot\paren{\E\varphi} \ud \x + \int \sigma_a\, u \varphi \ud \x \\ 
		= \int u\, Q_0 \ud \x + \int \nablah u \cdot \frac{\vec{Q}_1}{\sigma_t} \ud \x - \int_{\Gamma_0} \jump{u} \avg{\frac{\vec{Q}_1\cdot\hat{n}}{\sigma_t}} \ud s - 2\int_{\Gamma_b} u\, g \ud s \,, \quad \forall u \in Y_p \,. 
	\end{multline}
\begin{response}[orange][brexp]
Observe that the BR2 and IP discretizations differ only in the stabilization term. The BR2 stabilization bilinear form, given by Eq.~\ref{eq:br2_stab}, is similar in function to the IP penalty bilinear form in Eq.~\ref{eq:penalty_bilin} in that it ensures the resulting algebraic system is positive definite, has the effect of regularizing toward continuous solutions, and increases the condition number of the algebraic system such that the specialized preconditioner discussed in \S6 is required. However, due to the use of the more expensive local lifting operators, the BR2 stabilization parameter, $\eta$, does not need to scale with the mesh size, polynomial order, or material parameters. Instead, $\eta$ can be prescribed by the geometric properties of the element types in the mesh alone. 
In particular, it has been shown for the model problem that stability is guaranteed when, on each $\mathcal{F} = K_1 \cap K_2 \in \Gamma_0$,
	\begin{equation}
		\eta \geq \max_{K \in [K_1,K_2]} n(K) \,, 
	\end{equation}
where $n(K)$ is the number of faces in element $K$ \cite[Prop.~1]{brstab}. 
For example, $\eta = 3$ and $\eta = 4$ lead to stable discretizations on meshes composed of triangular and quadrilateral elements, respectively. 
Thus, the BR2 discretization avoids the ambiguity associated with tuning the penalty parameter. This comes at the cost of a more expensive assembly procedure compared to IP. However, we stress that the BR2 stabilization term can still be assembled locally on each face in the mesh. Implementation details associated with the BR2 local lifting operators are provided in \ref{sec:lifting}. 
\end{response}

\subsubsection{Local Discontinuous Galerkin}
Finally, we consider the local Discontinuous Galerkin (LDG) method. In general, LDG uses the following numerical fluxes:
	\begin{subequations}
	\begin{equation}
		\widehat{\vec{J}}\cdot\hat{n} = \avg{\vec{J}\cdot\hat{n}} + \beta\jump{\vec{J}\cdot\hat{n}} + \kappa \jump{\varphi} \,,
	\end{equation}
	\begin{equation}
		\widehat{\E\varphi}\hat{n} = \avg{\E\varphi\hat{n}} - \beta\jump{\E\varphi\hat{n}} \,,
	\end{equation}
	\end{subequations}
where $\vec{J}$ is defined as the discrete elimination of the current derived in Eq.~\ref{eq:current_form}. The scalar parameter $\beta$ can be defined as 
	\begin{equation} \label{eq:ldgbeta}
		\beta = \begin{cases}
			1/2 \,, & \vec{w}\cdot\hat{n} > 0 \\ 
			-1/2\,, & \vec{w} \cdot\hat{n} < 0 
		\end{cases} \,,
	\end{equation}
where $\vec{w}$ is any constant, non-zero vector. This choice imposes an arbitrary upwinding on the current that is balanced by an opposing choice for the scalar flux.
With this choice of $\beta$, the LDG method is stable for any $\kappa \geq 0$; if $\kappa\equiv 0$, the method is referred to as the minimal dissipation LDG (MDLDG) method \cite{10.1007/s10915-007-9130-3}.
Using the numerical flux for the scalar flux, the discrete current simplifies to 
	\begin{equation}
		\vec{J} = \frac{1}{\sigma_t}\paren{\vec{Q}_1 - \nablah \cdot\paren{\E\varphi}} - \vec{r}_0\!\paren{\jump{\E\varphi\hat{n}}} - \vec{\ell}\!\paren{\beta\jump{\E\varphi\hat{n}}} \,,
	\end{equation}
where $\vec{r}_0(\vec{\tau}) \in W_p$ is another lifting operator defined by 
	\begin{equation}
	 	\int \sigma_t\, \vec{v}\cdot\vec{r}_0(\vec{\tau}) \ud \x = -\int_{\Gamma_0} \avg{\vec{v}}\cdot \vec{\tau} \ud s \,, \quad \forall \vec{v} \in W_p \,, 
	\end{equation} 
that differs from $\vec{r}(\vec{\tau})$ only in the region of integration on the right hand side. The LDG method is then equivalent to setting 
	\begin{subequations}
	\begin{multline}
		\alpha(\varphi) = -\avg{\vec{r}_0\!\paren{\jump{\E\varphi\hat{n}}}\cdot\hat{n} + \vec{\ell}\!\paren{\beta\jump{\E\varphi\hat{n}}}\cdot\hat{n}} \\+ \beta\jump{\frac{1}{\sigma_t}\paren{\vec{Q}_1 - \nablah\cdot\paren{\E\varphi}}\cdot\hat{n} - \vec{r}_0\!\paren{\jump{\E\varphi\hat{n}}}\cdot\hat{n} - \vec{\ell}\!\paren{\beta\jump{\E\varphi\hat{n}}}\cdot\hat{n}} + \kappa\jump{\varphi} \,,
	\end{multline}
	\begin{equation}
		\vec{\theta}(\varphi) = -\beta\jump{\E\varphi\hat{n}} \,.
	\end{equation}
	\end{subequations}
We then have that 
	\begin{multline}
		\int_{\Gamma_0} \jump{u} \alpha(\varphi) \ud s = \int_{\Gamma_0} \beta \jump{u}\jump{\frac{\vec{Q}_1\cdot\hat{n}}{\sigma_t}} \ud s - \int_{\Gamma_0}\beta\jump{u}\jump{\frac{1}{\sigma_t}\nablah\cdot\paren{\E\varphi}\cdot\hat{n}} \ud s \\
		+ \int \paren{\vec{\rho}_0\!\paren{\jump{u}} + \vec{\lambda}\!\paren{\beta\jump{u}}}\cdot\paren{\vec{r}_0\!\paren{\jump{\E\varphi\hat{n}}} + \vec{\ell}\!\paren{\beta\jump{\E\varphi\hat{n}}}} \ud \x + \int_{\Gamma_0} \kappa \jump{u}\jump{\varphi} \ud s 
	\end{multline}
where $\vec{\rho}_0(\omega), \vec{\lambda}(\upsilon) \in W_p$ such that 
	\begin{equation}
		\int \vec{v} \cdot \vec{\rho}_0(\omega) \ud \x = -\int_{\Gamma_0} \avg{\vec{v}\cdot\hat{n}} \omega \ud s \,, \quad \forall \vec{v} \in W_p \,,
	\end{equation}
	\begin{equation}
		\int \vec{v} \cdot \vec{\lambda}(\upsilon) \ud \x = -\int_{\Gamma_0} \jump{\vec{v}\cdot\hat{n}} \upsilon \ud s \,, \quad \forall \vec{v} \in W_p \,, 
	\end{equation}
are analogs of $\vec{r}_0(\vec{\tau})$ and $\vec{\ell}(\vec{\chi})$, respectively, that do not include the total interaction cross section in the left hand side mass matrices and have scalar arguments. The LDG VEF discretization is then: find $\varphi \in Y_p$ such that 
	\begin{multline} \label{eq:ldg}
		\int_{\Gamma_b} E_b\, u \varphi \ud s + \int_{\Gamma_0} \kappa\, \jump{u}\jump{\varphi} \ud s - \int_{\Gamma_0} \jump{u} \avg{\frac{1}{\sigma_t}\nablah\cdot\paren{\E\varphi}\cdot\hat{n}} \ud s - \int_{\Gamma_0} \avg{\frac{\nablah u}{\sigma_t}} \cdot \jump{\E\varphi\hat{n}} \ud s \\
		+ \int \paren{\vec{\rho}_0\!\paren{\jump{u}} + \vec{\lambda}\!\paren{\beta\jump{u}}}\cdot\paren{\vec{r}_0\!\paren{\jump{\E\varphi\hat{n}}} + \vec{\ell}\!\paren{\beta\jump{\E\varphi\hat{n}}}} \ud \x \\+ \int \nablah u \cdot \frac{1}{\sigma_t}\nablah\cdot\paren{\E\varphi} \ud \x + \int \sigma_a\, u \varphi \ud \x 
		\\= \int u\, Q_0 \ud \x + \int \nablah u \cdot \frac{\vec{Q}_1}{\sigma_t} \ud \x - \int_{\Gamma_0} \jump{u} \paren{\avg{\frac{\vec{Q}_1\cdot\hat{n}}{\sigma_t}} + \beta\jump{\frac{\vec{Q}_1\cdot\hat{n}}{\sigma_t}}} \ud s - 2\int_{\Gamma_b} u\, g \ud s \,, \quad \forall u \in Y_p \,. 
	\end{multline}
\begin{response}[orange][ldgexp]
The advantage of LDG (with the choice of $\beta$ given in Eq.~\ref{eq:ldgbeta}) is that any $\kappa \geq 0$, including $\kappa = 0$, results in a stable discretization, avoiding the need to tune a penalty parameter. Additionally, LDG offers the flexibility to control the amount of solution regularization that occurs. For example, setting $\kappa \propto \sigma_t^{-1} p^2/h$ would provide numerical diffusion comparable to IP and BR2 whereas setting $\kappa=0$ provides the so-called minimally dissipative solution. If $\kappa$ is chosen independent of the mesh size and polynomial order, standard preconditioners for discretizations of elliptic problems, such as AMG, will be effective. Otherwise, the specialized preconditioner in \S \ref{sec:subspace} must be used.
These advantages come with the cost that the LDG stabilization term has a non-compact stencil that connects neighbors of neighbors, leading to less sparsity compared to the linear systems associated with the IP and BR2 methods. 
The details of assembling the LDG stabilization terms are provided in \ref{sec:lifting}. 
\end{response}

\subsection{Continuous Finite Element Discretization of VEF}
We now show how a continuous finite element (CG) discretization of the VEF drift-diffusion equation can be extracted from the DG framework presented above. An approximate inversion of this operator is one stage of the subspace correction preconditioner described in \S \ref{sec:subspace} that is used to efficiently solve the IP and BR2 VEF discretizations. This CG operator is also a VEF method itself and represents an extension of the algorithm in \cite{two-level-independent-warsa} to multiple dimensions, high-order, and curved meshes. A CG VEF method has fewer unknowns than an analogous DG method and requires simpler methods to solve the resulting linear system. We will show that this CG discretization has similar accuracy to DG and does not degrade convergence of the fixed-point iteration even in the asymptotic thick diffusion limit. However, it is unclear if using a continuous finite element space would negatively impact robustness and stability in the larger radiation-hydrodynamics multiphysics setting. 

Let $u,\varphi \in V_p$, the degree-$p$ continuous finite element space, then 
	\begin{equation}
		\jump{u} = 0 \,, \quad \jump{\varphi} = 0 \,, \quad \text{on} \ \mathcal{F}\in\Gamma_0 \,. 
	\end{equation}
However, since the Eddington tensor is still discontinuous, we have that 
	\begin{equation}
		\jump{\E\varphi\hat{n}} = \jump{\E\hat{n}} \varphi \,. 
	\end{equation}
Note that, for $u\in V_p$, $\nabla u = \nablah u \in W_p$. In other words, while $u\in V_p$ is continuous $\nabla u$ is not and, due to the continuity properties of functions in $V_p$, the gradient and broken gradient are equivalent \cite[Prop.~3.2.1]{quateroni}. Thus, by starting from the DG VEF discretization and assembling onto $V_p$, we arrive at a CG VEF discretization of the form: find $\varphi \in V_p$ such that 
	\begin{multline} \label{eq:cg}
		\int_{\Gamma_b} E_b\, u \varphi \ud s - \int_{\Gamma_0} \avg{\frac{\nabla u}{\sigma_t}} \cdot \jump{\E\hat{n}}\!\varphi \ud s 
		+ \int \nabla u \cdot \frac{1}{\sigma_t}\nablah\cdot\paren{\E\varphi} \ud \x + \int \sigma_a\, u \varphi \ud \x \\ 
		= \int u\, Q_0 \ud \x + \int \nabla u \cdot \frac{\vec{Q}_1}{\sigma_t} \ud \x - 2\int_{\Gamma_b} u\, g \ud s \,, \quad \forall u \in V_p \,. 
	\end{multline}
Observe that in the thick diffusion limit, where $\E = \frac{1}{3}\I$ and $E_b = 1/2$, a CG discretization of radiation diffusion with Marshak boundary conditions arises since $\jump{\E\hat{n}} = 0$ and $\frac{1}{\sigma_t}\nablah\cdot\paren{\E\varphi} = \frac{1}{3\sigma_t}\nabla \varphi$. 

\section{Subspace Correction Preconditioners} \label{sec:subspace}

In this section, we develop effective and efficient preconditioners for the linear systems resulting from the DG discretizations of the VEF equations developed in \S \ref{sec:dg-vef}.
These preconditioners are built using the additive Schwarz or parallel subspace correction framework \cite{Xu1992,Xu2002}.
We will first discuss the preconditioning of symmetric positive-definite DG discretizations of diffusion equations, and then extend the results to the non-symmetric VEF discretizations.
We begin by reviewing some preliminary results from the domain decomposition literature \cite{Toselli2005}.

\begin{rem}
   In what follows, we will be interested in proving estimates that are independent of discretization parameters such as mesh size $h$, polynomial degree $p$, and penalty parameter $\kappa$.
   For simplicity of notation, we will write $a \lesssim b$ to mean $a \leq C b$, for some constant $C$, independent of $h$, $p$, and $\kappa$.
   Similarly, $a \gtrsim b$ is used to mean $b \lesssim a$, and $a \approx b$ means that both $a \lesssim b$ and $b \lesssim a$.
\end{rem}

We consider a decomposition of the DG finite element space $Y_p$ as the sum of subspaces
\begin{equation}
   Y_p = Y_1 + Y_2 + \cdots + Y_J.
\end{equation}
Let $\A(u,v)$ denote a symmetric positive definite bilinear form, and let $A$ denote the corresponding operator, i.e.
\begin{equation}
   \A(u,v) = (Au, v),
\end{equation}
where $(\cdot\,,\cdot)$ is the standard $L^2(\D)$ inner product.
For example, we can take $\A(u,v)$ to be one of the standard DG discretizations of the diffusion equation as presented in \cite{Arnold2002}.
Let $A_j$ denote the restriction of $A$ to the subspace $Y_j$, and let $P_j$ be the elliptic projections onto $Y_j$.
That is,
\begin{equation}
   \A( P_j u, v_j ) = \A (u, v_j)
      \quad\text{for all $v_j\in Y_j$}.
\end{equation}
Similarly, define the $L^2$ projections onto $Y_j$ by
\begin{equation}
   (Q_j u, v_j ) = (u, v_j)
      \quad\text{for all $v_j\in Y_j$}.
\end{equation}
It can be seen that
\[
   A_j P_j = Q_j A,
\]
and so $P_j = A_j^{-1} Q_j A$.
Inverting the local problems $A_j$ exactly may be computationally infeasible, and so we can replace $A_j^{-1}$ with an approximate inverse $\tilde{A}_j^{-1}$ such that $\tilde{A}_j^{-1} A_j$ is uniformly well-conditioned.
Then, we make use of the operators $T_j = \tilde{A}_j^{-1} Q_j A$.
The \emph{preconditioned operator} $T$ is defined as the sum of the subspace operators, $T = \sum_{j=1}^J T_j$.
The corresponding preconditioner is given by $\sum_{j=1}^J T_j^{-1} Q_j$.
Under certain conditions on the subspaces $Y_j$, the preconditioned system $T = \sum_{j=1}^J T_j^{-1} Q_j A$ is well-conditioned.

\subsection{Decomposition into Conforming and Interface Subspaces}

At this point, we consider the particular decomposition of $Y_p$ into the sum of two subspaces (cf.~\cite{Antonietti2016}),
\begin{equation}
   Y_p = Y_B + V_p,
\end{equation}
where we recall that $V_p \subset Y_p$ consists of functions that are globally continuous, i.e.\ $V_p = Y_p \cap H^1(\mathcal{D})$.
$Y_B$ consists of functions that vanish at all \emph{element-interior} Gauss-Lobatto points (but which may take arbitrary values at element-boundary Gauss-Lobatto points).
This decomposition is closely related to the idea of preconditioning discontinuous Galerkin discretizations with a related continuous Galerkin discretization \cite{Dobrev2006,OMalley2017,Warsa2003}.
It is easy to see that an arbitrary function $w \in Y_p$ has a (non-unique) decomposition as $w = w_b + v$, $w_b \in Y_B$, $v \in V_p$.
Adopting the above notation, let $P_B$ and $P_V$ denote the elliptic projections onto $Y_B$ and $V_p$ respectively.

We recall some results concerning this space decomposition from \cite{Antonietti2016,Pazner2020}.
Let $\A$ denote here the standard interior penalty DG discretization of the diffusion equation.
\begin{prop}[Cf.~\cite{Antonietti2016}, Theorem 1]
   The space decomposition $Y_p = Y_B + V_p$ is stable, i.e.\ for any $w\in Y_p$, there exist a decomposition $w = w_b + v$, $w_b \in Y_B, v \in V_p$ such that
   \begin{equation}
      \A(w_b, w_b) + \A(v, v) \lesssim \A(w, w).
   \end{equation}
   As a consequence of Lions' lemma \cite{Lions1988}, we have
   \begin{equation}
      \A(w, w) \lesssim \A(P_B w, w) + \A(P_V w, w).
   \end{equation}
\end{prop}

An upper bound on $\A(P v_h, v_h)$, where $P = P_B + P_V$ is obtained by noting that the operators $P_B$ and $P_V$ are projections.

\begin{cor}
   The preconditioned operator $P = P_B + P_V$ is uniformly well-conditioned.
\end{cor}

Notice that the operator $A$ restricted to the continuous space $V_p$ corresponds to a standard $H^1$ discretization of the diffusion equation.
As a result, the local solver $A_V^{-1}$ can be replaced with any uniform preconditioner $\tilde{A}_V^{-1}$ for diffusion problems to obtain the approximate operator $T_V$.
For instance, we can take $\tilde{A}_V^{-1}$ to be one V-cycle of \textit{hypre}'s BoomerAMG \cite{hypre}.

It remains to find an approximate solver for the operator $A_B$.
Suppose the mesh $\mathcal{T}$ is conforming, and the space $Y_p$ has constant polynomial degree.
Let $\tilde{A}_B^{-1}$ be the simple point Jacobi preconditioner applied to $A_B$.
Then, we have the following result from \cite{Antonietti2016}.

\begin{thm} \label{thm:as-precond}
   Let $T_B = \tilde{A}_B^{-1} Q_B A$ and let $T_V = \tilde{A}_V^{-1} Q_V A$, where $\tilde{A}_B^{-1}$ is the point Jacobi preconditioner for $A_B$, and $\tilde{A}_V^{-1}$ represents one V-cycle of BoomerAMG (or any other uniform preconditioner for the $H^1$-conforming discretization of diffusion).
   Then, the preconditioned operator $T = T_B + T_V$ is uniformly well-conditioned.
\end{thm}

\begin{rem}
	When the mesh $\mathcal{T}$ is nonconforming (e.g.~as the result of adaptive mesh refinement), or when the DG finite element space $Y_p$ has variable polynomial degrees, then a more sophisticated subspace decomposition is required \cite{Pazner2021}.
	In this case, the boundary subspace $Y_B$ is decomposed into a collection of smaller subspaces defined on each non-conforming edge.
	Each of these small subspaces is solved independently, giving rise to a block Jacobi-type method.
	In the case that the mesh is conforming and the polynomial degree is constant, this construction reduces to the point Jacobi approximate solver described above.
\end{rem}

\subsection{Symmetric VEF Discretizations}

We extend the analysis of the above preconditioners to the family of DG discretizations of the VEF equations given by Eq.~\ref{eq:family_alpha}.
We first treat the simple case where $\E = \frac{1}{3}\I$.
In this case, the system defined by Eq.~\ref{eq:family_alpha} is symmetric and positive-definite.
These results can also be extended to the more general case of constant Eddington tensor; in this case, the results will depend on the spectrum of $\E$.
Let $\B(u,v)$ denote the bilinear form defined by Eq.~\ref{eq:family_alpha}.
We consider the subspace correction preconditioner defined above, and seek to extend Theorem \ref{thm:as-precond} to this modified system.
In order to do this, we must first show that the decomposition $Y_p = Y_B + V_p$ is stable with respect to the modified bilinear form $\B$.
To do this, it suffices to show that the norm induced by $\B$ is equivalent to the norm induced by $\A$.
We first note that the standard interior penalty DG discretization of the definite Helmholtz operator $\sigma_a u - \nabla \cdot ( \sigma_t^{-1} \nabla u)$ satisfies the following bounds (cf.\ \cite{Antonietti2016,Antonietti2010})
\begin{align*}
	\A(u, v) &\lesssim \viii{u} \viii{v}, \\
	\A(u, u) &\gtrsim \viii{u}^2,
\end{align*}
where the mesh-dependent DG norm $\viii{\cdot}$ is defined by
\[
	\viii{u}^2 = \| \sigma_a u \|_0^2 + \| \sigma_t^{-1/2} \nabla_h u \|_0^2 + \frac{p^2}{h} \| \sigma_t^{-1/2} \llb u \rrb \|_{0,\Gamma}^2.
\]
We first consider the interior penalty version of the VEF discretization, given by Eq.~\ref{eq:ip}.
It is straightforward to see that $\B$ satisfies the same inequalities,
\begin{align*}
	\B(u, v) &\lesssim \viii{u} \viii{v}, \\
	\B(u, u) &\gtrsim \viii{u}^2.
\end{align*}
The extension to BR2 and LDG discretizations follows from estimates of the lifting operators $\bm \rho_g$, $\bm r$, and $\bm \ell$, which are considered in \cite{Brezzi2000,Antonietti2010,Pazner2020}.

As a consequence of this equivalence in norms, we expect the parallel subspace preconditioner described above to result in a uniformly well-conditioned operator, independent of mesh size $h$, polynomial degree $p$ (as well as the size of the interior penalty stabilization penalty parameter $\kappa$).

\subsection{Non-Symmetric VEF Discretizations}

The case of more general Eddington tensor $\E$ is more difficult to treat because the resulting bilinear form $\B$ is no longer symmetric.
We analyze the convergence of the preconditioned GMRES iterative method, with the preconditioner defined by the parallel subspace correction procedure described above.
The rate of convergence of the GMRES method applied to a non-symmetric, but positive definite operator is controlled by the ratio of the minimal eigenvalue of the symmetric part of the operator to the norm of the operator \cite{Eisenstat1983}.
We wish to show that this ratio remains independent of the discretization parameters, and therefore that the number of GMRES iterations required to converge remains uniformly bounded.
To do this, recalling the literature on additive Schwarz methods for non-symmetric problems (cf.\ \cite{Cai1989,Cai1990}), we must show that the non-symmetric part of the operator is small in some sense.

In order to simplify the analysis, we consider a slightly modified VEF discretization that results from iteratively lagging certain non-symmetric terms.
In particular, we write $\nablah\cdot\paren{\E\varphi} = \E\nablah \varphi + \paren{\nablah\cdot\E}\!\varphi$, and iteratively lag the second term on the right-hand side, replacing $\paren{\nablah\cdot\E}\!\varphi$ with $\paren{\nablah\cdot\E}\!\hat{\varphi}$, where $\hat{\varphi}$ is given from the previous iteration.
The iteratively lagged version of Eq.~\ref{eq:ip} then gives rise to
\begin{equation}
\begin{multlined}[c][.9\displaywidth] \label{eq:lagged}
	\B(u, \varphi) =
	\int_{\Gamma_b} E_b\, u \varphi \ud s + \int_{\Gamma_0} \kappa \jump{u} \jump{\varphi} \ud s - \int_{\Gamma_0} \jump{u} \avg{\frac{1}{\sigma_t}\E\nablah\varphi\cdot\hat{n}} \ud s \\
	- \int_{\Gamma_0} \avg{\frac{\nablah u}{\sigma_t}} \cdot \jump{\E\varphi\hat{n}} \ud s
	+ \int \nablah u \cdot \frac{1}{\sigma_t} \E \nablah\varphi \ud \x + \int \sigma_a\, u \varphi \ud \x.
\end{multlined}
\end{equation}

We decompose $\B$ into its symmetric and skew-symmetric parts, $\B(u,\varphi) = \SS(u,\varphi) + \N(u,\varphi)$, where
\begin{gather*}
   \SS(u,\varphi) = \frac{1}{2} \left( \B(u,\varphi) + \B(\varphi,u) \right), \\
   \N(u,\varphi) = \frac{1}{2} \left( \B(u,\varphi) - \B(\varphi,u) \right).
\end{gather*}
Cf.\ Theorem 1.3 from \cite{Cai1989}, preconditioned GMRES will converge uniformly with respect to the discretization parameters if there exists a constant $0 \leq \delta < 1$ such that
\begin{equation}
  \label{eq:skew-sym-bound}
   \left| \N(u, Pu) \right| \leq \delta \B(u, P u),
\end{equation}
where $P = P_B + P_V$ is the preconditioned operator.
We see that the skew-symmetric part of Eq.~\ref{eq:lagged} is given by
\begin{equation}
  \label{eq:skew-sym}
  \begin{multlined}[c][.9\displaywidth]
	\N(u, \varphi) =
	\frac{1}{2}\Bigg(
		- \int_{\Gamma_0} \jump{u} \avg{\frac{1}{\sigma_t}\E\nablah\varphi\cdot\hat{n}} \ud s
		+ \int_{\Gamma_0} \jump{\varphi} \avg{\frac{1}{\sigma_t}\E\nablah u\cdot\hat{n}} \ud s \\
		- \int_{\Gamma_0} \avg{\frac{\nablah u}{\sigma_t}} \cdot \jump{\E\varphi\hat{n}} \ud s
		+ \int_{\Gamma_0} \avg{\frac{\nablah \varphi}{\sigma_t}} \cdot \jump{\E u\hat{n}} \ud s
	\Bigg).
  \end{multlined}
\end{equation}
Applying the identity
$\jump{ ac } \avg{ b } - \jump{ a } \avg{ bc } = \frac{1}{2} \jump{c} \paren{ a_1 b_2 + a_2 b_1 }$
to the above expression yields the following boundedness property
\[
	\left| \N(u, \varphi) \right| \lesssim \jump{\E} \viii{u} \, \viii{\varphi},
\]
where $\jump{\E}$ represents an upper bound on the jump of $\E$ over all element interfaces in the mesh.
Using that $Y_p = Y_B + V_p$ is a stable decomposition, we have
\begin{align*}
  \B(u, Pu)
  = \B(u, P_B u) + \B(u, P_V u)
  &= \B(P_B u, P_B u) + \B(P_V u, P_V u) \\
  &= \SS(P_B u, P_B u) + \SS(P_V u, P_V u)
  \gtrsim \viii{u}^2.
\end{align*}
Furthermore, since $P_B$ and $P_V$ are projections,
\[
  \viii{Pu} = \viii{P_B u + P_V u} \leq 2 \viii{u}.
\]
Combining the above estimates, we obtain
\begin{align*}
  \left|\N(u, Pu)\right|
  \lesssim \jump{\E} \viii{ u } \, \viii{ Pu }
  \lesssim \jump{\E} \B(u, Pu).
\end{align*}
Therefore, in order to obtain the bound \eqref{eq:skew-sym-bound} with $0 \leq \delta < 1$, according to the size of the jumps $\llb \E \rrb$, we may choose $\kappa$ sufficiently large in the symmetric penalty term
\[
	\int_{\Gamma_0} \kappa\,\llb u \rrb \llb \varphi \rrb \, \ud s.
\]
Having chosen $\kappa$ to satisfy this bound, preconditioned GMRES applied to this system will converge uniformly, independent of the discretization parameters.

\begin{rem}
	While the GMRES convergence estimates shown in this section apply in the case of the modified (iteratively lagged) VEF discretization with sufficiently large penalty parameter, in practice we observe uniform convergence for the non-lagged VEF discretizations, without additional conditions on the size of the penalty parameter $\kappa$.
	This behavior is typical of domain decomposition algorithms applied to non-symmetric and indefinite problems, for which the theoretical convergence estimates tend to be pessimistic \cite{Toselli2005}.
\end{rem}

\begin{rem}[AMG convergence] \label{rem:amg}
	The practical subspace correction preconditioner is obtained by replacing $A_V^{-1}$ (the inverse of the continuous discretization, which is infeasible to compute for large problems) with a tractable approximation $\tilde{A}_V^{-1}$, such as one V-cycle of algebraic multigrid, cf.~Theorem \ref{thm:as-precond}.
	This procedure relies on $\tilde{A}_V^{-1}$ well approximating $A_V^{-1}$ (i.e.~spectrally equivalent in the symmetric case).
	AMG performance may suffer on highly non-symmetric problems, and so in the following sections, we consider also choosing $\tilde{A}_V^{-1}$ to be one V-cycle of AMG built with a symmetrized version of the continuous operator $A_V$.
\end{rem}

\section{Results}
\resp[red]{We now present numerical results concerning the iterative efficiency and computational performance of the outer fixed-point iteration and inner preconditioned iterative solvers for each of the discretizations of the VEF equations discussed above. The outer iteration refers to the evaluation of the fixed-point operator $G(\varphi)$ defined in Eq.~\ref{eq:fp_op} which includes inverting the streaming and collision terms in the transport equation and solving the discrete VEF equations. The inner iteration refers to solving the discrete VEF equations iteratively. Each inner iteration requires applying the matrix operator and preconditioner corresponding to the VEF discretization. }

The VEF algorithms described in this paper were implemented using the MFEM \cite{MFEM,mfem-web} finite element framework. The stabilized bi-conjugate gradient (BiCGStab) and Jacobi solvers from MFEM were used to solve the VEF discretizations along with BoomerAMG, the AMG solver from the sparse linear algebra library \emph{hypre} \cite{hypre}. \resp[orange]{Note that BiCGStab performed equivalently to GMRES and thus we elect to use BiCGStab because it does not require storage of a Krylov space.} KINSOL, from the Sundials package \cite{hindmarsh2005sundials}, provided the Anderson-accelerated fixed-point solver. \resp[blue][kinsol]{As described in \citet[\S 2]{hindmarsh2005sundials}, the fixed-point and Anderson-accelerated fixed-point iteration is terminated when the max norm of the difference between successive iterates is below the iterative tolerance.} When iterative solver results are not presented, the parallel implementation of the sparse direct solver SuperLU \cite{lidemmel03} was used. We use the high-order DG \Sn transport solver from \cite{graph_sweeps}.

Unless otherwise specified, we set the penalty parameter to 
	\begin{equation}
		\kappa = \avg{\frac{(p+1)^2}{\sigma_t h}} 
	\end{equation}
and the BR2 stabilization parameter to $\eta = 4$. \resp[orange][penstand]{These choices are standard in the literature for the model elliptic problem \cite{Bassi2000,Arnold1982} and are the default choices implemented for discretizations of the Poisson equation in MFEM.} We use the MDLDG method, the variant of the LDG method where $\kappa \equiv 0$ and set the upwinding vector $\vec{w}$ to be a unit vector at a $45^\circ$ angle from the $x$-axis. The VEF discretizations all use the element-local basis defined using the Gauss-Lobatto points to enable the use of the subspace correction preconditioner where required. The transport discretization is always solved with the same finite element order as the VEF scalar flux. However, we use the positive Bernstein basis \cite{doi:10.1137/11082539X} for the transport discretization. A summary of the properties associated with each VEF discretization is presented in Table \ref{tab:disc_sum}. 
\begin{table}
\centering
\begin{response}[orange]
\caption{A summary of the key algorithmic properties of the VEF discretizations.}
\label{tab:disc_sum}
\begin{tabular}{c cccc}
	\toprule 
	& IP & BR2 & MDLDG & CG \\ 
	\midrule 
	Solution Space & $Y_p$ & $Y_p$ & $Y_p$ & $V_p$ \\
	Penalty scales with mesh size & Yes & Yes & No & -- \\
	Local Stencil & Yes & Yes & No & Yes \\
	Requires specialized preconditioner & Yes & Yes & No & No \\ 
	\bottomrule
\end{tabular}
\end{response}
\end{table}

\subsection{Method of Manufactured Solutions} \label{sec:mms}
The accuracy of the methods are ascertained with the Method of Manufactured Solutions (MMS). The solution is set to 
	\begin{equation}
		\psi = \frac{1}{4\pi}\paren{\sin\paren{\pi x}\sin\paren{\pi y} + \Omegahat_x\Omegahat_y \sin\paren{2\pi x} \sin\paren{2\pi y} + \Omegahat_x^2 \sin\paren{\frac{3\pi(x+\delta)}{1+2\delta}}\sin\paren{\frac{3\pi(y+\delta)}{1+2\delta}} + \gamma} 
	\end{equation}
where the parameters $\delta = 0.1$ and $\gamma = 0.5$ control the amount of spatially varying, quadratically anisotropic inflow and uniform, isotropic inflow, respectively. The computational domain is $\D = [0,1]^2$. With this solution, the Eddington tensor varies in space and has non-zero off-diagonal components. Trigonometric functions are used so that the solution cannot be exactly represented by polynomials. The scalar flux is then 
	\begin{equation}
		\phi = \sin\paren{\pi x} \sin\paren{\pi y} + \frac{1}{3} \sin\paren{\frac{3\pi(x+\delta)}{1+2\delta}}\sin\paren{\frac{3\pi(y+\delta)}{1+2\delta}} + \gamma \,. 
	\end{equation}
These MMS angular and scalar flux solution functions are substituted into the transport equation to solve for the MMS source function.

The accuracy of the VEF discretizations can be investigated in isolation by computing the VEF data from the MMS angular flux and setting the sources $Q_0$ and $\vec{Q}_1$ to the moments of the transport MMS source. This is accomplished by computing the VEF data from the MMS angular flux projected onto a finite element space of equal order to the VEF finite element space. An open, Gauss-Legendre basis is used for the angular flux so that the Eddington tensor has discontinuities of magnitude $\mathcal{O}(h^{p+1})$ on interior mesh faces. The VEF data and source moments are computed using level symmetric $S_4$ angular quadrature. The VEF equations are then solved as if $\E$, $E_b$, $Q_0$, and $\vec{Q}_1$ are provided data. 

\begin{figure}
\centering
\begin{subfigure}{.49\textwidth}
	\centering
	\includegraphics[height=2in]{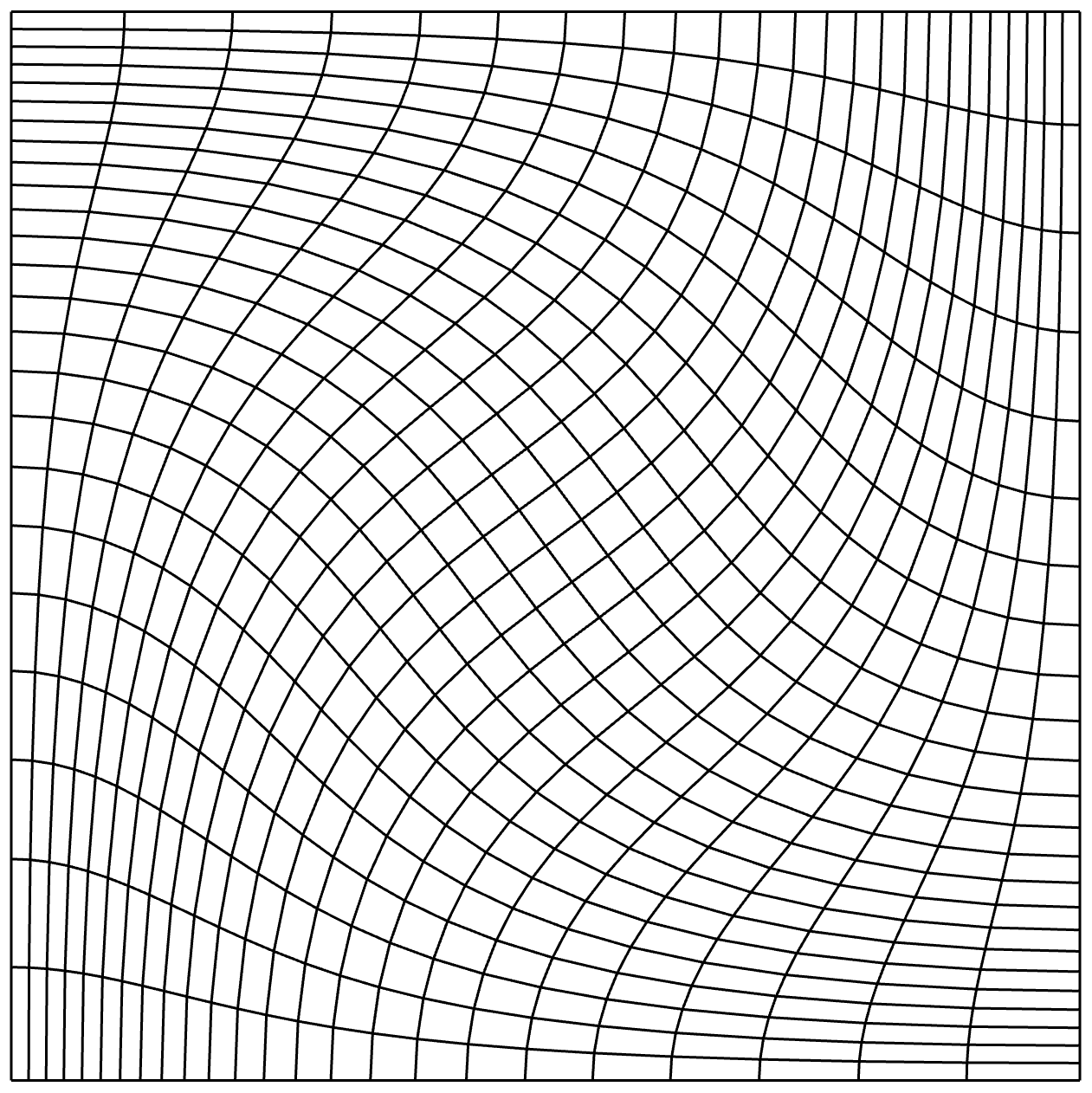}
	\caption{}
	\label{fig:tgmesh}
\end{subfigure}
\begin{subfigure}{.49\textwidth}
	\centering
	\includegraphics[width=\textwidth]{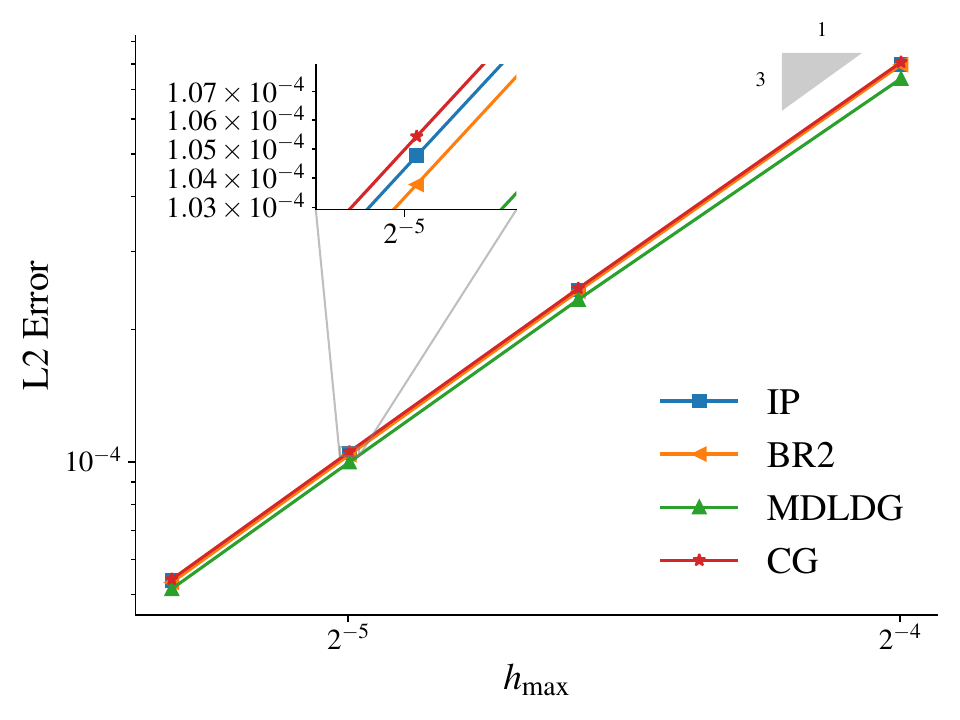}
	\caption{}
	\label{fig:mms2}
\end{subfigure}
\caption{(a) An example third-order mesh distorted according to the Taylor Green hydrodynamics solution. (b) The MMS error as the mesh is refined for each VEF method when $p=2$. Each method converges with third-order accuracy. }
\end{figure}

We use refinements of a third-order curved mesh created by distorting an orthogonal mesh according to the velocity field of the Taylor Green vortex. This mesh distortion is generated by advecting the mesh control points with 
	\begin{equation}
		\x = \int_0^T \mat{v} \ud t \,,
	\end{equation}
where the final time $T = 0.3 \pi$ and 
	\begin{equation}
		\mat{v} = \begin{bmatrix} 
			\sin(x_1) \cos(x_2) \\
			-\cos(x_1)\sin(x_2) 
		\end{bmatrix}
	\end{equation}
is the analytic solution of the Taylor Green vortex. The time integration is calculated with 300 forward Euler time steps. An example mesh is shown in Fig.~\ref{fig:tgmesh}. 

Figure \ref{fig:mms2} shows the $L^2(\D)$ error between the VEF solution and the exact MMS scalar flux solution as the mesh is refined for the IP, BR2, MDLDG, and CG VEF discretizations when quadratic basis functions are used. Here, $h_\text{max}$ is the maximum value of the characteristic element length in the mesh. All methods have nearly identical error behavior and converge with third-order accuracy as expected. This experiment is repeated with $p=3$ in Table \ref{tab:mms3}. Logarithmic regression is used to compute the exponent and constant of the error function $E = C h_\text{max}^{\bar{p}}$ with $C$ the constant and $\bar{p}$ the method's experimentally observed order of accuracy. The standard deviation of the four error values for each mesh is also provided to quantify the variance in the error behavior. Accuracy of $\mathcal{O}(h^{p+1})$ is observed and the four variants are shown to have variance below the discretization error. 

\begin{table}
\centering
\caption{MMS error for each method as a function of the maximum characteristic mesh size, $h_\text{max}$, with $p=3$. The standard deviation of the four error values in each row is also provided showing that differences between each method are below the discretization error. The order of accuracy and error constant were computed with logarithmic regression. }
\label{tab:mms3}
\begin{tabular}{cc cccc cc}
	\toprule
	$h_\mathrm{max}$ && IP & BR2 & MDLDG & CG && Deviation \\ 
	\midrule 
	\num{8.345e-02} &  & \num{2.678e-04} & \num{2.676e-04} & \num{2.598e-04} & \num{2.688e-04} &  & \num{3.622e-06} \\
\num{5.564e-02} &  & \num{5.163e-05} & \num{5.158e-05} & \num{4.837e-05} & \num{5.169e-05} &  & \num{1.415e-06} \\
\num{4.173e-02} &  & \num{1.631e-05} & \num{1.630e-05} & \num{1.505e-05} & \num{1.632e-05} &  & \num{5.463e-07} \\
\num{3.338e-02} &  & \num{6.684e-06} & \num{6.680e-06} & \num{6.115e-06} & \num{6.686e-06} &  & \num{2.459e-07}  \\\addlinespace
	Order &  & 4.028 & 4.028 & 4.092 & 4.032 \\
Constant &  & 5.885 & 5.882 & 6.678 & 5.966  \\
	\bottomrule
\end{tabular}
\end{table}

\subsection{Thick Diffusion Limit} \label{sec:tdl}
Next, we investigate the iterative convergence properties of the VEF methods in the regime known as the asymptotic thick diffusion limit \cite{diflim}. The material data are set to:
	\begin{equation} \label{eq:tdl_scaling}
		\sigma_t = \frac{1}{\epsilon} \,, \quad \sigma_a = \epsilon \,, \quad \sigma_s = \frac{1}{\epsilon} - \epsilon \,, \quad q = \epsilon 
	\end{equation}
with $\epsilon \in (0,1]$ and the thick diffusion limit corresponding to the limit $\epsilon \rightarrow 0$. A coarse mesh that does not adequately resolve the mean free path is used to stress the convergence of the VEF algorithm. This is a numerically challenging, but common in practice, regime where robust performance is crucial. 

We first demonstrate robust convergence on an $8\times 8$ linear mesh with $\D = [0,1]^2$. Convergence was identical for linear, quadratic, and cubic basis functions so we present results for $p=2$ only. Level symmetric $S_4$ angular quadrature is used. Fixed-point iteration without Anderson acceleration is used to solve the coupled transport-VEF system. 

\begin{figure}
	\centering
	\begin{subfigure}{.24\textwidth}
		\centering
		\includegraphics[width=\textwidth]{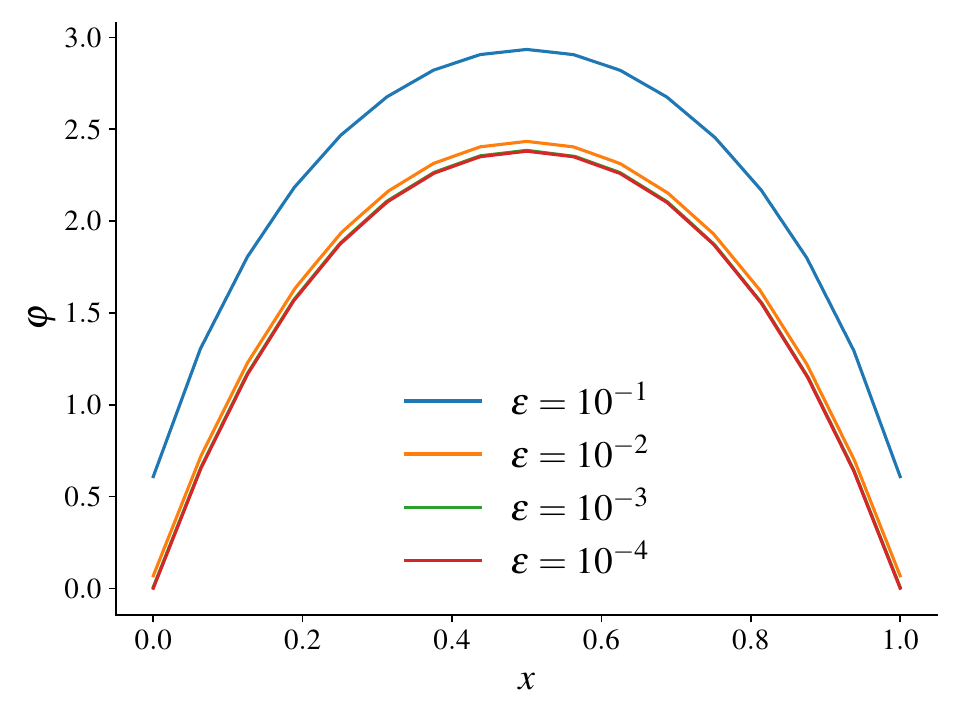}
		\caption{}
	\end{subfigure}
	\begin{subfigure}{.24\textwidth}
		\centering
		\includegraphics[width=\textwidth]{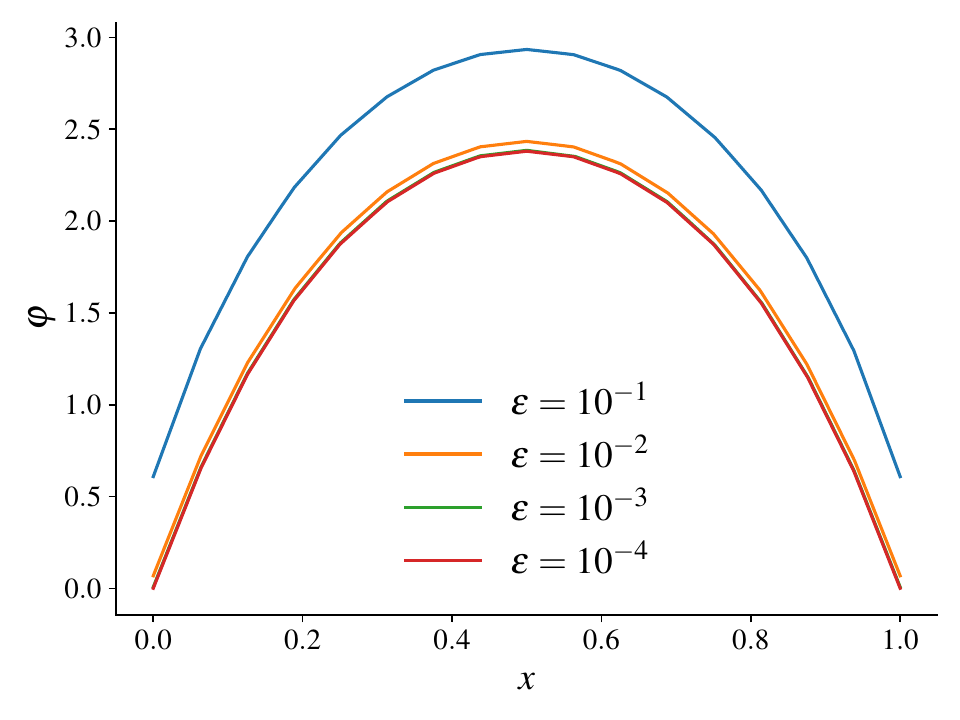}
		\caption{}
	\end{subfigure}
	\begin{subfigure}{.24\textwidth}
		\centering
		\includegraphics[width=\textwidth]{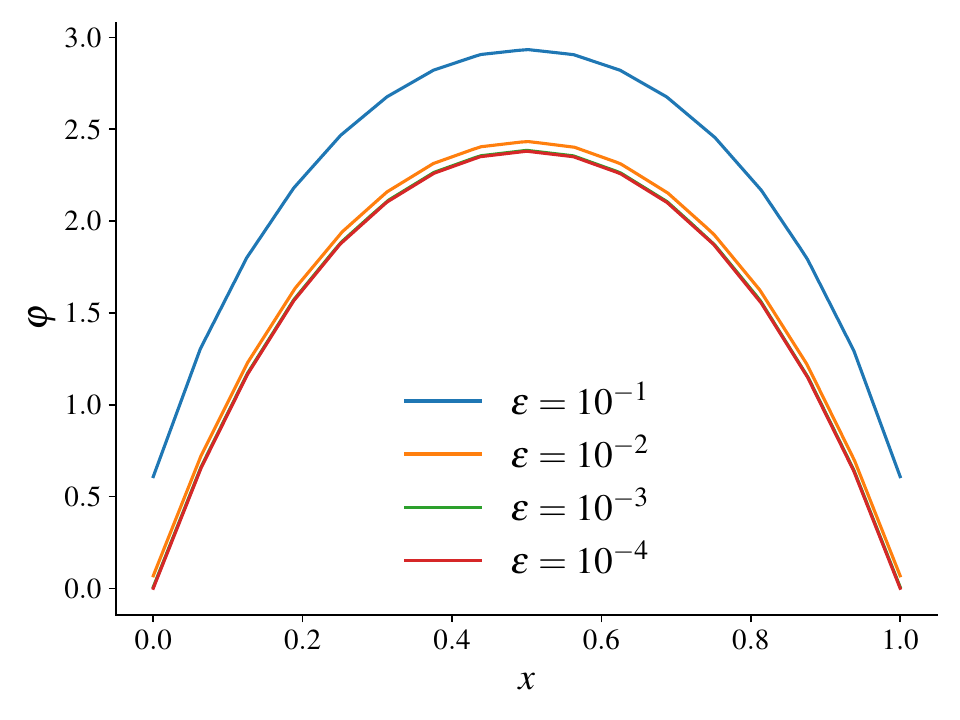}
		\caption{}
	\end{subfigure}
	\begin{subfigure}{.24\textwidth}
		\centering
		\includegraphics[width=\textwidth]{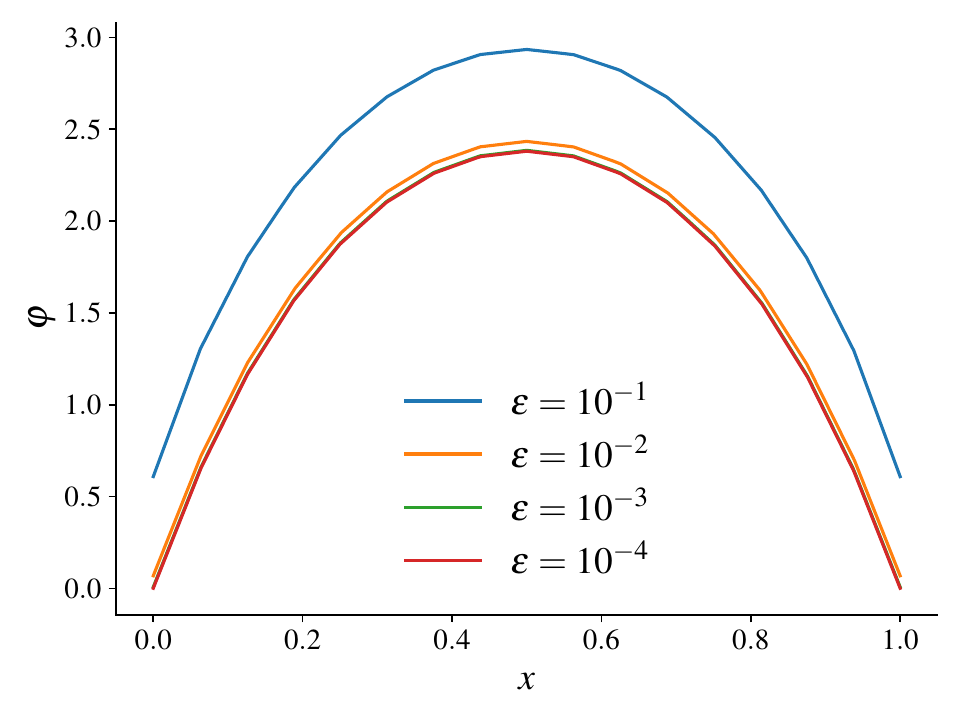}
		\caption{}
	\end{subfigure}
	\caption{Lineouts of the 2D thick diffusion limit solutions taken at $y=\frac{1}{2}$ for the (a) IP, (b) BR2, (c) MDLDG, and (d) CG methods. All methods converge to the non-trivial, diffusion solution as $\epsilon \rightarrow 0$.}
	\label{fig:eps_lineout}
\end{figure}

Table \ref{tab:tdl_orthog} shows the number of fixed-point iterations required to converge to a tolerance of $10^{-6}$ as $\epsilon \rightarrow 0$. All four VEF variants converged robustly and in an identical number of iterations for each value of $\epsilon$. Lineouts of the 2D solutions are shown in Fig.~\ref{fig:eps_lineout} to demonstrate that the non-trivial, diffusion solution is obtained by each method. Note that even the continuous finite element discretization paired with the discontinuous finite element transport discretization is robust in the thick diffusion limit. 

\begin{table}
\centering
\caption{Number of \resp[blue]{fixed-point} iterations to convergence in the thick diffusion limit on a coarse, orthogonal mesh.}
\label{tab:tdl_orthog}
\begin{tabular}{c cccc}
	\toprule 
	$\epsilon$ & IP & BR2 & MDLDG & CG \\ 
	\midrule 
	$10^{-1}$ & 8 & 8 & 8 & 8 \\
$10^{-2}$ & 6 & 6 & 6 & 6 \\
$10^{-3}$ & 4 & 4 & 4 & 4 \\
$10^{-4}$ & 3 & 3 & 3 & 3 \\
	\bottomrule 
\end{tabular}
\end{table}

\begin{figure}
	\centering
	\includegraphics[width=.65\textwidth]{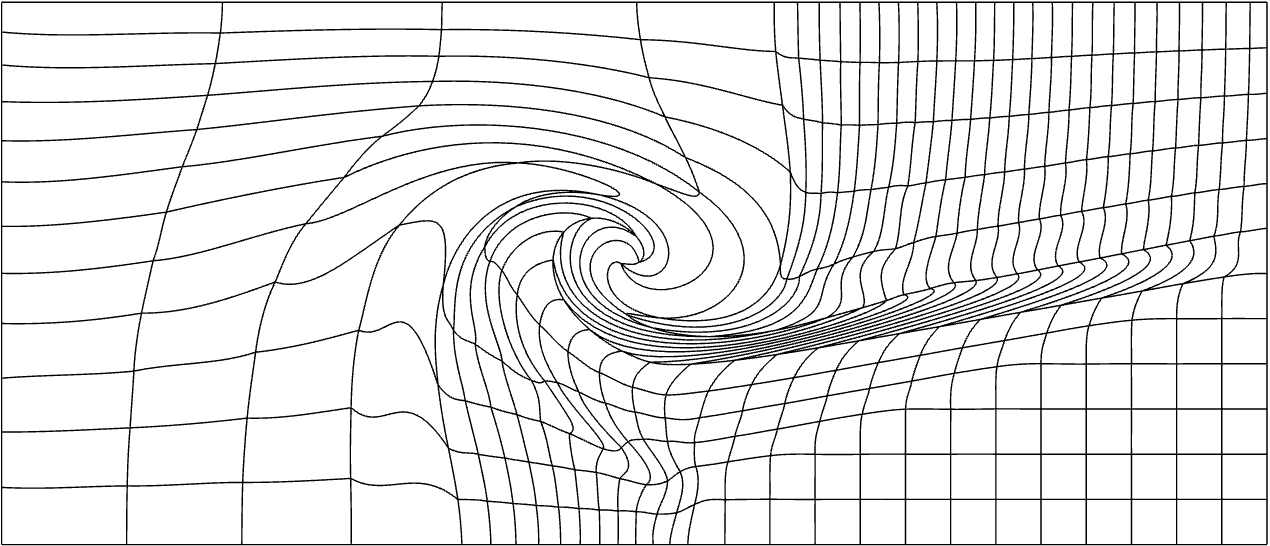}
	\caption{A depiction of the triple point mesh used to stress test the VEF algorithms on a severely distorted, third-order mesh. The mesh was generated with a purely Lagrangian hydrodynamics simulation. }
	\label{fig:3point_mesh}
\end{figure}

This experiment is repeated on the triple point mesh shown in Fig.~\ref{fig:3point_mesh}. This mesh was generated by running a purely Lagrangian hydrodynamics simulation on a third-order mesh. The mesh contains concave/reentrant interior faces meaning the matrix corresponding to the transport discretization cannot be reordered to be strictly lower block triangular. The pseudo-optimally reordered sweep from \cite{graph_sweeps}, which lags the incoming angular flux on reentrant faces, is used to enable an element-by-element transport solve. Since the incoming fluxes on reentrant faces are lagged, the angular flux on these faces is not linearly eliminated. In other words, the presence of reentrant faces means that the transport equation is not fully inverted at every fixed-point iteration. In addition, the mesh elements in the ``swirl'' at the center are severely distorted and thus have poor approximation ability. In practice, the mesh would be remapped before this level of distortion were present. Due to this severe distortion, stability of the IP VEF discretization required scaling the penalty parameter according to 
	\begin{equation}
		\kappa = C \avg{\frac{(p+1)^2}{\sigma_t h}} \,,
	\end{equation}	
where $C = \max_{K_e \in \meshT} C_e$ with $C_e$ the condition number of the Jacobian matrix for element $K_e$. For the triple point mesh, $C=169$. \resp[blue][br2stab]{Note that the BR2 method was stable on the triple point mesh without modifying the parameter $\eta$. This is an example where the increased assembly cost of the BR2 method provides additional robustness compared to the IP method.} \resp[orange][jcond]{However, we have found that this heuristic for scaling the penalty parameter is effective for ensuring stability of the IP VEF method on many meshes with varying levels of distortion.}

Table \ref{tab:tdl_3point} shows the number of fixed-point iterations without Anderson acceleration required to converge to a tolerance of $10^{-6}$ for the four VEF variants as $\epsilon \rightarrow 0$. Fixed-point convergence is shown when one, two, and three lagged transport sweeps are applied per fixed-point iteration. Here, lagged transport sweep refers to inverting the streaming and collision operator using lagged information on reentrant faces. While one lagged sweep per fixed-point iteration required more iterations than the equivalent orthogonal-mesh problem, especially for large values of $\epsilon$, the three lagged sweeps per fixed-point iteration option had similar convergence properties to its orthogonal-mesh counterpart. This suggests the iterative slow-down can be attributed to the approximate sweep. \resp[blue][effsweeps]{While performing more lagged sweeps per fixed-point iteration did improve iterative efficiency, efficiency was not improved to the point that the total number of lagged sweeps was reduced. That is, the three-sweep option, which converged in the fewest iterations, performed the most lagged sweeps.}

\begin{figure}
	\centering
	\begin{subfigure}{.24\textwidth}
		\centering
		\includegraphics[width=\textwidth]{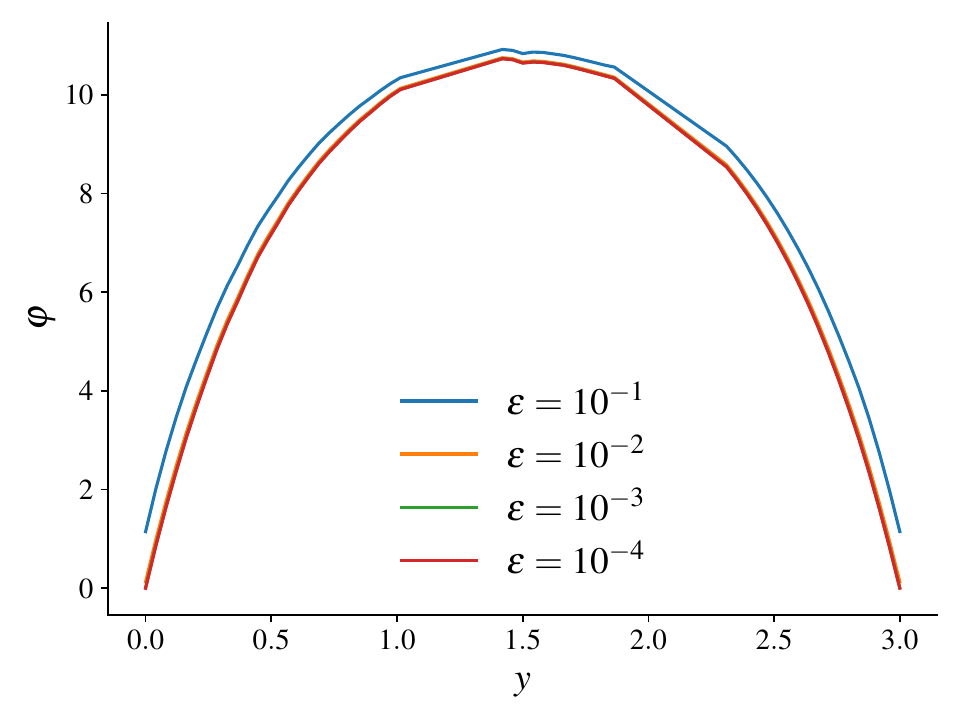}
		\caption{}
	\end{subfigure}
	\begin{subfigure}{.24\textwidth}
		\centering
		\includegraphics[width=\textwidth]{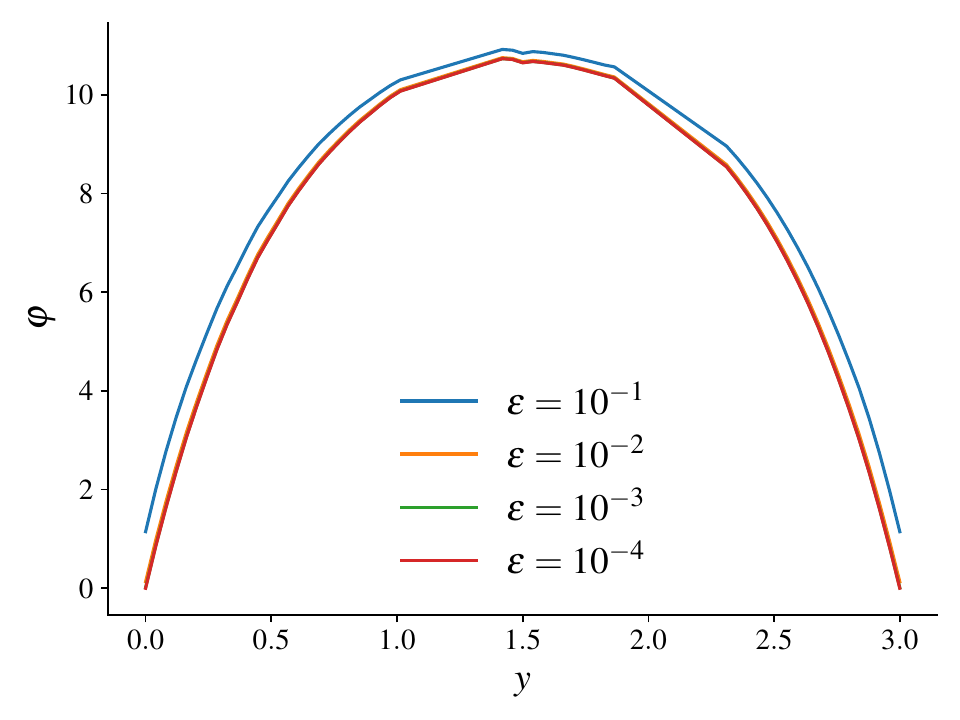}
		\caption{}
	\end{subfigure}
	\begin{subfigure}{.24\textwidth}
		\centering
		\includegraphics[width=\textwidth]{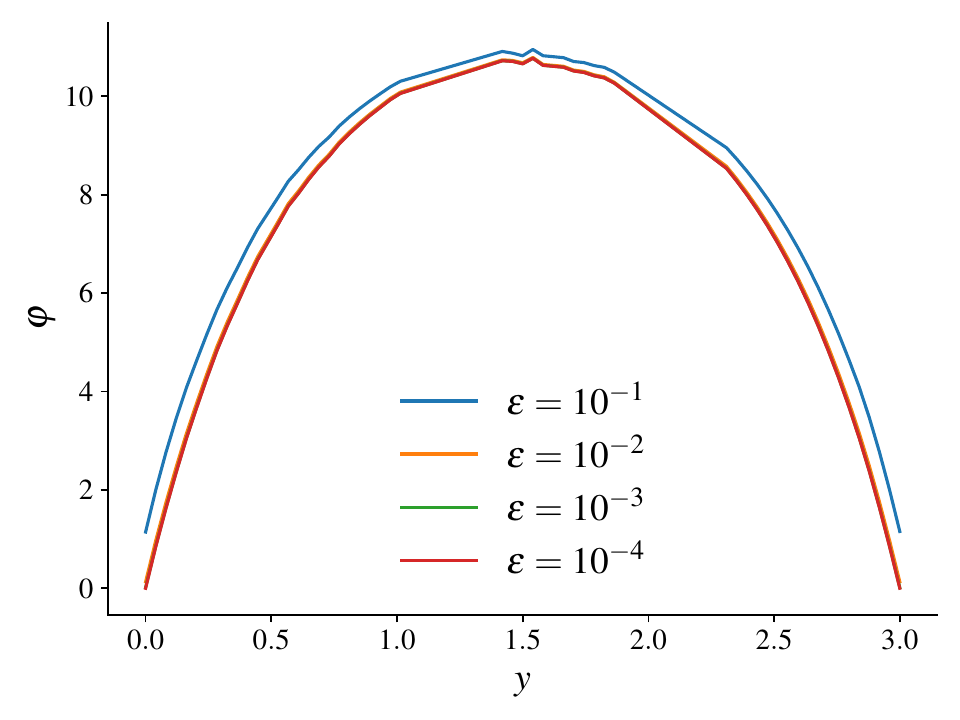}
		\caption{}
	\end{subfigure}
	\begin{subfigure}{.24\textwidth}
		\centering
		\includegraphics[width=\textwidth]{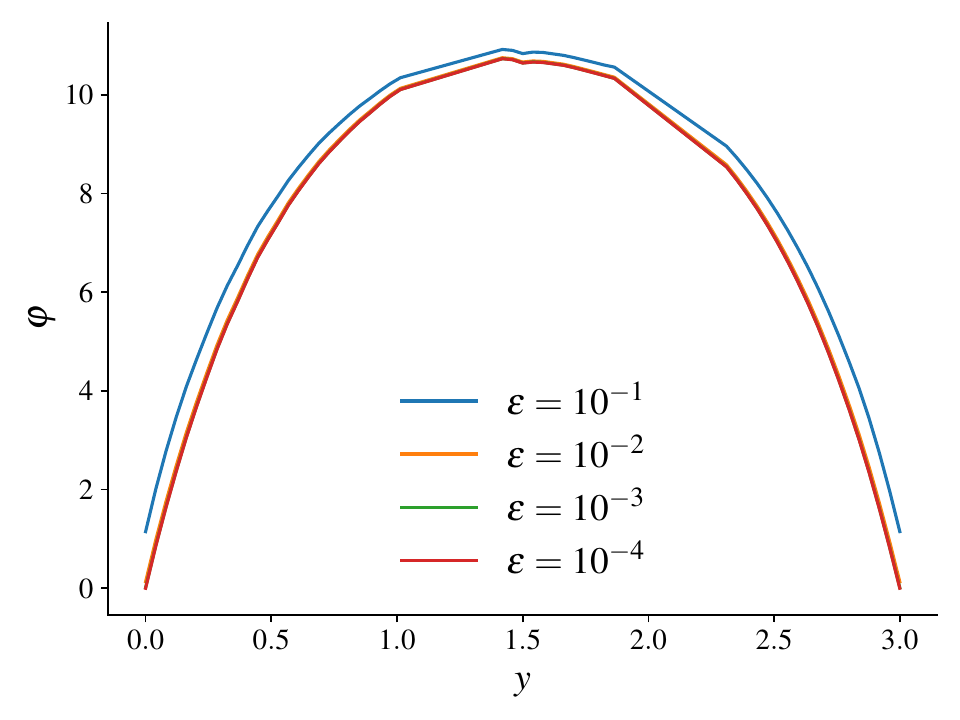}
		\caption{}
	\end{subfigure}
	\caption{\resp[orange]{Lineouts of the 2D thick diffusion limit solutions on the triple point mesh taken at $x=3.5$ for the (a) IP, (b) BR2, (c) MDLDG, and (d) CG methods. Non-monotonic oscillations are observed due to imprinting from the severely distorted mesh.}}
	\label{fig:3p_lineout}
\end{figure}

\resp[orange]{Lineouts of the solutions are provided in Fig.~\ref{fig:3p_lineout} to demonstrate that a non-trivial solution was obtained even on the distorted triple point mesh. The solutions have non-physical, non-monotonic oscillations due to imprinting of the severely distorted mesh.
We present the solutions generated by the one sweep per fixed-point iteration option only as the two and three sweep options converged to equivalent solutions. }

Table \ref{tab:tdl_anderson} shows the diffusion scaling on the triple point mesh for the IP VEF method with Anderson acceleration. An Anderson space of size $a$ is used where $a=0$ is equivalent to fixed-point iteration. We compare convergence when the Anderson space is built from the scalar flux only and when it is built from the scalar and angular fluxes. These variants are referred to as ``low memory'' and ``augmented'', respectively. Note that to simplify the implementation, the augmented Anderson space is built from the entire angular flux and not just the subset of angular flux unknowns corresponding to reentrant faces. The augmented variant saw improvement for $\epsilon = 10^{-1}$ but otherwise converged equivalently to fixed-point iteration. The low memory option was not improved by Anderson acceleration and actually took 1-3 more iterations to converge. Since convergence is primarily hindered by the inexact transport inversion, it is expected that Anderson cannot improve convergence when the Anderson space is not augmented with the angular flux. 

\begin{table}
\centering
\caption{Number of fixed-point iterations required for convergence on the triple point mesh as $\epsilon \rightarrow 0$. On the triple point mesh, the presence of reentrant faces mean the streaming and collision operator is not fully inverted at each iteration. Each method is tested with 1, 2, and 3 approximate transport inversions per fixed-point iteration. }
\label{tab:tdl_3point}
\begin{tabular}{c c ccc c ccc c ccc c ccc}
	\toprule
	&& \multicolumn{3}{c}{IP} && \multicolumn{3}{c}{BR2} && \multicolumn{3}{c}{MDLDG} && \multicolumn{3}{c}{CG} \\
	\cmidrule{3-5} \cmidrule{7-9} \cmidrule{11-13} \cmidrule{15-17}
	$\epsilon$ && 1 & 2 & 3 && 1 & 2 & 3 && 1 & 2 & 3 && 1 & 2 & 3\\
	\midrule 
	$10^{-1}$ &  & 19 & 11 & 10 &  & 19 & 11 & 10 &  & 23 & 14 & 12 &  & 19 & 11 & 10 \\
$10^{-2}$ &  & 11 & 8 & 7 &  & 11 & 8 & 7 &  & 19 & 9 & 8 &  & 11 & 8 & 7 \\
$10^{-3}$ &  & 8 & 5 & 4 &  & 8 & 5 & 4 &  & 9 & 5 & 5 &  & 8 & 5 & 4 \\
$10^{-4}$ &  & 6 & 4 & 3 &  & 6 & 4 & 3 &  & 6 & 4 & 3 &  & 6 & 4 & 3 \\
	\bottomrule
\end{tabular}
\end{table}

\begin{table}
\centering
\caption{The number of Anderson-accelerated fixed-point iterations, with Anderson space of size $a$, required for convergence on the triple point mesh. The interior penalty VEF method is used. The low memory variant builds the Anderson space from the VEF scalar flux only while the augmented version builds the Anderson space from the VEF scalar flux and the angular flux. }
\label{tab:tdl_anderson}
\begin{tabular}{c c cc c cc}
	\toprule 
	&& \multicolumn{2}{c}{Low Memory} && \multicolumn{2}{c}{Augmented} \\ 
	\cmidrule{3-4} \cmidrule{6-7} 
	$\epsilon$ && $a=0$ & $a=5$ && $a=0$ & $a=5$ \\
	\midrule
	$10^{-1}$ &  & 19 & 20 &  & 19 & 14 \\
$10^{-2}$ &  & 11 & 13 &  & 11 & 11 \\
$10^{-3}$ &  & 8 & 11 &  & 8 & 8 \\
$10^{-4}$ &  & 6 & 8 &  & 6 & 6 \\
	\bottomrule
\end{tabular}
\end{table}

\subsection{Linearized Crooked Pipe} \label{sec:lcp}
\begin{figure}
	\centering
	\includegraphics[width=.65\textwidth]{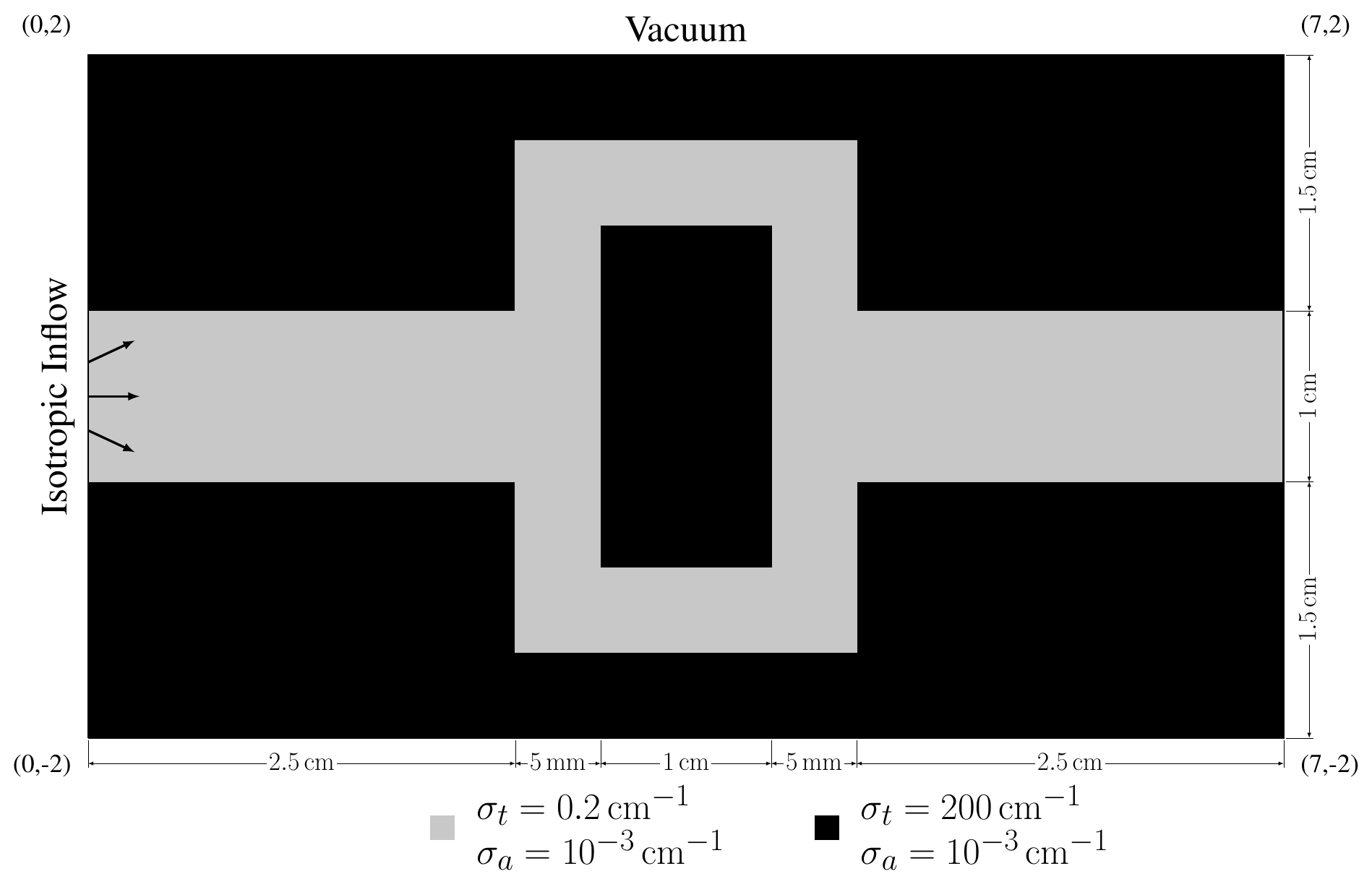}
	\caption{Geometry, material data, and boundary conditions for the linearized crooked pipe problem.}
	\label{fig:cp_diag}
\end{figure}
We now demonstrate the efficacy of the methods on a more realistic, multi-material problem. A common benchmark is the crooked pipe problem. The geometry and materials are shown in Fig.~\ref{fig:cp_diag}. The problem consists of two materials, the wall and the pipe, which have an 1000x difference in total interaction cross section. We mock the time-dependent benchmark as a steady-state problem by adding artificial absorption and fixed-source terms corresponding to backward Euler time integration. A large time step such that $c\Delta t = 10^3$ is used with an initial condition $\psi_0 = 10^{-4}$ for all $(\x,\Omegahat) \in \D \times \mathbb{S}^2$. Thus, the absorption and source terms are
	\begin{subequations}
	\begin{equation}
		\sigma_a = \frac{1}{c\Delta t} = 10^{-3} \,\si{\per\cm}\,,  
	\end{equation}
	\begin{equation}
		q = \frac{1}{c\Delta t} \psi_0 = 10^{-1}\,\si{\per\cm\cubed\per\s\per\str} \,. 
	\end{equation}
	\end{subequations}
The boundary conditions are
	\begin{equation}
		f = \begin{cases}
			\frac{1}{2\pi}\,, & x = 0 \ \mathrm{and}\ y \in [-1/2,1/2] \\ 
			0 \,, & \mathrm{otherwise}
		\end{cases} \,,
	\end{equation}
so that radiation enters the pipe at the left side of the domain. We use a Level Symmetric $S_{12}$ angular quadrature set. The zero and scale \cite{hamilton2009negative} negative flux fixup -- a sweep-compatible method that zeros out negativity and rescales so that particle balance is preserved -- is used inside the inversion of the streaming and collision operator to ensure positivity. 

The efficiency of the outer fixed-point and inner linear iterations is investigated by refining in $h$ and $p$ on a uniform mesh of quadrilateral elements that is aligned with the materials. The outer solver is Anderson-accelerated fixed-point iteration with two Anderson vectors. Anderson acceleration is not required for convergence on this problem but does provide more uniform convergence in $h$. Since the mesh is orthogonal, the transport equation is fully inverted at each outer iteration. This allows use of the low memory variant so that the storage cost of Anderson acceleration is two scalar flux-sized vectors. The outer and inner tolerances $10^{-6}$ and $10^{-8}$, respectively. The uniform subspace correction (USC) preconditioner with one Jacobi iteration and one AMG V-cycle per application is used for the IP and BR2 discretizations. The CG and MDLDG discretizations use one V-cycle of AMG as a preconditioner. The previous outer iteration is used as an initial guess for BiCGStab so that the initial guess becomes progressively better as the outer iteration converges. 
\resp[orange]{For runtime data, each method, refinement, and polynomial order was computed five times with the presented time the minimum runtime across the repeated runs.}

\begin{table}
\centering
\caption{The number of Anderson-accelerated fixed-point iterations until convergence to a tolerance of $10^{-6}$ for the IP, BR2, MDLDG, and CG discretizations of VEF on the linearized crooked pipe problem refined in $h$ and $p$. An Anderson space of size two is used. }
\label{tab:cphp_outer}
\begin{tabular}{ccccccc}
\toprule
 & $N_e$ & IP & BR2 & MDLDG & CG \\
\midrule
\multirow{4}{*}{\rotatebox{90}{$p=1$}} & 112 & 10 & 10 & 13 & 10 \\
 & 448 & 11 & 11 & 14 & 11 \\
 & 1792 & 13 & 13 & 16 & 13 \\
 & 7168 & 14 & 14 & 16 & 14 \\
\addlinespace
\multirow{4}{*}{\rotatebox{90}{$p=2$}} & 112 & 13 & 13 & 15 & 13 \\
 & 448 & 14 & 14 & 16 & 14 \\
 & 1792 & 15 & 15 & 16 & 15 \\
 & 7168 & 15 & 15 & 17 & 15 \\
\addlinespace
\multirow{4}{*}{\rotatebox{90}{$p=3$}} & 112 & 14 & 14 & 16 & 14 \\
 & 448 & 15 & 16 & 16 & 16 \\
 & 1792 & 15 & 15 & 17 & 15 \\
 & 7168 & 15 & 15 & 17 & 16 \\
\bottomrule
\end{tabular}
\end{table}
Table \ref{tab:cphp_outer} shows the number of outer Anderson-accelerated fixed-point iterations until convergence for each of the four VEF methods. 
\resp[red]{The convergence in outer iterations is identical for the IP, BR2, and CG methods aside from a few deviations by one iteration for the case of $p=3$. The MDLDG method took between 1-3 iterations more to converge than the IP, BR2, and CG methods.} \resp[orange][stabexp]{The IP and BR2 methods both have stabilization terms that regularize toward the CG solution whereas the MDLDG method does not. MDLDG's slower convergence indicates that the numerical diffusion induced by using stabilization terms or a continuous solution representation may mildly increase convergence of the outer iteration. Furthermore, the identical convergence rates exhibited by the IP/BR2 and CG methods suggests that the stabilization terms cause the overall algorithm to behave as if a continuous solution representation were used.}

\begin{table}
\centering
\caption{The maximum, minimum, and average number of inner BiCGStab iterations until convergence to an inner tolerance of $10^{-8}$ across all the outer iterations for each of the VEF methods. The previous outer iterate is used as the initial guess for the inner solver so that the number of inner iterations decreases as the outer iteration converges. }
\label{tab:cphp_inner}
\begin{tabular}{cccccccccccccccccc}
\toprule
 &  & \multicolumn{3}{c}{IP}  &  & \multicolumn{3}{c}{BR2}  &  & \multicolumn{3}{c}{MDLDG}  &  & \multicolumn{3}{c}{CG} \\
\cmidrule{3-5}\cmidrule{7-9}\cmidrule{11-13}\cmidrule{15-17}
 & $N_e$ & Max & Min & Avg. & & Max & Min & Avg. & & Max & Min & Avg. & & Max & Min & Avg. \\
\midrule
\multirow{4}{*}{\rotatebox{90}{$p=1$}} & 112 & 15 & 6 & 12.40 & & 14 & 6 & 12.30 & & 10 & 4 & 7.46 & & 7 & 3 & 5.70 \\
 & 448 & 17 & 6 & 12.82 & & 16 & 6 & 12.36 & & 11 & 4 & 8.21 & & 7 & 3 & 5.82 \\
 & 1792 & 17 & 6 & 12.54 & & 17 & 6 & 12.23 & & 11 & 4 & 8.06 & & 8 & 2 & 5.77 \\
 & 7168 & 18 & 6 & 12.79 & & 17 & 6 & 12.21 & & 12 & 4 & 8.12 & & 8 & 2 & 5.50 \\
\addlinespace
\multirow{4}{*}{\rotatebox{90}{$p=2$}} & 112 & 16 & 5 & 11.77 & & 16 & 5 & 11.77 & & 16 & 4 & 8.93 & & 9 & 3 & 7.08 \\
 & 448 & 17 & 7 & 12.57 & & 16 & 5 & 12.57 & & 12 & 5 & 9.19 & & 10 & 3 & 7.00 \\
 & 1792 & 17 & 5 & 12.87 & & 16 & 5 & 12.73 & & 14 & 6 & 10.50 & & 10 & 3 & 7.13 \\
 & 7168 & 17 & 6 & 12.87 & & 18 & 6 & 13.00 & & 14 & 6 & 10.71 & & 10 & 3 & 7.20 \\
\addlinespace
\multirow{4}{*}{\rotatebox{90}{$p=3$}} & 112 & 21 & 6 & 14.71 & & 18 & 7 & 14.00 & & 30 & 7 & 14.44 & & 11 & 4 & 8.57 \\
 & 448 & 22 & 7 & 15.40 & & 21 & 6 & 14.44 & & 17 & 7 & 13.38 & & 14 & 4 & 9.19 \\
 & 1792 & 22 & 9 & 16.33 & & 22 & 9 & 15.93 & & 18 & 8 & 14.35 & & 15 & 5 & 10.00 \\
 & 7168 & 22 & 9 & 16.73 & & 20 & 9 & 16.73 & & 20 & 8 & 14.76 & & 14 & 4 & 10.50 \\
\bottomrule
\end{tabular}
\end{table}
\begin{table}
\centering
\begin{response}[orange]
\caption{The average time spent per outer iteration assembling and solving the discrete VEF system on $hp$ refinements of the crooked pipe problem. Times are presented in milliseconds.}
\label{tab:cphp_inner_times}
\begin{tabular}{cccccccccccc}
\toprule
 &  & \multicolumn{4}{c}{VEF Assembly Time (ms)}  &  & \multicolumn{4}{c}{VEF Solve Time (ms)} \\
\cmidrule{3-6}\cmidrule{8-11}
 & $N_e$ & IP & BR2 & MDLDG & CG & & IP & BR2 & MDLDG & CG \\
\midrule
\multirow{4}{*}{\rotatebox{90}{$p=1$}} & 112 & 13.05 & 14.68 & 13.76 & 13.07 & & 2.15 & 2.15 & 2.71 & 1.62 \\
 & 448 & 49.88 & 56.78 & 52.67 & 49.79 & & 7.90 & 7.87 & 11.15 & 6.00 \\
 & 1792 & 193.29 & 220.92 & 205.48 & 194.45 & & 30.55 & 30.36 & 43.61 & 23.11 \\
 & 7168 & 766.81 & 874.59 & 818.62 & 784.98 & & 124.33 & 121.41 & 174.15 & 93.34 \\
\addlinespace
\multirow{4}{*}{\rotatebox{90}{$p=2$}} & 112 & 24.88 & 30.84 & 31.04 & 24.77 & & 4.40 & 4.44 & 5.96 & 2.81 \\
 & 448 & 95.89 & 120.08 & 120.35 & 96.90 & & 17.38 & 17.56 & 23.58 & 10.77 \\
 & 1792 & 377.29 & 478.98 & 490.49 & 389.00 & & 70.61 & 71.08 & 101.59 & 42.97 \\
 & 7168 & 1504.68 & 1915.12 & 1986.40 & 1566.33 & & 280.57 & 286.26 & 409.00 & 170.49 \\
\addlinespace
\multirow{4}{*}{\rotatebox{90}{$p=3$}} & 112 & 47.08 & 67.40 & 70.94 & 46.73 & & 11.66 & 11.49 & 18.25 & 6.49 \\
 & 448 & 184.06 & 265.64 & 288.93 & 186.42 & & 47.00 & 44.61 & 73.13 & 26.52 \\
 & 1792 & 737.56 & 1083.37 & 1171.01 & 750.89 & & 199.79 & 195.44 & 298.32 & 110.74 \\
 & 7168 & 3064.93 & 4510.38 & 4872.59 & 3098.09 & & 891.92 & 897.10 & 1301.06 & 463.27 \\
\bottomrule
\end{tabular}
\end{response}
\end{table}
The maximum, minimum, and average number of preconditioned BiCGStab iterations to solve the VEF system at each outer iteration are shown in Table \ref{tab:cphp_inner}. \resp[red]{The use of the previous outer iteration's solution as the initial guess for the inner iteration allows BiCGStab to take fewer and fewer iterations as the outer iteration converges. This can be seen by the discrepancy between the maximum and minimum iterations required to converge.} The CG method required the fewest iterations of all the methods, followed by MDLDG, and then IP and BR2. These results show that BiCGStab preconditioned with the USC preconditioner for the IP and BR2 methods and AMG for the MDLDG and CG methods is a scalable solver for the inner iteration in both $h$ and $p$. 

\resp[orange]{The average assembly and solve times are provided in Table \ref{tab:cphp_inner_times}. Here, the costs have been normalized by the number of outer iterations to facilitate their direct comparison. Note that in our implementation, the linear system for the CG method is formed by building the linear system for the IP VEF method (over the space $Y_p$) and then assembling it onto the continuous finite element space $V_p$. That is, the CG assembly cost includes the cost of assembling the IP VEF linear system and an additional step where entries corresponding to shared degrees of freedom in the space $V_p$ are accumulated. A more optimal implementation would not assemble the bilinear forms over $\Gamma_0$ that ultimately cancel when assembled on a continuous finite element space. MDLDG and BR2 were the most expensive to assemble followed by IP and CG. 
Both the BR2 and MDLDG methods have lifting operators which require factorizing the block-diagonal-by-element $W_p$ total interaction mass matrix, an expense that the IP and CG methods avoid.}

\resp[orange]{The CG method has the fewest linear unknowns to solve for and only applies AMG to the smaller continuous finite element operator. Through the USC preconditioner, IP and BR2 also only apply AMG to the continuous operator but the USC preconditioner includes an additional Jacobi iteration on the interfacial unknowns. Further, USC preconditioned BiCGStab applied to the IP and BR2 VEF systems required $\approx\!50\%$ more iterations to converge compared to applying AMG preconditioned BiCGStab to the CG VEF system. While MDLDG typically took fewer iterations to converge than IP or BR2, AMG is applied to a larger system of equations corresponding to the space $Y_p$ instead of $V_p$, increasing the expense of each AMG V-cycle. 
The sparse matrix operations associated with solving the MDLDG system are also more expensive than the IP, BR2, and CG methods due to MDLDG's non-compact stencil which decreases the sparsity of the system. 
Thus, Table \ref{tab:cphp_inner_times} shows MDLDG as the most expensive to solve, followed by IP and BR2, with the CG linear system the fastest to solve.}

\begin{table}
\centering
\begin{response}[orange]
\caption{The total runtime along with the total time spent in the transport sweep and in forming and solving the VEF equations on the crooked pipe problem under $hp$ refinement. All times are presented in seconds.}
\label{tab:cphp_outer_times}
\begin{adjustbox}{max width=\textwidth}
\begin{tabular}{ccccccccccccccccc}
\toprule
 &  & \multicolumn{4}{c}{Total Time (s)}  &  & \multicolumn{4}{c}{Sweep Time (s)}  &  & \multicolumn{4}{c}{VEF Time (s)} \\
\cmidrule{3-6}\cmidrule{8-11}\cmidrule{13-16}
 & $N_e$ & IP & BR2 & MDLDG & CG & & IP & BR2 & MDLDG & CG & & IP & BR2 & MDLDG & CG \\
\midrule
\multirow{4}{*}{\rotatebox{90}{$p=1$}} & 112 & 2.08 & 2.09 & 2.60 & 2.06 & & 1.82 & 1.82 & 2.28 & 1.81 & & 0.15 & 0.17 & 0.22 & 0.15 \\
 & 448 & 7.75 & 7.85 & 9.62 & 7.72 & & 6.72 & 6.74 & 8.35 & 6.72 & & 0.65 & 0.72 & 0.90 & 0.62 \\
 & 1792 & 34.19 & 34.55 & 41.43 & 34.12 & & 29.73 & 29.74 & 35.94 & 29.75 & & 2.98 & 3.33 & 4.05 & 2.88 \\
 & 7168 & 144.94 & 145.97 & 164.66 & 145.43 & & 126.14 & 125.73 & 142.68 & 126.93 & & 12.88 & 14.34 & 16.27 & 12.64 \\
\addlinespace
\multirow{4}{*}{\rotatebox{90}{$p=2$}} & 112 & 4.30 & 4.37 & 4.97 & 4.27 & & 3.76 & 3.75 & 4.26 & 3.75 & & 0.39 & 0.47 & 0.56 & 0.36 \\
 & 448 & 16.77 & 17.05 & 19.39 & 16.70 & & 14.57 & 14.52 & 16.49 & 14.60 & & 1.63 & 1.97 & 2.33 & 1.53 \\
 & 1792 & 70.16 & 71.67 & 77.15 & 69.93 & & 60.69 & 60.88 & 65.20 & 61.03 & & 6.95 & 8.48 & 9.66 & 6.65 \\
 & 7168 & 278.40 & 285.69 & 323.74 & 279.80 & & 241.33 & 242.30 & 272.90 & 243.85 & & 27.91 & 34.18 & 41.67 & 26.91 \\
\addlinespace
\multirow{4}{*}{\rotatebox{90}{$p=3$}} & 112 & 9.21 & 9.50 & 10.86 & 9.13 & & 8.08 & 8.09 & 9.13 & 8.08 & & 0.84 & 1.13 & 1.44 & 0.76 \\
 & 448 & 37.88 & 41.58 & 42.36 & 40.47 & & 33.16 & 35.36 & 35.28 & 35.85 & & 3.57 & 5.08 & 5.86 & 3.47 \\
 & 1792 & 150.60 & 155.33 & 178.69 & 150.10 & & 131.28 & 130.94 & 148.44 & 132.15 & & 14.63 & 19.74 & 25.38 & 13.30 \\
 & 7168 & 626.69 & 653.19 & 736.71 & 659.28 & & 544.86 & 550.04 & 609.53 & 580.13 & & 61.90 & 83.70 & 106.64 & 58.64 \\
\bottomrule
\end{tabular}
\end{adjustbox}
\end{response}
\end{table}
\newcommand{\sweepratioip}{\num{9.5}\xspace} 
\newcommand{\sweepratiobr}{\num{7.9}\xspace} 
\newcommand{\sweepratiomdldg}{\num{7.4}\xspace} 
\newcommand{\sweepratiocg}{\num{10.2}\xspace} 
\newcommand{\rstdmax}{\num{10.2}\xspace} 
\newcommand{\rstdmin}{\num{4.0}\xspace} 

\resp[orange]{Next, we compare the total runtime to find the fixed-point $\varphi = G(\varphi)$ on the crooked pipe problem along with the relative costs of the two major components of evaluating the fixed-point operator $G(\varphi)$: the transport inversion, referred to as the transport sweep, and forming and solving the discrete VEF equations. This timing data is shown in Table \ref{tab:cphp_outer_times}. The sweep and VEF costs are averaged over the number of outer iterations. 
The ratio of the sweep to VEF costs averaged over four refinements in $h$ and three refinements in $p$ was \sweepratioip, \sweepratiobr, \sweepratiomdldg, and \sweepratiocg for the IP, BR2, MDLDG, and CG methods, respectively. 
The relative standard deviation in total runtime across the four methods ranged from 4\% to 10\%. 
In other words, the sweep dominates the cost of the algorithm and thus total runtime was largely insensitive to the choice of VEF discretization. }

\begin{table}
\centering
\begin{response}
\caption{The average percent of elements requiring the application of the negative flux fixup during the transport sweep on the linearized crooked pipe problem.}
\label{tab:cphp_fixup}
\begin{tabular}{ccccccc}
\toprule
 & $N_e$ & IP & BR2 & MDLDG & CG \\
\midrule
\multirow{4}{*}{\rotatebox{90}{$p=1$}} & 112 & 13.936 & 14.332 & 10.315 & 15.194 \\
 & 448 & 5.065 & 5.405 & 4.112 & 6.731 \\
 & 1792 & 2.255 & 2.249 & 2.017 & 2.503 \\
 & 7168 & 1.221 & 1.216 & 1.221 & 1.229 \\
\addlinespace
\multirow{4}{*}{\rotatebox{90}{$p=2$}} & 112 & 14.749 & 16.522 & 12.592 & 16.791 \\
 & 448 & 6.501 & 6.857 & 5.492 & 7.344 \\
 & 1792 & 3.111 & 3.109 & 2.832 & 3.192 \\
 & 7168 & 1.934 & 1.930 & 1.870 & 1.934 \\
\addlinespace
\multirow{4}{*}{\rotatebox{90}{$p=3$}} & 112 & 19.633 & 20.133 & 15.137 & 20.499 \\
 & 448 & 8.149 & 8.431 & 6.874 & 8.652 \\
 & 1792 & 4.555 & 4.556 & 4.283 & 4.565 \\
 & 7168 & 2.674 & 2.677 & 2.583 & 2.696 \\
\bottomrule
\end{tabular}
\end{response}
\end{table}
\begin{response}[orange]
The variance in the average sweep time is due to the methods varying use of the negative flux fixup. Table \ref{tab:cphp_fixup} provides the average percentage of the elements in the transport sweep where the flux fixup was applied. 
Generally, refining the mesh reduced the reliance on the fixup since the solutions converge to the true solution that is positive as $h\rightarrow 0$ while increasing the polynomial order caused an increase in fixup usage due to the increased oscillations caused by high-order interpolation. 
From least to most reliance on the flux fixup, the methods were ordered: MDLDG, IP, BR2, CG.
The IP, BR2, and CG methods had similar usage of the fixup whereas MDLDG required significantly less. For example, on the coarsest meshes MDLDG differed from IP, BR2, and CG by between three and five percentage points. This discrepancy indicates that the more numerically diffusive VEF discretizations create scattering sources that are more likely to induce negativities in the transport sweep. This effect may be caused by an increase in numerical oscillations near material discontinuities produced by the methods that use a continuous solution representation or have a stabilization term that regularizes towards the continuous solution compared to the minimally dissipative MDLDG VEF solution. 

Finally, we note that the IP and CG methods were the fastest in overall runtime. Aside from the case of $p=3$ with one and three refinements, where the CG algorithm required one more outer iteration than the IP method, IP and CG had nearly equivalent run times. While the CG VEF solve was faster than the IP VEF solve, this speedup was balanced by longer sweep times due to CG VEF's increased reliance on the negative flux fixup. The next fastest was BR2 which was slowed down by longer assembly times compared to IP VEF. Finally, MDLDG was the slowest in overall runtime due to both its more expensive assembly and solve times and its larger number of outer iterations required for convergence compared to the other methods.
\end{response}

\begin{response}[blue]
\subsection{Spatial Convergence to a Reference Transport Method} \label{sec:snconv}
In this section, we compare the solutions generated by the IP VEF method and a reference transport method taken to be the high-order DG \Sn method and DSA preconditioner of \citet{ldrd_dsa} as the mesh is refined. 
Convergence between the VEF and \Sn solutions is shown on the thick diffusion limit problem from \S \ref{sec:tdl} and the crooked pipe problem from \S \ref{sec:lcp}. 
Given a fixed angular quadrature rule, let the asymptotic spatial error for the VEF and \Sn methods in isolation be written: 
	\begin{subequations} \label{eq:asym_err}
	\begin{equation}
		\| \varphi_\text{VEF} - \varphi_\text{ex} \| = C_\text{VEF} h^{p+1} \,,
	\end{equation}
	\begin{equation}
		\| \varphi_\text{SN} - \varphi_\text{ex} \| = C_\text{SN} h^{p+1} \,,
	\end{equation}
	\end{subequations}
where $\varphi_\text{ex}$ is the true solution of the problem, $\varphi_\text{VEF}$ and $\varphi_\text{SN}$ the VEF and \Sn numerical solutions, respectively, the $C_i$ the error constants, $h$ the mesh size, and $p$ the finite element polynomial degree. We use the same mesh and polynomial degree for both the VEF and \Sn methods so that the value of $h^{p+1}$ is the same in both methods.  
Comparison of the VEF and \Sn solutions is facilitated by the following bound that makes use of the triangle inequality: 
	\begin{equation}
	\begin{aligned}
		\| \varphi_\text{VEF} - \varphi_\text{SN} \| &= \| (\varphi_\text{VEF} - \varphi_\text{ex}) + (\varphi_\text{ex} - \varphi_\text{SN}) \| \\
		&\leq \| \varphi_\text{VEF} - \varphi_\text{ex} \| + \| \varphi_\text{SN} - \varphi_\text{ex}\| \\
		&= (C_\text{VEF} + C_\text{SN}) h^{p+1} \,. 
	\end{aligned}
	\end{equation}
Thus, we expect that the solutions produced by VEF and \Sn will converge with order $p+1$ on smooth problems. Note, however, that this bound relies on the assumption that the numerical solutions are resolved enough to exhibit the asymptotic error behavior characterized by Eqs.~\ref{eq:asym_err} (i.e.~$\sigma_t h \ll 1)$. In particular, this means that convergence between the solutions produced by the VEF and \Sn methods can only be expected when boundary layers are resolved. 
Non-uniform meshes are used to reduce the expense of capturing the steep gradients present in the boundary layers of the thick diffusion limit and crooked pipe problems. 

\subsubsection{Thick Diffusion Limit}
The single-material thick diffusion limit problem is used to demonstrate convergence between the VEF and \Sn solutions on a problem with a smooth solution. In particular, this problem has no angular discontinuities and, for $\epsilon$ small enough, has a solution that is linearly anisotropic in angle. This means that the error due to \Sn angular quadrature is much smaller than the spatial error, allowing spatial convergence to be seen. We present results for $p=2$. From Eq.~\ref{eq:asym_err}, we expect to see third-order convergence. 
The material parameters are defined in Eq.~\ref{eq:tdl_scaling} with the domain $\D = [0,1]^2$ and $\epsilon = 10^{-2}$. To capture the boundary layer, meshes built from the Chebyshev points defined on the interval $[0,1]$ are used. The Chebyshev points cluster at the endpoints leading to a mesh that grades elements toward the boundary. An example of a mesh built from 21 Chebyshev points is shown in Fig.~\ref{fig:chebmesh}. 
\begin{figure}
\centering
\begin{subfigure}{.49\textwidth}
	\centering
	\includegraphics[height=2.5in]{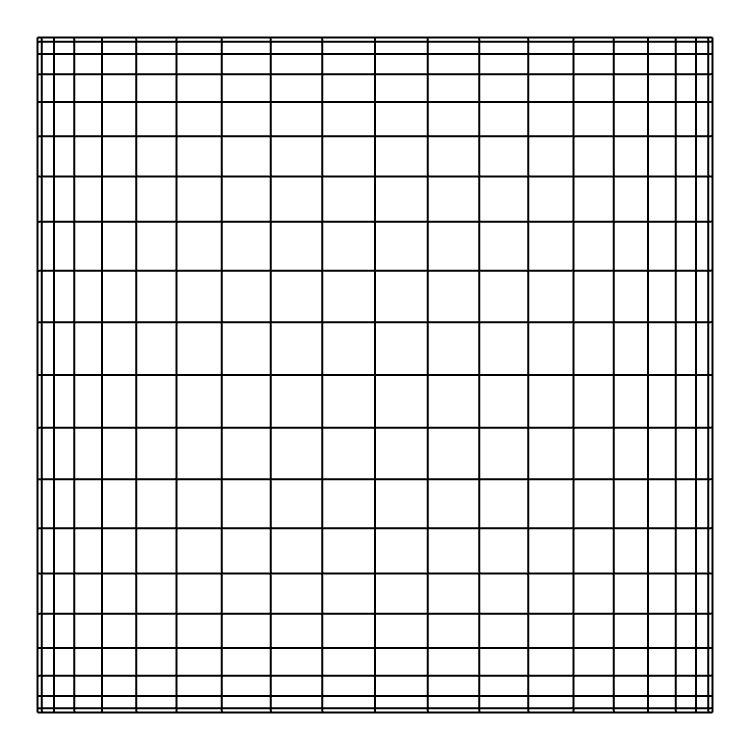}
	\caption{}
	\label{fig:chebmesh}
\end{subfigure}
\begin{subfigure}{.49\textwidth}
	\centering
	\includegraphics[height=2.5in]{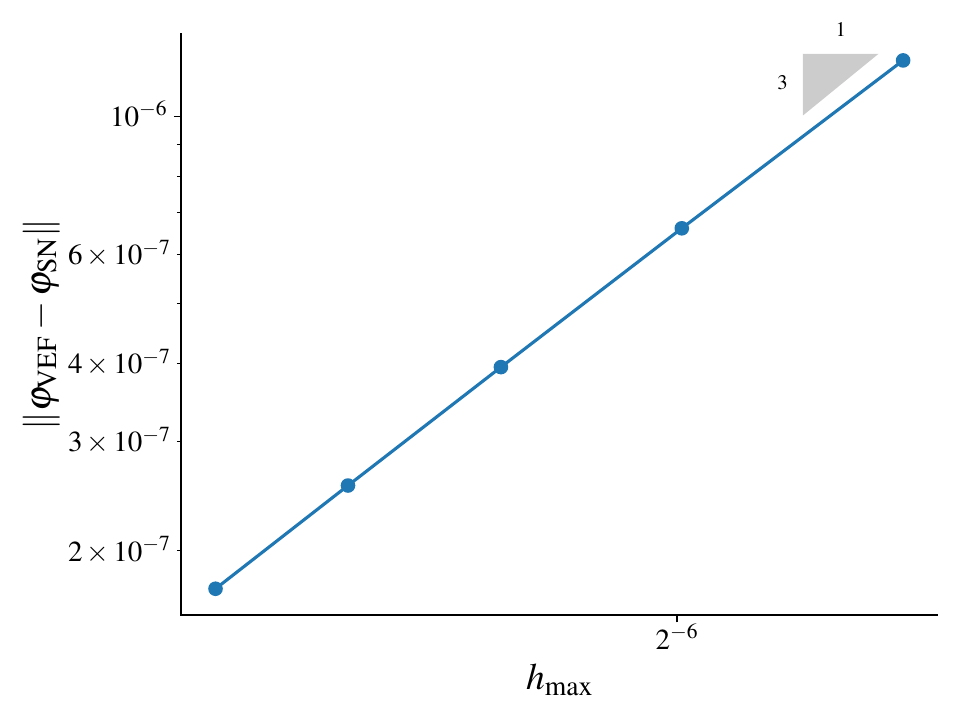}
	\caption{}
	\label{fig:dsacomp_tdl}
\end{subfigure}
\caption{\resp[blue]{(a) An example of a mesh built from a tensor product of Chebyshev points in the interval $[0,1]$ used to resolve the steep gradients in the solution at the boundary of the domain on the thick diffusion limit problem. (b) A plot of the $L^2(\D)$ norm difference between the solutions generated by the IP VEF method and a DG \Sn transport method preconditioned with DSA on the thick diffusion limit problem with $\epsilon = 10^{-2}$. Both the VEF and \Sn methods used $p=2$. The solutions are compared on four meshes generated from a tensor product of 61, 81, 101, and 121 Chebyshev points in each direction. The norms are presented as a function of the maximum characteristic mesh length, $h_\text{max}$. Convergence is compared to a reference third-order line to show that the VEF and \Sn solutions converge at the expected order. }}
\end{figure}

Figure \ref{fig:dsacomp_tdl} shows the convergence of the VEF and \Sn solutions in the $L^2(\D)$ norm as a function of the maximum characteristic element length in the mesh. The solutions are compared on four meshes generated from 61, 81, 101, and 121 Chebyshev points. The ratio of the maximum to minimum mesh size increases as more Chebyshev points are used: for the mesh built from 61 Chebyshev points $h_\text{max}/h_\text{min} \approx 51$ while for 121 Chebyshev points $h_\text{max}/h_\text{min} \approx 89$. The experimentally observed order of convergence on these four meshes was determined to be 2.825 using logarithmic regression. Note that $h_\text{max}$ was used in the logarithmic regression calculation. This result demonstrates that VEF methods do in fact produce the transport solution and can converge with high-order accuracy on a problem that is smooth in space and angle. 

\subsubsection{Crooked Pipe}
\begin{figure}
	\centering 
	\begin{subfigure}{.65\textwidth}
		\centering
		\includegraphics[width=\textwidth]{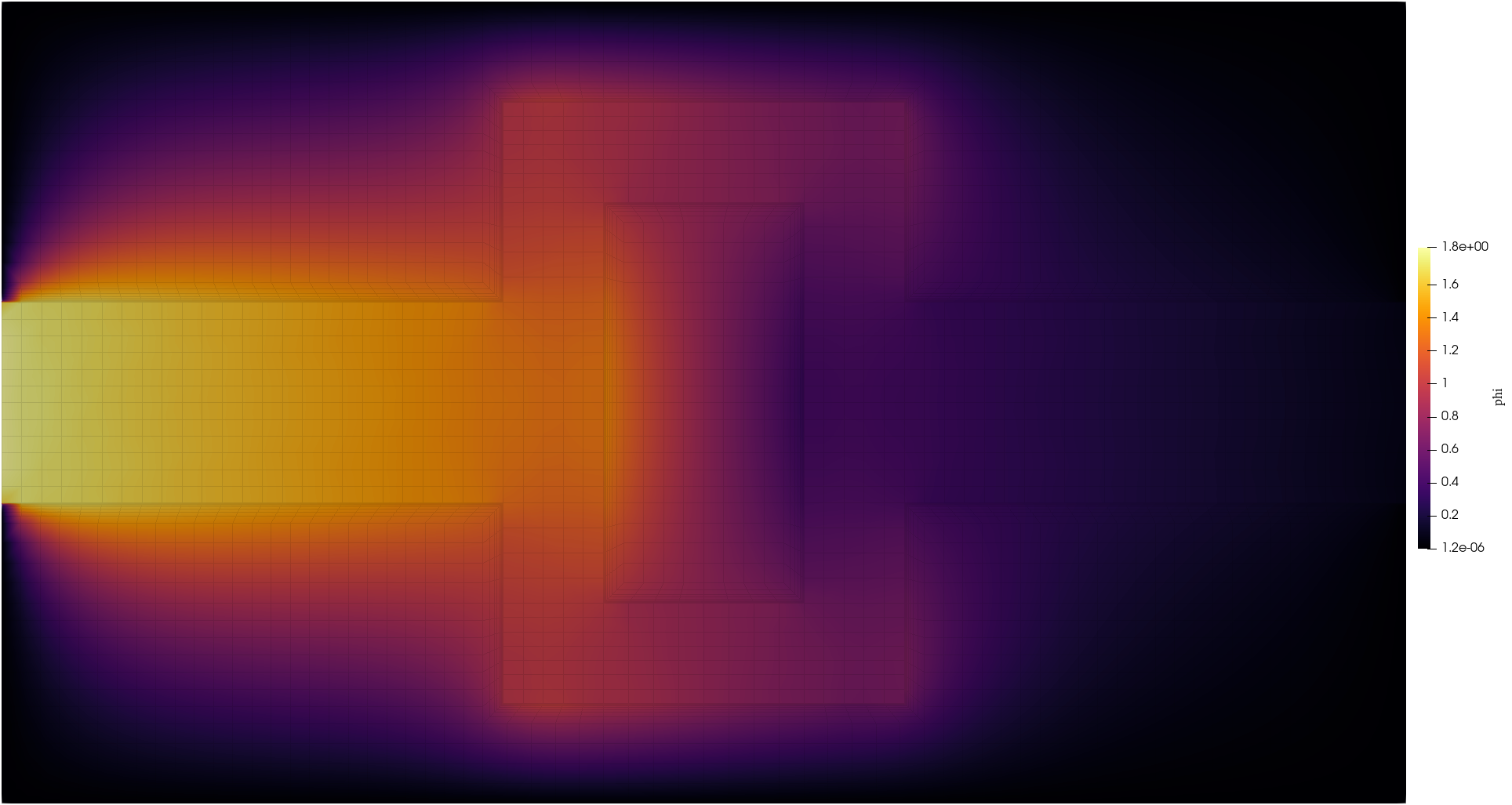}
		\caption{}	
		\label{fig:cpr_vef}	
	\end{subfigure}
	\begin{subfigure}{.65\textwidth}
		\centering
		\includegraphics[width=\textwidth]{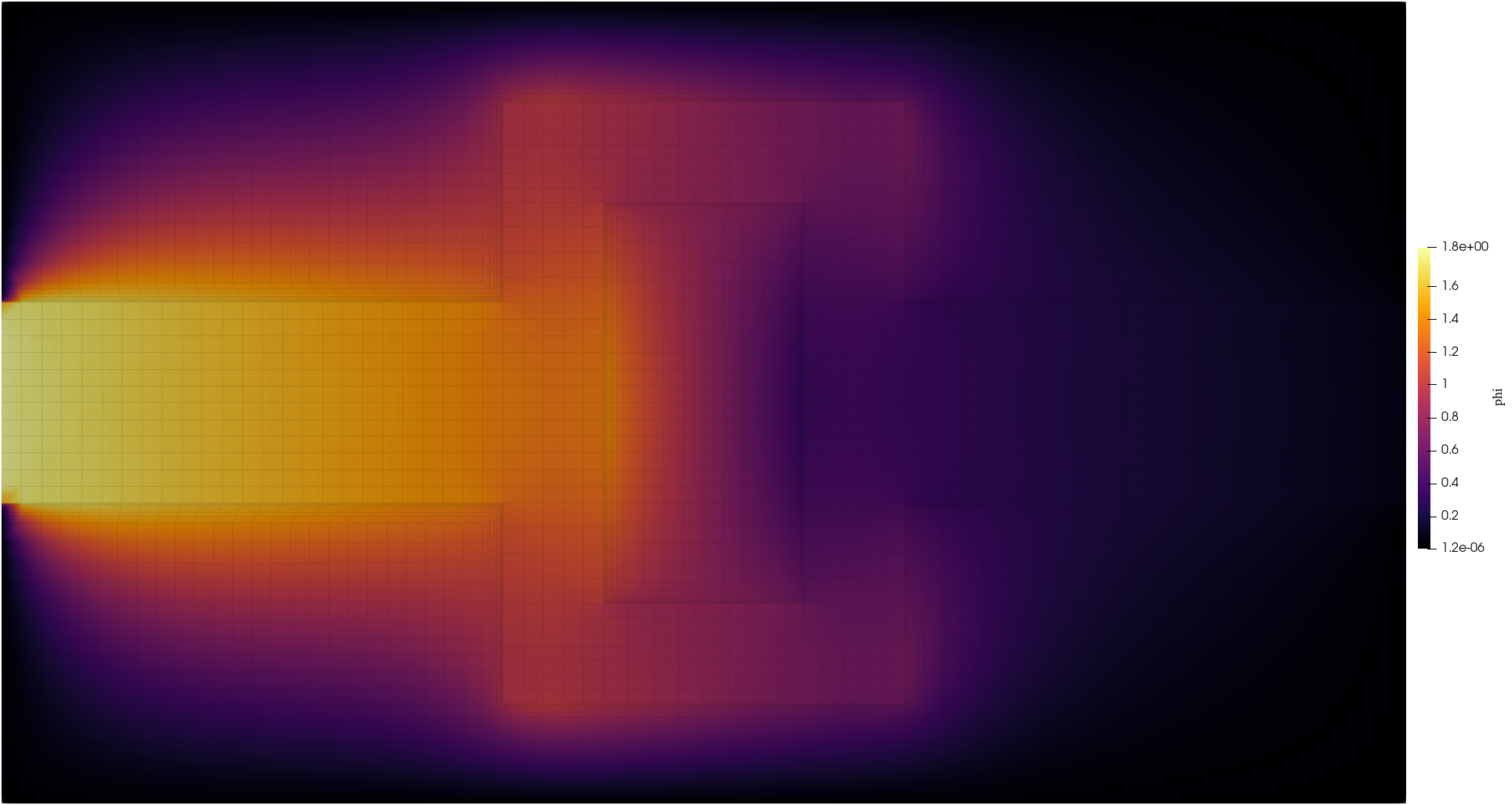}
		\caption{}
		\label{fig:cpr_dsa}
	\end{subfigure}
	\begin{subfigure}{.65\textwidth}
		\centering
		\includegraphics[width=\textwidth]{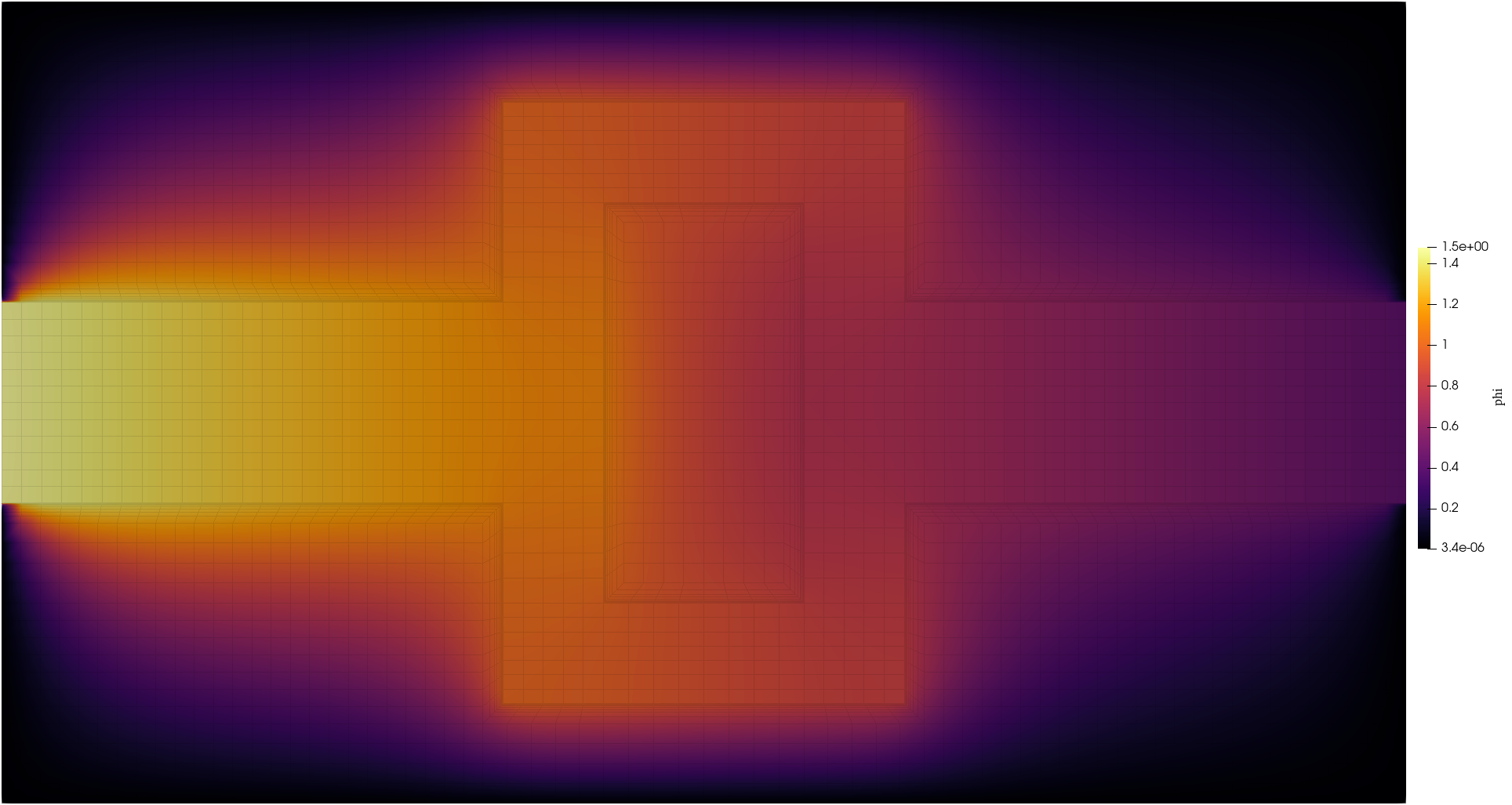}
		\caption{}
		\label{fig:cpr_diff}		
	\end{subfigure}
	\caption{Solutions to the crooked pipe problem using (a) the IP VEF method, (b) the \Sn method from \cite{ldrd_dsa}, and (c) an IP radiation diffusion model. The mesh is refined at the interface between the thick and thin regions. All three methods used $p=2$. The transport models both used $S_{12}$ level symmetric angular quadrature. The VEF and \Sn methods both show the shadow induced by the inner wall forcing the radiation to flow around the pipe that the diffusion model misses. }
	\label{fig:cp_sol}
\end{figure}
\begin{figure}
	\centering
	\begin{subfigure}{.32\textwidth}
		\centering
		\includegraphics[width=\textwidth]{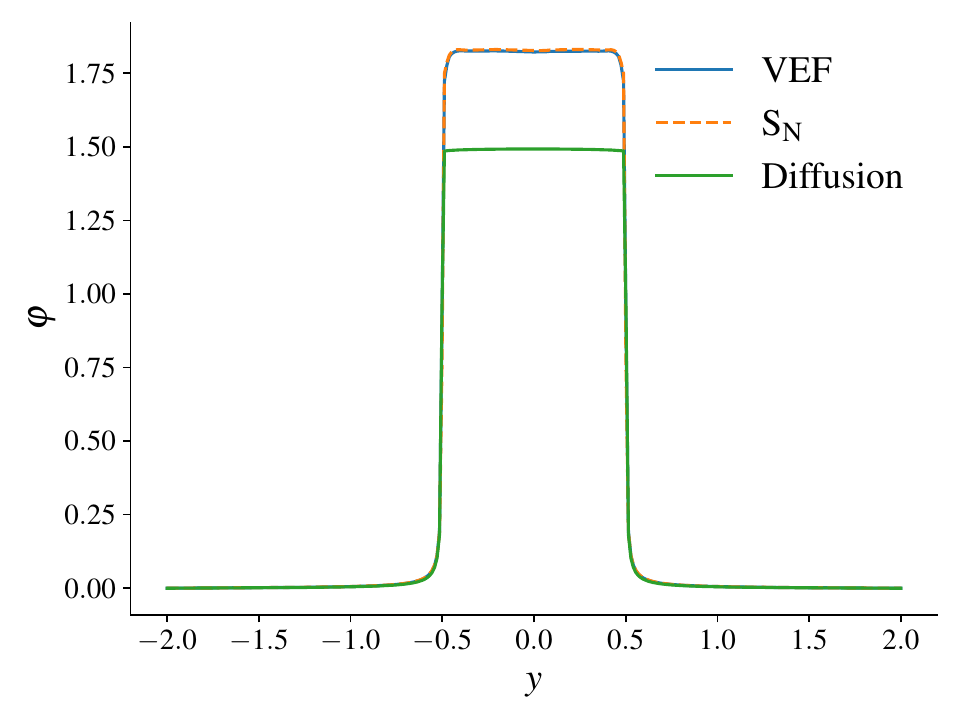}
		\caption{}
	\end{subfigure}
	\begin{subfigure}{.32\textwidth}
		\centering
		\includegraphics[width=\textwidth]{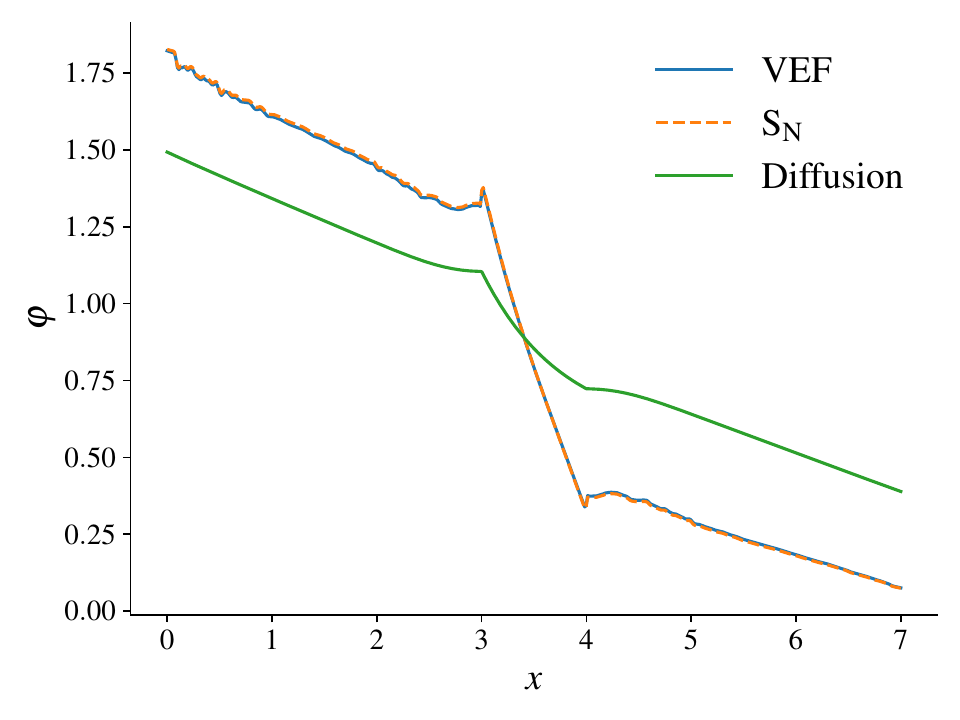}
		\caption{}
		\label{fig:cp_centerline}
	\end{subfigure}
	\begin{subfigure}{.32\textwidth}
		\centering
		\includegraphics[width=\textwidth]{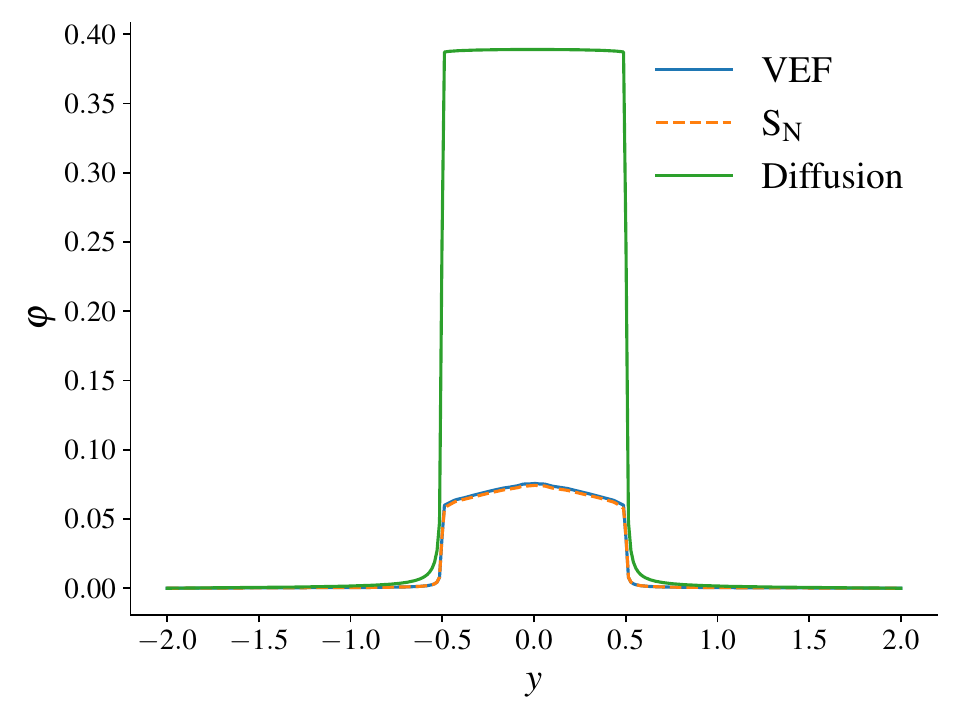}
		\caption{}
	\end{subfigure}
	\caption{\resp{Lineouts of the VEF, \Sn, and diffusion crooked pipe solutions from Fig.~\ref{fig:cp_sol}. The solutions are compared along (a) the vertical line $x=0$, (b) the horizontal line $y=0$, and (c) the vertical line $x=7$. The horizontal lineout shows the diffusion model incorrectly under and over heating the front and back side of the inner wall, respectively, when compared to the VEF and \Sn transport models. The \Sn and VEF solutions are visually close with small discrepancies due to the numerical errors present in both the \Sn and VEF solutions. }} 
	\label{fig:cp_lineouts}
\end{figure}
We now show that VEF converges to the reference transport method on a multi-material problem. The geometry and material data are defined in \S\ref{sec:lcp}. We use three refinements of a non-uniform mesh that grades the elements along the optically thick, wall side of the interface between the two materials. The base mesh contains \num{5216} elements and has minimum and maximum characteristic mesh lengths of \num{8e-3} and \num{1e-1}, respectively. We use $p=2$ and S$_{12}$ angular quadrature in this comparison. This leads to $\approx\!4$ million and $\approx\!250$ million angular flux unknowns on the base mesh and the mesh with three uniform refinements, respectively. 

Figure \ref{fig:cp_sol} shows the scalar flux solutions to the crooked pipe problem on the base mesh with $S_{12}$ angular quadrature computed using the IP VEF method, the \Sn method of \citet{ldrd_dsa}, and an IP radiation diffusion model derived by setting the VEF data to their asymptotic, diffusive values of $\E = \frac{1}{3}\I$ and $E_b = 1/2$. All methods were solved on the same mesh and use $p=2$. 
Compared to the transport models, the diffusion model under and over heats the front and back of the inner wall, respectively, and over predicts the outflow at the end of the pipe. This behavior is evident in the lineouts of the solutions provided in Fig.~\ref{fig:cp_lineouts} which show the solutions along the vertical lines at the left and right edges of the domain and along the center line defined by $y=0$. 

Using DG \Sn as a reference solution, we see that the VEF solution does capture transport effects such as the ``shadow'' behind the inner wall induced by the radiation turning the corner as well as the oscillations in the solution seen in Fig.~\ref{fig:cp_centerline} likely induced by ray effects. On the base mesh using $S_{12}$ angular quadrature, the \Sn and VEF solutions differ in the $L^2(\D)$ norm by 6.29\%. We stress that the discrepancy between \Sn and VEF is due to the numerical errors present in both the VEF solution and the \Sn solution we consider as the reference transport solution. By contrast, the diffusion model differs from \Sn by 101\% in the $L^2(\D)$ norm.

\begin{figure}
\centering
\includegraphics[width=.5\textwidth]{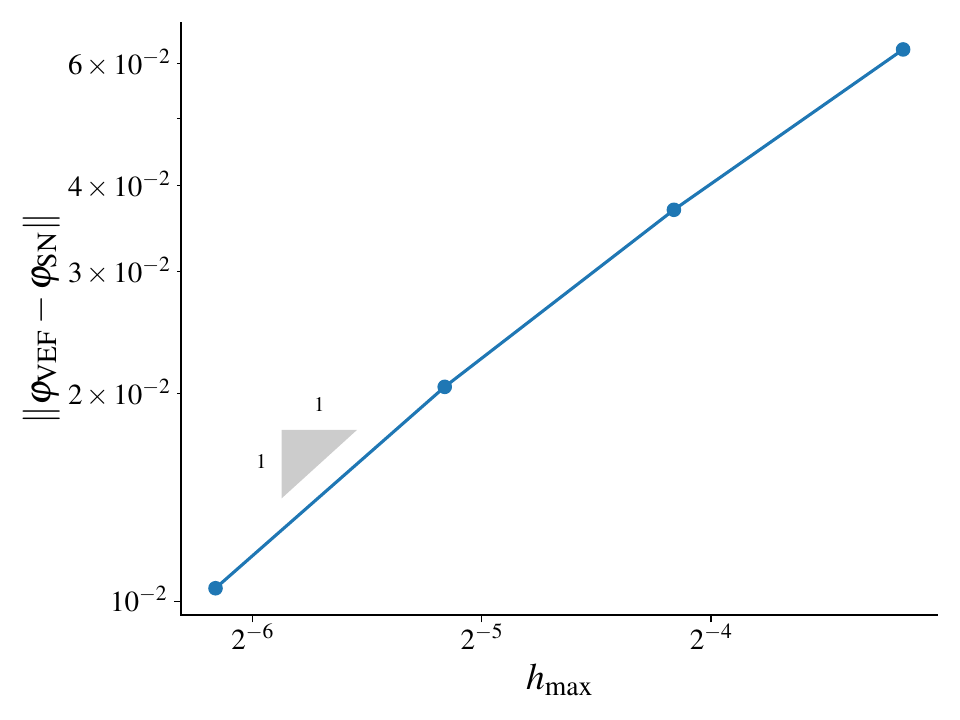}
\caption{\resp[blue]{The difference between the IP VEF and DG \Sn solutions in the $L^2(\D)$ norm on the crooked pipe problem as a function of the maximum characteristic mesh length, $h_\text{max}$. Both methods used $p=2$ and S$_{12}$ angular quadrature. The solutions are solved on three uniform refinements of a base mesh that grades elements along the optically thick side of the interface between the wall and pipe. First-order convergence in space is observed once the problem is spatially resolved enough.}}
\label{fig:dsacomp_cp}
\end{figure}
The $L^2(\D)$ difference between the VEF and \Sn solutions computed on the base mesh and with three uniform refinements is plotted in Fig.~\ref{fig:dsacomp_cp} as a function of the maximum mesh length. Third-order convergence is predicted since we use $p=2$. However, we have observed only first order convergence after the problem is spatially resolved enough. This sub-optimal convergence may be due to the VEF method's use of a negative flux fixup in the sweep where the \Sn method does not use a fixup or under resolution in space and/or angle. We stress that this is a difficult problem where the solution has angular dependence with discontinuous derivatives. In such case, angular quadrature is expected to converge slowly as the number of angles is increased. Furthermore, the solutions shown in Figs.~\ref{fig:cpr_vef} and \ref{fig:cpr_dsa} have visually obvious ray effects indicating under resolution in angle. Thus, it is unclear whether high-order convergence in space between these two numerically disparate schemes is possible on this problem. 
\end{response}

\subsection{Weak Scaling}
Here, we present a weak scaling study of the IP VEF method with $p=2$ for both the inner iteration in isolation and the full fixed-point solve. \resp{Uniform refinements are used in combination with increasing the parallel partitioning by a factor of four so that the degrees of freedom per processor remains constant.} The following results were generated on 29 nodes of the \texttt{rztopaz} machine at LLNL which has two 18-core Intel Xeon E5-2695 CPUs and 128GB of memory per node. 
\resp[orange]{Timing data is presented as the minimum time measured across three repeated runs.}

\subsubsection{Inner Solve on Problem with Mock VEF Data} \label{sec:weak_mock}
First, we investigate weak scaling on a mock problem where the VEF data are provided as inputs to the problem (as opposed to being solved for through fixed-point iteration). This allows the VEF system to be solved in isolation from the transport equation. We use the materials, geometry, and boundary conditions from the crooked pipe problem shown in Fig.~\ref{fig:cp_diag} but set the Eddington tensor and boundary factor to 
	\begin{subequations}
	\begin{equation}
		\E = \begin{cases}
			\begin{bmatrix} 
				9/11 & 0 \\ 0 & 1/11 
			\end{bmatrix}\,, & \x \in \mathrm{pipe} \\[1em] 
			\begin{bmatrix} 
				1/3	& 0 \\ 0 & 1/3
			\end{bmatrix}\,, & \x \in \mathrm{wall} 
		\end{cases} \,, 
	\end{equation}
	\begin{equation}
		E_b = \begin{cases}
			9/10 \,, & \x \in \partial(\mathrm{pipe})\\
			1/2 \,, & \x \in \partial(\mathrm{wall}) 
		\end{cases} \,. 
	\end{equation}
	\end{subequations}
This corresponds to a linearly anisotropic (i.e.~diffusive) angular flux in the wall and an extremely forward peaked solution 
	\begin{equation}
		\psi = \Omegahat_x^8 
	\end{equation}
in the pipe. The motivation for this choice is that the solvers are predicted to struggle when the Eddington tensor is discontinuous. We stress that this setup does not correspond to a physically realistic problem. 
\resp[red]{In particular, the angular flux has an $\mathcal{O}(1)$ jump along the interface between the thick wall and thin pipe.}

\begin{table}
\centering
\caption{Weak scaling the linear solve for the IP VEF method with $p=2$ on a non-physically difficult problem with mock VEF data. The columns parameterize the method used to approximately invert the continuous operator in the USC preconditioner. The standard USC preconditioner with AMG on the continuous operator did not converge due to the large discontinuity in the VEF data. Convergence is recovered when a sparse direct solver is used to invert the continuous operator indicating that AMG is struggling to accurately invert the continuous operator. Efficient iterative solvers are found by applying AMG to a symmetrized continuous operator. The USC-S method applies one AMG V-cycle to this symmetrized operator while USC-S3 uses three AMG V-cycles in an attempt to better approximate the inverse of the original non-symmetric continuous operator. Solve times are also provided. }
\label{tab:weak_mock2}
\begin{tabular}{cccccccccccc}
\toprule
 &  & \multicolumn{4}{c}{Iterations}  &  & \multicolumn{4}{c}{Solve Time (s)} \\
\cmidrule{3-6}\cmidrule{8-11}
Processors & DOF & USC & USC-Direct & USC-S & USC-S3 & & USC & USC-Direct & USC-S & USC-S3 \\
\midrule
1 & \num{12348} & 87 & 21 & 28 & 25 & & 0.24 & 0.25 & 0.08 & 0.15 \\
4 & \num{49392} & 205 & 26 & 30 & 23 & & 0.72 & 1.08 & 0.11 & 0.19 \\
16 & \num{197568} & -- & 24 & 30 & 26 & & -- & 4.38 & 0.16 & 0.31 \\
64 & \num{790272} & -- & 25 & 34 & 24 & & -- & 16.30 & 0.30 & 0.39 \\
256 & \num{3161088} & -- & 27 & 33 & 26 & & -- & 63.20 & 0.32 & 0.48 \\
1024 & \num{12644352} & -- & 28 & 33 & 25 & & -- & 272.24 & 0.37 & 0.52 \\
\bottomrule
\end{tabular}
\raggedright
-- indicates solve did not converge within 250 iterations 
\end{table}
Table \ref{tab:weak_mock2} shows the number of BiCGStab iterations to convergence for the USC-preconditioned IP VEF system on this mock problem. The columns of the table parameterize the solver used for the continuous stage of the USC preconditioner. The standard USC preconditioner used in the previous results did not converge. However, when a sparse direct solver is applied to the CG operator instead of AMG, uniform convergence is recovered. 
\ifreview\linelabel{ln:remcallback}\fi
This suggests that the preconditioner is failing due to AMG's inability to adequately approximate the inverse of the continuous operator when jumps in the Eddington tensor are present (see Remark \ref{rem:amg}). 

Note that AMG is effective on the standard continuous finite element discretization of diffusion. It is then possible that AMG applied to a symmetrized VEF operator could be an effective preconditioner for the non-symmetric VEF bilinear form. 
A symmetric operator more amenable to accurate inversion via AMG is found by lagging the terms
	\begin{equation}
		-\int_{\Gamma_0} \avg{\frac{\nabla u}{\sigma_t}}\cdot \jump{\E\varphi\hat{n}} \ud s + \int \nabla u \cdot \frac{1}{\sigma_t} \paren{\nablah\cdot\E}\!\varphi \ud \x
	\end{equation}
in the CG VEF discretization (Eq.~\ref{eq:cg}). The symmetrized operator is then:
	\begin{equation} \label{eq:sym_bilin}
		\int_{\Gamma_b} E_b\, u \varphi \ud s + \int \nablah u \cdot \frac{1}{\sigma_t}\E \nablah \varphi \ud \x + \int \sigma_a\, u \varphi \ud \x \,, 
	\end{equation}
with $u,v \in V_p$. 
This is a CG discretization of 
	\begin{equation}
		-\nabla \cdot \frac{1}{\sigma_t}\E\nabla \varphi + \sigma_a \varphi 
	\end{equation} 
which corresponds to the VEF drift-diffusion equation where the advective term $\paren{\nabla\cdot\E}\!\varphi$ is lagged and moved to the right hand side. 

The remaining columns of Table \ref{tab:weak_mock2} present the use of AMG applied to the symmetrized VEF operator in Eq.~\ref{eq:sym_bilin} to precondition the original non-symmetric IP VEF system on the mock problem. 
The ``USC-S'' column shows convergence for a preconditioner where one AMG V-cycle is applied to the symmetrized CG operator in place of AMG applied to the non-symmetric operator. 
\resp[blue][symcl]{In other words, an approximate inverse of the symmetric operator given in Eq.~\ref{eq:sym_bilin} is used to precondition and solve the non-symmetric VEF system.} The method converges and is roughly uniform in iteration counts as the mesh is refined. 

The ``USC-S3'' column corresponds to the use of a preconditioner that uses three iterations of an inner Richardson iteration to approximate the inverse of the non-symmetric CG operator. The Richardson iteration is preconditioned using one V-cycle of AMG applied to the symmetrized CG operator. In this way, an approximation to the inverse of the original non-symmetric operator is computed while supporting the use of AMG on the symmetrized operator. \resp[blue][symcltwo]{Note that we use a fixed number of iterations and thus do not attempt to converge the inner iteration. The intent of the inner iteration is only to provide a preconditioner that more closely approximates the inverse of the non-symmetric operator compared to the USC-S option.}
For this preconditioner, iterative efficiency generally fell between that of the sparse direct solver and using only AMG on the symmetrized CG operator. This is expected because more computational work is performed at each iteration compared to the USC-S preconditioner. Inner iterations do reduce the number of total iterations to convergence but, since three V-cycles are performed per preconditioner application, not to the degree that fewer V-cycles are performed. \resp[orange]{This is corroborated by the solve times also presented in Table \ref{tab:weak_mock2}: on the largest problem size, the USC-S3 method was 40\% more expensive despite requiring 8 fewer iterations than USC-S. }

\begin{figure}
\centering
\includegraphics[width=.55\textwidth]{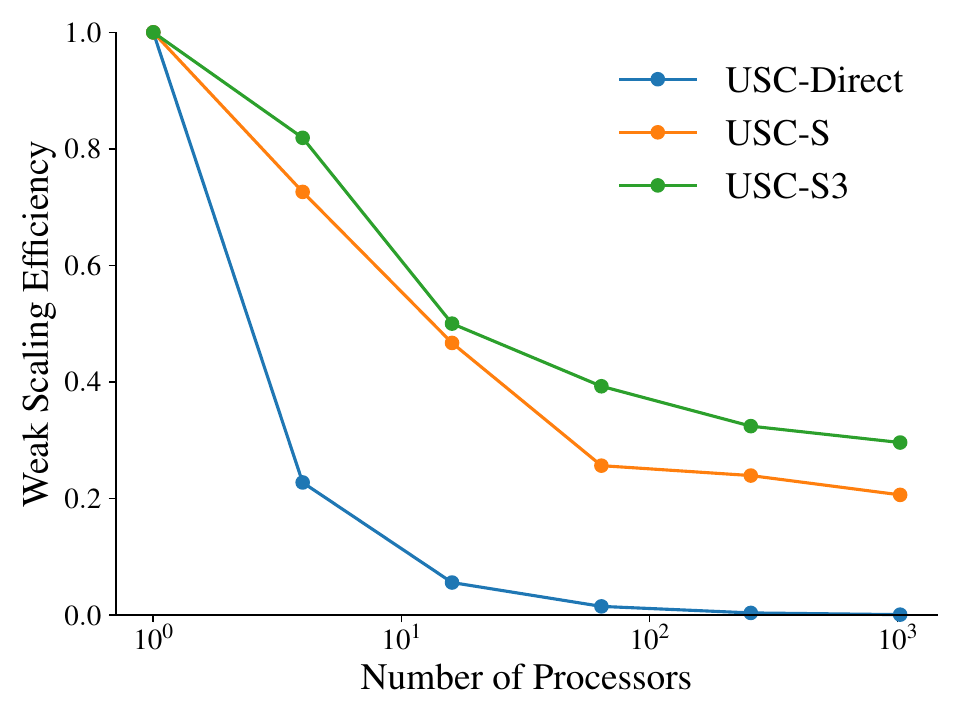}
\caption{Weak scaling efficiency for solving the IP VEF linear system on the mock problem with $p=2$. Three methods for approximating the inverse of the CG operator used in the USC preconditioner are compared. On the mock problem, applying AMG to the CG operator was not convergent. Scaling was recovered by applying AMG to a symmetrized CG operator. The USC-S3 option applies AMG three times per iteration leading to more expensive solve times but overall better scaling compared to USC-S which applies AMG to the symmetrized CG operator only once per iteration. The scaling of the direct method is provided to show the efficiency of a method that is uniform in iteration count but unscalable due to the poor scaling inherent to sparse direct methods. }
\label{fig:mock_weak_efficiency}
\end{figure}
\begin{response}[orange]
Figure \ref{fig:mock_weak_efficiency} plots the weak scaling efficiency defined as
	\begin{equation} \label{eq:scaling_eff}
		\varepsilon_n = \frac{\text{solve time with one processor}}{\text{solve time with $n$ processors}}
	\end{equation}
for the three convergent preconditioning schemes. The ideal weak scaling is $\varepsilon_n = 1$. 
As expected, the sparse direct solver does not weak scale. 
Weak scaling efficiency for the iterative techniques is not expected to be ideal due to the unavoidable communication costs inherent to distributed sparse matrix operations. In particular, weak scaling efficiency is expected to degrade when intra-node communication is required. For the \texttt{rztopaz} architecture, each node has 36 processors meaning intra-node communication is required for the problems where 64, 256, and 1024 processors are used. The weak scaling efficiency of the USC-S and USC-S3 methods appears to saturate at 20\% and 30\%, respectively. The increased efficiency of USC-S3 is due to its more consistent required number of iterations to convergence. Thus, while USC-S3 is more expensive than USC-S it may scale more predictably and robustly.
\end{response}

\subsubsection{Inner Solve on First Iteration of Crooked Pipe with Parallel Block Jacobi Transport Sweep} \label{sec:weak_first_it}
Next, we show weak scaling of the IP VEF linear solve on the first outer iteration of the linearized crooked pipe problem from \S \ref{sec:lcp} with $p=2$. One parallel block Jacobi transport sweep is performed to provide angular fluxes to compute the VEF data. 
In other words, each processor performs a local sweep on its processor-local domain using incoming angular flux information that is lagged and taken from the previous iteration. Thus, each angle can be computed independently on each processor. However, the transport sweep is no longer exact and the convergence of the outer fixed-point problem will now depend on the parallel decomposition. Our use of parallel block Jacobi is motivated by a desire to avoid the communication costs and idle times associated with a full parallel upwind sweep and that problems of interest typically have large regions of optically thick materials that allow parallel block Jacobi to converge quickly. 

\begin{figure}
\centering
\includegraphics[width=.55\textwidth]{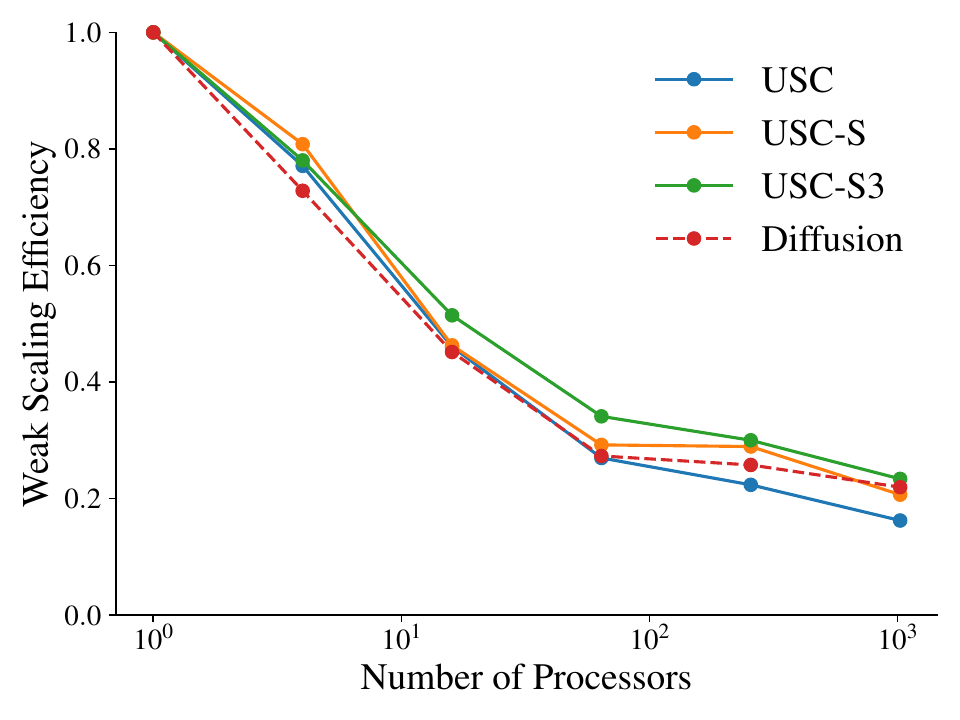}
\caption{Weak scaling efficiency on the first iteration of the crooked pipe problem. Three preconditioners for the IP VEF system are compared. Scaling on the non-symmetric VEF system is compared to solving an IP radiation diffusion system preconditioned with the standard USC preconditioner. Here, it can be seen that all three preconditioners for the VEF system scale similarly to solving the symmetric radiation diffusion system.}
\label{fig:real_weak_efficiency}
\end{figure}

\begin{table}
\centering
\caption{Weak scaling data for the linear solve for the IP VEF method with $p=2$ on the first iteration of the linearized crooked pipe problem. A parallel block Jacobi sweep is used to generate the VEF data needed to form the VEF system. On this physically realistic problem, both the standard USC and the two USC preconditioners that utilize a symmetrized CG operator converged uniformly. Iterative efficiency and solve times are compared to solving the symmetric positive definite linear system associated with an IP discretization of radiation diffusion using the USC preconditioner.}
\label{tab:weak_cp2}
\begin{tabular}{cccccccccccc}
\toprule
 &  & \multicolumn{4}{c}{Iterations}  &  & \multicolumn{4}{c}{Solve Time (s)} \\
\cmidrule{3-6}\cmidrule{8-11}
Processors & DOF & USC & USC-S & USC-S3 & Diffusion & & USC & USC-S & USC-S3 & Diffusion \\
\midrule
1 & \num{42588} & 17 & 19 & 13 & 14 & & 0.16 & 0.18 & 0.27 & 0.14 \\
4 & \num{170352} & 19 & 20 & 14 & 16 & & 0.21 & 0.22 & 0.34 & 0.19 \\
16 & \num{681408} & 20 & 22 & 15 & 17 & & 0.36 & 0.39 & 0.52 & 0.31 \\
64 & \num{2725632} & 21 & 22 & 14 & 17 & & 0.61 & 0.62 & 0.78 & 0.50 \\
256 & \num{10902528} & 25 & 21 & 16 & 18 & & 0.73 & 0.63 & 0.89 & 0.53 \\
1024 & \num{43610112} & 28 & 24 & 15 & 17 & & 1.01 & 0.88 & 1.15 & 0.63 \\
\bottomrule
\end{tabular}
\end{table}
Table \ref{tab:weak_cp2} shows the number of BiCGStab iterations to solve the IP VEF system to a tolerance of $10^{-8}$. BiCGStab is preconditioned with the USC, USC-S, and USC-S3 preconditioners.  
\ifreview\linelabel{ln:diffdef}\fi
In addition, the number of iterations to solve the corresponding IP diffusion problem (by setting $\E = \frac{1}{3}\I$ and $E_b = 1/2$) with the standard USC preconditioner are shown. 
These results indicate that on a physically realistic problem the USC, USC-S, and USC-S3 options are all effective. 
\resp[red]{Compared to IP diffusion, the standard USC preconditioner took 65\% more iterations to solve the non-symmetric VEF equations on the largest problem size. However, when USC-S was used, this discrepancy was reduced to 41\%.}
\resp[orange]{The USC-S method led to the fastest solve times, followed by USC, and then USC-S3. In particular, solving the IP VEF linear system using the USC-S preconditioner was only 12\% more expensive than solving the symmetric IP radiation diffusion system on the finest mesh with 1024 processors.
Weak scaling efficiency is plotted in Fig.~\ref{fig:real_weak_efficiency}. Efficiency of the USC-S and USC-S3 methods is comparable to that of solving the IP radiation diffusion system.}


\begin{response}[red]
\subsubsection{Full Algorithm on Crooked Pipe with Parallel Block Jacobi Transport Sweep}
The single iteration test above is now repeated on the full algorithm. Performance is presented for the IP VEF method with $p=2$ coupled to a parallel block Jacobi sweep. Due to the lagging of incoming angular fluxes on processor boundaries, it is not expected that the number of outer iterations will be independent of the parallel decomposition. The fixed-point tolerance is $10^{-6}$ with the inner BiCGStab tolerance $10^{-8}$. Fixed-point iteration without Anderson acceleration is used. Performance is compared when the inner BiCGStab iteration is preconditioned by the USC and USC-S methods. 
Note that the number of degrees of freedom per processor is significantly lower in this section than the scaling study for the first iteration of the crooked pipe in order to enable solving the full problem at scale. 

\begin{table}
\centering
\begin{response}
\caption{A weak scaling study of the full IP VEF fixed-point solve with $p=2$ on the crooked pipe problem. A parallel block Jacobi transport sweep was used to approximate the inverse of the streaming and collision operator with minimal communication cost. Outer refers to the number of fixed-point iterations required to converge to a tolerance of $10^{-6}$. The maximum, minimum, and average number of inner preconditioned BiCGStab iterations are shown for two types of preconditioners: the standard USC preconditioner which applies AMG to the continuous operator and the USC-S preconditioner which applies AMG to a symmetrized continuous operator. The inner tolerance was $10^{-8}$. The maximum and minimum percentage of elements requiring the negative flux fixup within the parallel block Jacobi transport sweep are presented along with the maximum and minimum values of $\| \jump{\E\hat{n}}\|_{L^2(\Gamma_0)}$. Due to the parallel block Jacobi sweep and small problem size per processor, the outer iteration does not converge independent of the processor count. In addition, the maximum number of USC-preconditioned BiCGStab iterations increases with the processor count due to the increasing maximum discontinuity in the Eddington tensor caused by inaccurate angular fluxes computed by the parallel block Jacobi sweep in the early stages of the outer iteration. This scaling is less pronounced when the USC-S preconditioner is used. Note that both preconditioners have an average number of iterations that is uniform with respect to the processor count and problem size. }
\label{tab:fweak}
\begin{adjustbox}{max width=\textwidth}
\begin{tabular}{ccccccccccccccccc}
\toprule
 &  &  & \multicolumn{3}{c}{Inner It.~(USC)}  &  & \multicolumn{3}{c}{Inner It.~(USC-S)}  &  & \multicolumn{2}{c}{\% Fixed Up}  &  & \multicolumn{2}{c}{Eddington Jump} \\
\cmidrule{4-6}\cmidrule{8-10}\cmidrule{12-13}\cmidrule{15-16}
Processors & DOF & Outer & Max & Min & Avg. & & Max & Min & Avg. & & Max & Min & & Max & Min \\
\midrule
1 & \num{4032} & 30 & 18 & 2 & 9.267 & & 17 & 2 & 9.167 & & 7.39 & 1.61 & & \num{6.5e-02} & \num{4.9e-02} \\
4 & \num{16128} & 51 & 22 & 2 & 10.451 & & 21 & 2 & 10.549 & & 11.08 & 1.14 & & \num{1.7e-01} & \num{4.3e-02} \\
16 & \num{64512} & 76 & 30 & 2 & 10.724 & & 25 & 2 & 10.618 & & 13.97 & 0.74 & & \num{2.0e-01} & \num{2.5e-02} \\
64 & \num{258048} & 135 & 35 & 1 & 10.881 & & 27 & 1 & 10.963 & & 12.93 & 0.39 & & \num{1.4e-01} & \num{1.4e-02} \\
256 & \num{1032192} & 261 & 50 & 1 & 10.590 & & 27 & 1 & 10.508 & & 14.84 & 0.22 & & \num{1.8e-01} & \num{6.8e-03} \\
1024 & \num{4128768} & 540 & 95 & 1 & 10.528 & & 36 & 1 & 10.155 & & 16.22 & 0.13 & & \num{2.5e-01} & \num{2.4e-03} \\
\bottomrule
\end{tabular}
\end{adjustbox}
\end{response}
\end{table}
Table \ref{tab:fweak} presents the outer and inner iteration data as the mesh is refined and the processor count increased. Due to the parallel block Jacobi sweep, the outer iteration is not robust to the processor count. The maximum required inner BiCGStab iterations increases with higher processor counts for both preconditioning schemes. However, the average number of inner iterations per outer iteration remains constant. Table \ref{tab:fweak} also includes information on the percentage of elements that required the negative flux fixup and a measure of the discontinuity of the Eddington tensor. This measure is computed with the $L^2(\Gamma_0)$ norm of the jump of the Eddington tensor applied to the normal. In other words, we tabulate the maximum and minimum values of 
	\begin{equation}
		\| \jump{\E\hat{n}} \|_{L^2(\Gamma_0)} = \sqrt{\int_{\Gamma_0} \jump{\E\hat{n}} \cdot \jump{\E\hat{n}} \ud s} \,, 
	\end{equation}
across each of the outer iterations. Note that the maximum percentage of elements requiring the fixup and the maximum measure of the discontinuity of the Eddington tensor both increase as the processor count is increased and the mesh is refined. Conversely, the minimum values both decrease with processors and mesh refinements. This behavior is likely due to the use of the parallel block Jacobi sweep. At the beginning of the iteration, parallel block Jacobi provides a poor approximation to the inverse of the streaming and collision operator, especially at parallel boundaries, inducing jumps beyond the spatial discretization error in the angular flux and thus the Eddington tensor. 
As the outer iteration converges, parallel block Jacobi provides a better approximation to the inverse of the streaming and collision operator leading to jumps in the angular flux on the order of the discretization error. Thus, the first few iterations will have discontinuities in the Eddington tensor that increase as more processors are used while the last few iterations will have discontinuities that decrease as the discretization error is reduced through mesh refinements. 
We note that the solve requiring the most inner iterations does not occur at the first iteration and thus the maximum inner iterations to convergence deviates from the behavior observed in \S\ref{sec:weak_first_it} which only considers the first outer iteration of the problem under consideration in this section. 

As observed in \S\ref{sec:weak_mock}, the standard USC method degrades when strong jumps in the Eddington tensor are present. This explains the scaling of the maximum number of inner iterations with processors for the standard USC preconditioner: the Eddington tensor in the first few outer iterations becomes more discontinuous along processor boundaries leading to higher and higher maximum inner iteration counts. However, these discontinuities decrease in magnitude as the iteration converges leading to the uniform scaling in minimum and average inner iterations to convergence. The USC-S preconditioner was more robust to these discontinuities, having either equivalent or much better convergence as the processor count increased. 
In particular, for the 1024 processor case, USC-S converged in 37\% of the iterations that USC did.

\begin{table}
\centering
\begin{response}[red]
\caption{Timing data for the VEF components of the full algorithm weak scaling study. The USC and USC-S preconditioners for the inner BiCGStab iteration are compared. All times are presented in milliseconds as the average cost per outer iteration. The USC-S preconditioner is cheaper to solve (for high processor counts) but is more expensive overall compared to the USC preconditioner due to its additional cost of assembling the symmetrized continuous operator. }
\label{tab:fweak_timing}
\begin{tabular}{cccccccccccccc}
\toprule
 &  & \multicolumn{2}{c}{Total VEF}  &  & \multicolumn{2}{c}{VEF Assembly}  &  & \multicolumn{2}{c}{Prec.~Setup}  &  & \multicolumn{2}{c}{Solve} \\
\cmidrule{3-4}\cmidrule{6-7}\cmidrule{9-10}\cmidrule{12-13}
Processors & DOF & USC & USC-S & & USC & USC-S & & USC & USC-S & & USC & USC-S \\
\midrule
1 & \num{4032} & 105.20 & 123.16 & & 88.35 & 89.50 & & 1.10 & 17.70 & & 14.01 & 14.21 \\
4 & \num{16128} & 122.46 & 140.02 & & 96.73 & 97.98 & & 1.94 & 18.15 & & 21.88 & 21.97 \\
16 & \num{64512} & 156.68 & 177.42 & & 116.66 & 117.57 & & 3.08 & 23.09 & & 34.36 & 34.26 \\
64 & \num{258048} & 204.02 & 228.52 & & 151.40 & 151.32 & & 3.75 & 29.16 & & 46.58 & 45.72 \\
256 & \num{1032192} & 212.15 & 237.02 & & 151.20 & 151.58 & & 3.47 & 29.28 & & 54.56 & 53.14 \\
1024 & \num{4128768} & 232.67 & 250.01 & & 152.29 & 152.76 & & 4.06 & 30.56 & & 72.71 & 63.77 \\
\bottomrule
\end{tabular}
\end{response}
\end{table}
Table \ref{tab:fweak_timing} shows timing data for the inner iteration preconditioned by the two USC methods. The data are averaged across all outer iterations. For the USC preconditioner, the preconditioner setup cost includes linear algebra operations that build the CG operator without rediscretizing and the setup costs associated with AMG. For USC-S, this cost additionally includes forming the symmetrized operator. Due to the symmetrization of the operator, it must be assembled independently from the original IP VEF operator, incurring additional assembly costs in the setup phase. Thus, the standard USC preconditioner was faster despite the USC-S preconditioner requiring fewer total iterations to converge. 

\begin{figure}
\centering
\includegraphics[width=.55\textwidth]{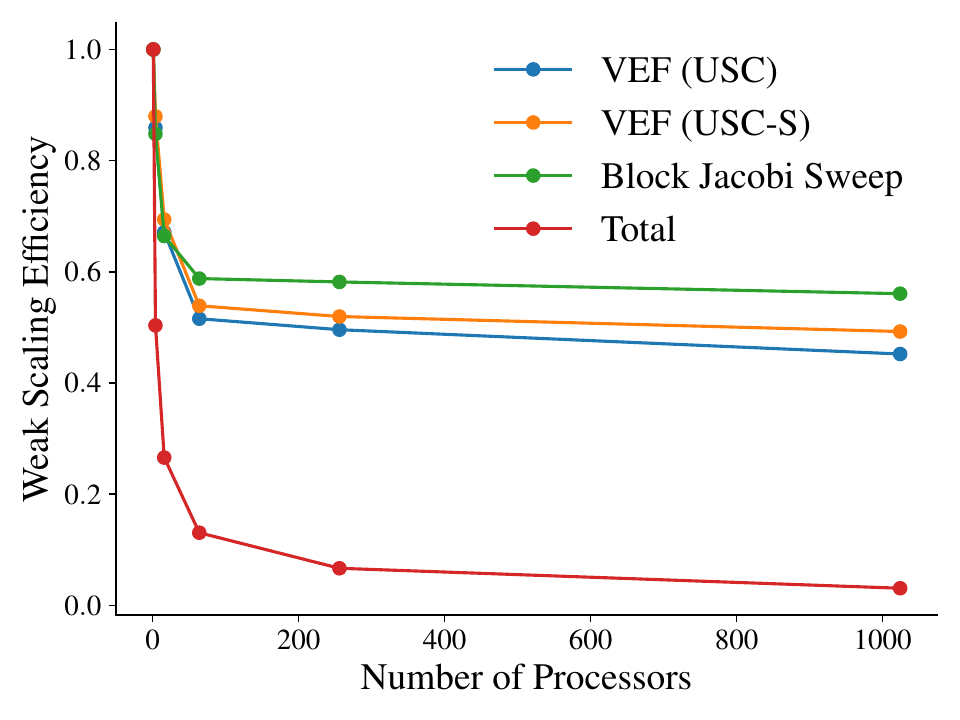}
\caption{\resp[red]{Weak scaling efficiency for the IP VEF fixed-point solve on the crooked pipe problem. The scaling of the average cost per iteration for the VEF solve (which includes assembly, preconditioner setup, and solve) using the USC and USC-S preconditioners are compared. In addition, the scaling of the average cost per iteration for the parallel block Jacobi sweep and the scaling of the total runtime are shown. The VEF and parallel block Jacobi sweeps weak scale but the total runtime does not due to the use of the parallel block Jacobi sweep which causes the total number of outer iterations to increase with the processor count. }}
\label{fig:fweak}
\end{figure}
Weak scaling efficiency is plotted in Fig.~\ref{fig:fweak}. 
Here, the VEF costs include assembling the IP VEF operator, setting up the preconditioner, and solving the IP VEF system with preconditioned BiCGStab. Note that the efficiencies reported in \S\ref{sec:weak_mock} and \S\ref{sec:weak_first_it} only considered the costs associated with the BiCGStab iteration and not assembling the operators or constructing the preconditioners. 
The VEF and sweep costs are averaged over the outer iteration. 
Using the USC-S preconditioner led to a VEF cost per iteration that scaled with an efficiency of 49\% when 1024 processors were used. For the same number of processors, using the USC preconditioner led to a lower efficiency of 45\%. This improved scaling is attributed to USC-S's more uniform convergence with respect to processors. 
While USC-S is more expensive due to its increased assembly costs compared to USC, assembly weak scales efficiently, meaning the improved convergence of USC-S leads to better overall weak scaling efficiency. 

The cost of a single block Jacobi sweep scales well by design from its low communication costs. 
However, the total runtime does not scale due to the increase in outer iterations to convergence as the processor count is increased. We stress that this poor scaling is due to the use of a parallel block Jacobi sweep and may be particularly poor due to the small number of degrees of freedom per processor used in this study. 
Although not investigated here, the use of multiple parallel block Jacobi sweeps per outer iteration may improve the parallel scaling of the full algorithm. 
Furthermore, in time-dependent calculations, the solution at the previous time step can often provide a good initial guess for the parallel block Jacobi sweep, especially in the context of multiphysics problems where the time step size is controlled by other, explicitly integrated physics components. In such case, the initial error in the parallel block Jacobi sweep is much lower, enabling more robust convergence for the outer iteration than observed in this time-independent study. 
\end{response}

\begin{response}[orange]
\subsection{Strong Scaling} \label{sec:strong}
Finally, we present a strong scaling study of the full fixed-point algorithm on the crooked pipe problem from \S\ref{sec:lcp}. The problem size was fixed at \num{28672} equally sized elements with S$_{12}$ angular quadrature and $p=2$ which corresponded to \num{258048} VEF scalar flux unknowns and \num{21676032} angular flux unknowns. Fixed-point iteration without Anderson acceleration is used. The fixed-point tolerance is $10^{-6}$ and the inner BiCGStab tolerance is $10^{-8}$. The streaming and collision operator is approximately inverted using a parallel block Jacobi sweep. Performance is characterized on a single node of the \texttt{rztopaz} machine. 

\begin{table}
\centering
\begin{response}[orange]
\caption{A single-node strong scaling study of the crooked pipe problem with \num{28672} elements, S$_{12}$ angular quadrature, and $p=2$. A parallel block Jacobi transport sweep was used to approximate the inverse of the streaming and collision operator. The IP VEF method was used with the USC preconditioner. Outer refers to the number of fixed-point iterations required to converge to a tolerance of $10^{-6}$. The maximum, minimum, and average number of USC-preconditioned inner BiCGStab iterations required to converge to a tolerance of $10^{-8}$ is shown. The parallel block Jacobi sweep's decreasing accuracy in the early stages of the outer iteration as the processor count increases led to an increase in the maximum percentage of elements requiring the negative flux fixup as well as an increase in the maximum discontinuity in the Eddington tensor. The scaling of the maximum discontinuity in the Eddington tensor induces a scaling of the maximum number of USC-preconditioned BiCGStab iterations. However, the average inner iterations to convergence is uniform with respect to the processor count. }
\label{tab:strong}
\begin{tabular}{cccccccccccc}
\toprule
 &  & \multicolumn{3}{c}{Inner Iterations}  &  & \multicolumn{2}{c}{\% Fixed Up}  &  & \multicolumn{2}{c}{Eddington Jump} \\
\cmidrule{3-5}\cmidrule{7-8}\cmidrule{10-11}
Processors & Outer & Max & Min & Avg. & & Max & Min & & Max & Min \\
\midrule
1 & 47 & 18 & 1 & 8.255 & & 1.36 & 0.50 & & 1.32 & \num{1.38e-02} \\
2 & 57 & 25 & 1 & 9.368 & & 2.42 & 0.50 & & 1.40 & \num{4.18e-02} \\
4 & 68 & 33 & 1 & 10.956 & & 6.19 & 0.61 & & 1.72 & \num{7.80e-02} \\
8 & 76 & 44 & 1 & 11.697 & & 6.43 & 0.45 & & 1.80 & \num{1.17e-01} \\
16 & 83 & 42 & 1 & 11.072 & & 8.55 & 0.55 & & 1.89 & \num{1.49e-01} \\
32 & 112 & 40 & 1 & 11.509 & & 10.82 & 0.42 & & 2.06 & \num{1.47e-01} \\
\bottomrule
\end{tabular}
\end{response}
\end{table}
Table \ref{tab:strong} shows the number of outer iterations and the maximum, minimum, and average number of USC-preconditioned BiCGStab iterations as well as the maximum and minimum percentage of elements that required the negative flux fixup and the maximum and minimum values of $\| \jump{\E\hat{n}} \|_{L^2(\Gamma_0)}$. Due to the parallel block Jacobi sweep, the outer iteration is not robust to increasing processors. 
The maximum discontinuity in the Eddington tensor increases as the processor count increases inducing degradation in the USC preconditioner as indicated by the scaling of the maximum number of inner BiCGStab iterations. 
The average number of iterations is uniform with respect to processor count. The reliance on the negative flux fixup in the first few outer iterations increases with processors as well due to the use of the parallel block Jacobi sweep. 

Strong scaling speedup, defined equivalently to $\varepsilon_n$ in Eq.~\ref{eq:scaling_eff}, is plotted in Fig.~\ref{fig:strong} for the average cost per iteration associated with VEF assembly, the VEF solve, and the parallel block Jacobi sweep along with the total cost of the full algorithm. The dashed line represents the ideal speedup of $n$ times faster when $n$ processors are used. VEF assembly and the block Jacobi sweep strong scale well due to their low communication costs showing speedups using 32 processors of 25x and 23x, respectively. The VEF solve requires communication and thus only achieves a 12x speedup when 32 processors are used. The full algorithm is primarily hindered by the scaling of the number of outer iterations required to converge as the processor count is increased. Overall, the IP VEF method with parallel block Jacobi transport sweep achieved a speedup of 10x using 32 processors. 
\begin{figure}
\centering
\includegraphics[width=.65\textwidth]{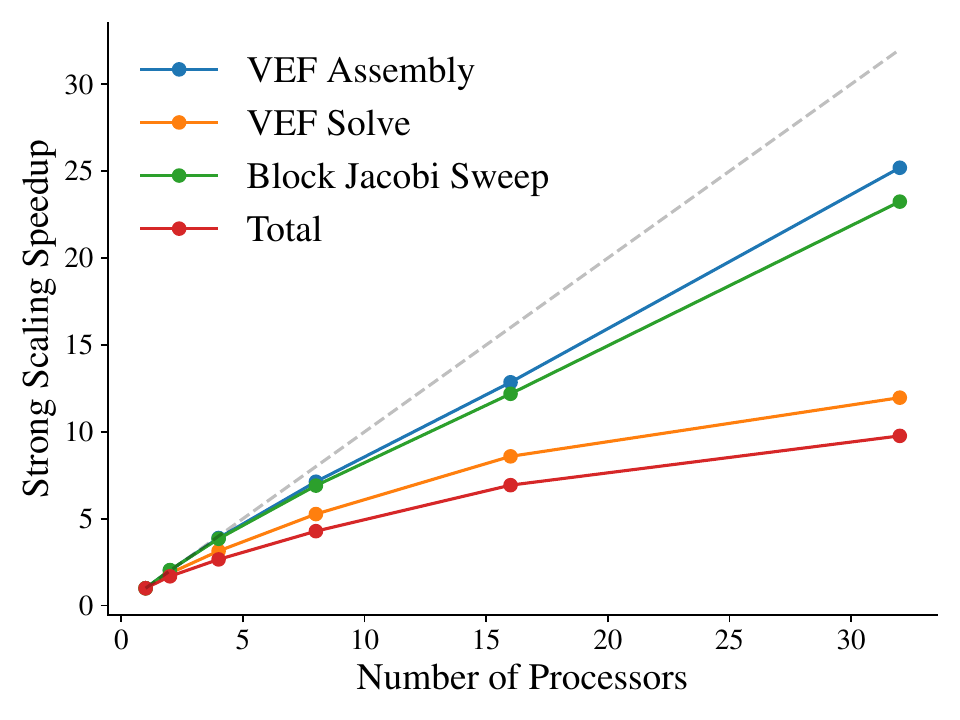}
\caption{\resp[orange]{Strong scaling speedup as a function of the number of processors on the crooked pipe problem with \num{28672} elements, S$_{12}$ angular quadrature, and $p=2$. The average cost per outer fixed-point iteration is shown for assembling and solving the IP VEF linear system of equations and performing the parallel block Jacobi transport sweep along with the total runtime cost. The dashed line represents the ideal scaling of $\varepsilon_n = n$. }}
\label{fig:strong}
\end{figure}
\end{response}

\section{Conclusions}
We have developed a family of high-order discretizations of the Variable Eddington Factor (VEF) equations that are compatible with curved meshes and have efficient preconditioned iterative solvers.
When combined with a high-order Discontinuous Galerkin (DG) discretization of Discrete Ordinates (\Sn) transport, the resulting VEF methods are efficient in both outer fixed-point iterations and inner linear iterations on a challenging proxy problem from thermal radiative transfer (TRT).
We adapted the unified framework for DG methods for elliptic problems presented in \cite{Arnold2002} to the VEF equations to derive analogs of the interior penalty (IP), second method of Bassi and Rebay (BR2), minimal dissipation local Discontinuous Galerkin (MDLDG), and continuous finite element (CG) methods.
\resp[red]{The uniform subspace correction (USC) preconditioner developed by \citet{Pazner2021}, originally designed for DG discretizations of the model Poisson problem, was extended to the IP and BR2 VEF systems and shown to be effective leading to iteration counts independent of the mesh size and polynomial order.}
GMRES convergence estimates for the preconditioned system were derived for the non-symmetric VEF system of equations under relatively mild assumptions.
The MDLDG and CG discretizations were effectively preconditioned by Algebraic Multigrid (AMG). 

The VEF methods were verified to converge with $\mathcal{O}(h^{p+1})$ on refinements of a third-order mesh using a quadratically anisotropic manufactured solution. They were also tested in the thick diffusion limit both on an orthogonal mesh and a severely distorted third-order mesh generated with a Lagrangian hydrodynamics code. In both cases, all the VEF methods preserved the thick diffusion limit and converged robustly. Convergence on the triple point mesh indicates that these methods are robust to extreme mesh distortions and inexact transport inversions arising from reentrant faces. 

The methods were also tested on a linearized crooked pipe problem. This problem had a 1000x difference in total cross section and was designed to emulate \resp[red]{the first time step of} a time-dependent TRT calculation. Using the stabilized bi-conjugate gradient method (BiCGStab), all of the VEF methods were efficiently solvable independent from the mesh size, polynomial order, and, if present, penalty parameter. Using a small Anderson space of size two, each VEF algorithm converged in a uniform number of outer Anderson-accelerated fixed-point iterations as well. 
\resp[orange]{The CG and IP methods were the fastest in overall runtime due to their lower assembly and solve costs compared to BR2 and MDLDG. 
However, we note that due to the relative dominance of the cost of the transport sweep, the choice of the VEF discretization led to a relative variance of only 10\% in total runtime. 
It was observed that the more numerically diffusive methods, namely the IP, BR2, and CG methods, produced scattering sources that led to negativities in the transport sweep more so than for the minimally dissipative MDLDG method. This led to fewer elements requiring the application of a negative flux fixup for the MDLDG method compared to IP, BR2, and CG. 
It was also observed that the IP, BR2, and CG methods converged nearly identically, indicating that the stabilization terms used by the IP and BR2 methods cause the overall algorithm to behave as if a continuous solution representation were used. }

\resp[blue]{The solution from the IP VEF method was compared to a DG \Sn method preconditioned with Diffusion Synthetic Acceleration (DSA) on the thick diffusion limit and crooked pipe problems. Using a fixed angular quadrature rule, the solutions from VEF and \Sn were compared as the mesh was refined. The VEF and \Sn solutions converged to each other with the optimal spatial order of accuracy on the thick diffusion limit problem. However, the two schemes converged to each other with only first order accuracy on the multiple material crooked pipe problem. Such behavior may be due to the use of a negative flux fixup for the VEF method or due to under resolution in space and/or angle. 
These comparisons suggest that independent VEF methods do in fact converge to the transport solution. 
However, the first-order convergence in space between an independent VEF method and an \Sn method observed on the multi-material crooked pipe problem should be investigated further in order to isolate the cause of the sub-optimal convergence. Ideally, convergence would be compared as both the spatial mesh is refined and the \Sn order is increased and, if possible, in a regime where the spatial mesh is resolved enough to not require a negative flux fixup.
}

\begin{response}[orange][ipdiffcomp]
Finally, weak and strong scaling studies were performed on the IP VEF method with $p=2$. A non-physically difficult mock problem with discontinuous VEF data was found to cause non-convergence of the USC-preconditioned linear solver. Uniform convergence was recovered by using a symmetrized variant of the USC preconditioner. On the physically realistic crooked pipe problem, both the standard USC and the symmetrized USC preconditioners performed well. Solving the non-symmetric IP VEF system using the symmetrized USC preconditioner on a problem with 43 million scalar flux unknowns and 1024 processors was only 12\% more expensive than solving the symmetric positive definite linear system corresponding to an IP discretization of radiation diffusion. Furthermore, the IP VEF solvers were shown to have weak scaling efficiency out to 1024 processors comparable to that of solving IP radiation diffusion.
\end{response}

\resp[red]{Parallel performance was investigated when the IP VEF method was coupled to a parallel block Jacobi transport sweep. Such a method performs a transport sweep on each processor domain independently using lagged angular flux information on each processor domain's inflow boundary. Since the streaming and collision operator is not inverted exactly at each fixed-point iteration, the number of fixed-point iterations until convergence grew with increasing processor counts, leading to poor weak scaling efficiency for the full fixed-point solve. However, the inner IP VEF computations weak scaled efficiently. It was observed that use of the parallel block Jacobi sweep led to increased need for the negative flux fixup in the sweep and an increase in non-physical discontinuities in the Eddington tensor in the beginning stages of the fixed-point iteration. Solver efficiency was mildly improved by using the symmetrized USC preconditioner in place of the standard USC preconditioner which was shown on the problem with mock VEF data to be more robust to discontinuities in the Eddington tensor.
However, due to the additional assembly costs associated with forming the symmetrized operator, use of the symmetrized preconditioner was overall more expensive than the standard USC preconditioner. 
Similar increases in reliance on the negative flux fixup and discontinuity in the Eddington tensor due to the parallel block Jacobi sweep were seen on a single-node strong scaling study as well. Using the parallel block Jacobi sweep, the full IP VEF algorithm achieved a 10x speedup on 32 processors. }


The primary takeaway from this work is that all of the VEF methods presented here are strong candidates for implementation in a TRT algorithm. All of the methods were robust to the thick diffusion limit, inexact sweeps from reentrant faces, and strongly heterogeneous materials. 
\resp[red]{While the MDLDG method is notable for its reduced reliance on the negative flux fixup, the four methods had only minor differences in runtime due to the large expense of the transport sweep relative to the cost of forming and inverting the discrete VEF equations. 
Thus, the methods were primarily differentiated by their ease of implementation in forming both the discrete linear system and its associated preconditioner. The CG and IP methods are the simplest to implement since they do not use the complicated lifting operators needed by the BR2 and MDLDG methods. 
On the other hand, the CG and MDLDG methods are simpler to solve since they only require a black box AMG solver whereas IP and BR2 require the additional complications of forming the continuous operator and performing a Jacobi iteration on the interfacial unknowns. 
Thus, the CG method is recommended for its simpler assembly and preconditioning requirements compared to the other methods. 
For solving the IP or BR2 linear system, both the USC and the symmetrized USC preconditioner were efficient on realistic problems with the standard USC method typically being the cheaper option and the symmetrized method being more robust in terms of iterations to convergence.}

\resp[red]{In the future, algorithmic aspects associated with extending the methods to the full thermal radiative transfer or radiation hydrodynamics problems need to be developed. In particular, it may be interesting to compare the solution quality of the four VEF methods on a problem such as the Marshak wave which has a solution with discontinuous derivatives in space. Such a problem could expose differences in solution quality not seen in this study, particularly for the CG VEF method paired with a DG transport discretization.}

\section{Acknowledgements}
This work was performed under the auspices of the U.S. Department of Energy by Lawrence Livermore National Laboratory under Contract DE-AC52-07NA27344 (LLNL-JRNL-829396). S.O. was supported by the U.S. Department of Energy, Office of Science, Office of Advanced Scientific Computing Research, and the Department of Energy Computational Science Graduate Fellowship under Award Number DE-SC0019323. 
W.P. was partially supported by the LLNL-LDRD Program under Project Number 20-ERD-002.

This report was prepared as an account of work sponsored by an agency of the United States Government. Neither the United States Government nor any agency thereof, nor any of their employees, makes any warranty, express or implied, or assumes any legal liability or responsibility for the accuracy, completeness, or usefulness of any information, apparatus, product, or process disclosed, or represents that its use would not infringe privately owned rights. Reference herein to any specific commercial product, process, or service by trade name, trademark, manufacturer, or otherwise does not necessarily constitute or imply its endorsement, recommendation, or favoring by the United States Government or any agency thereof. The views and opinions of authors expressed herein do not necessarily state or reflect those of the United States Government or any agency thereof.

\bibliographystyle{IEEEtranN}
\bibliography{references}

\begin{thebibliography}{69}
\providecommand{\natexlab}[1]{#1}
\providecommand{\url}[1]{#1}
\csname url@samestyle\endcsname
\providecommand{\newblock}{\relax}
\providecommand{\bibinfo}[2]{#2}
\providecommand{\BIBentrySTDinterwordspacing}{\spaceskip=0pt\relax}
\providecommand{\BIBentryALTinterwordstretchfactor}{4}
\providecommand{\BIBentryALTinterwordspacing}{\spaceskip=\fontdimen2\font plus
\BIBentryALTinterwordstretchfactor\fontdimen3\font minus
  \fontdimen4\font\relax}
\providecommand{\BIBforeignlanguage}[2]{{%
\expandafter\ifx\csname l@#1\endcsname\relax
\typeout{** WARNING: IEEEtranN.bst: No hyphenation pattern has been}%
\typeout{** loaded for the language `#1'. Using the pattern for}%
\typeout{** the default language instead.}%
\else
\language=\csname l@#1\endcsname
\fi
#2}}
\providecommand{\BIBdecl}{\relax}
\BIBdecl

\bibitem[Haut et~al.(2019)Haut, Maginot, Tomov, Southworth, Brunner, and
  Bailey]{graph_sweeps}
T.~S. Haut, P.~G. Maginot, V.~Z. Tomov, B.~S. Southworth, T.~A. Brunner, and
  T.~S. Bailey, ``An efficient sweep-based solver for the {$S_N$} equations on
  high-order meshes,'' \emph{Nuclear Science and Engineering}, 2019.

\bibitem[Arnold et~al.(2002)Arnold, Brezzi, Cockburn, and Marini]{Arnold2002}
D.~N. Arnold, F.~Brezzi, B.~Cockburn, and L.~D. Marini, ``Unified analysis of
  discontinuous {G}alerkin methods for elliptic problems,'' \emph{SIAM Journal
  on Numerical Analysis}, vol.~39, no.~5, pp. 1749--1779, 2002.

\bibitem[Pazner and Kolev(2021)]{Pazner2021}
W.~Pazner and T.~Kolev, ``Uniform subspace correction preconditioners for
  discontinuous {G}alerkin methods with $hp$-refinement,'' \emph{Communications
  on Applied Mathematics and Computation}, Jul. 2021.

\bibitem[{D. Mihalas}(1978)]{mihalas}
{D. Mihalas}, \emph{{Stellar Atmospheres}}.\hskip 1em plus 0.5em minus
  0.4em\relax W.~H.~Freeman and Co, 1978.

\bibitem[{V. Ya. Gol'din}(1964)]{goldin}
{V. Ya. Gol'din}, ``A quasi-diffusion method of solving the kinetic equation,''
  \emph{USSR Comp. Math. and Math. Physics}, vol.~4, pp. 136--149, 1964.

\bibitem[Aristova and Baydin(2013)]{airstova_eigenvalue}
E.~Aristova and D.~Baydin, ``Implementation of the quasidiffusion method for
  calculating the critical parameters of a fast neutron reactor in {3D}
  hexagonal geometry,'' \emph{Mathematical Models and Computer Simulations},
  vol.~5, pp. 145--155, 2013.

\bibitem[Tamang and Anistratov(2014)]{doi:10.13182/NSE13-42}
A.~Tamang and D.~Y. Anistratov, ``A multilevel projective method for solving
  the space-time multigroup neutron kinetics equations coupled with the heat
  transfer equation,'' \emph{Nuclear Science and Engineering}, vol. 177, no.~1,
  pp. 1--18, 2014.

\bibitem[Anistratov et~al.(1996)Anistratov, Aristova, and
  Gol'din]{anistratov1996nonlinear}
D.~Y. Anistratov, E.~N. Aristova, and V.~Y. Gol'din, ``A nonlinear method for
  solving the problems of radiation transfer in medium,'' \emph{Matematicheskoe
  modelirovanie}, vol.~8, no.~12, pp. 3--28, 1996.

\bibitem[Anistratov and Gol'din(1986)]{anistratov_fvm}
D.~Y. Anistratov and V.~Y. Gol'din, ``Comparison of difference schemes for the
  quasidiffusion method for solving the transport equation,'' \emph{Problems of
  Atomic Science and Engineering: Method and Codes for Numerical Solution
  Mathematical Physics Problems}, vol.~2, pp. 17--23, 1986.

\bibitem[Aristova and Gol'din(2000)]{ARISTOVA2000139}
E.~Aristova and V.~Gol'din, ``Computation of anisotropy scattering of solar
  radiation in atmosphere (monoenergetic case),'' \emph{Journal of Quantitative
  Spectroscopy and Radiative Transfer}, vol.~67, no.~2, pp. 139--157, 2000.

\bibitem[Alcouffe(1977)]{A}
R.~Alcouffe, ``Diffusion synthetic acceleration methods for the
  diamond-differenced discrete-ordinates equations,'' \emph{Nuclear Science and
  Engineering}, vol.~64, pp. 344--355, 1977.

\bibitem[Anistratov and Gol'din(1993)]{doi:10.1080/00411459308203810}
D.~Y. Anistratov and V.~Y. Gol'din, ``Nonlinear methods for solving particle
  transport problems,'' \emph{Transport Theory and Statistical Physics},
  vol.~22, no. 2-3, pp. 125--163, 1993.

\bibitem[Warsa and Anistratov(2018)]{two-level-independent-warsa}
J.~Warsa and D.~Anistratov, ``Two-level transport methods with independent
  discretization,'' \emph{Journal of Computational and Theoretical Transport},
  vol.~47, no. 4-6, pp. 424--450, 2018.

\bibitem[Larsen et~al.(1987)Larsen, Morel, and Warren F.~Miller]{diflim}
E.~W. Larsen, J.~Morel, and J.~Warren F.~Miller, ``Asymptotic solutions of
  numerical transport problems in optically thick, diffusive regimes,''
  \emph{Journal of Computational Physics}, vol.~69, pp. 283--324, 1987.

\bibitem[Ghassemi and Anistratov(2020)]{GHASSEMI2020109315}
P.~Ghassemi and D.~Y. Anistratov, ``Multilevel quasidiffusion method with
  mixed-order time discretization for multigroup thermal radiative transfer
  problems,'' \emph{Journal of Computational Physics}, vol. 409, p. 109315,
  2020.

\bibitem[Yee et~al.(2021)Yee, Olivier, Southworth, Holec, and Haut]{yee_mc21}
B.~Yee, S.~Olivier, B.~Southworth, M.~Holec, and T.~Haut, ``A new scheme for
  solving high-order {DG} discretizations of thermal radiative transfer using
  the variable {E}ddington factor method,'' in \emph{Proceedings of the
  International Conference on Mathematics and Computational Methods Applied to
  Nuclear Science and Engineering (M\&C 2021)}, 2021.

\bibitem[Anistratov and Coale(2021)]{anistratov2021implicit}
D.~Y. Anistratov and J.~M. Coale, ``Implicit methods with reduced memory for
  thermal radiative transfer,'' 2021.

\bibitem[Jiang et~al.(2012)Jiang, Stone, and Davis]{Jiang_2012}
Y.-F. Jiang, J.~M. Stone, and S.~W. Davis, ``A {G}odunov method for
  multidimensional radiation magnetohydrodynamics based on a variable
  {E}ddington tensor,'' \emph{The Astrophysical Journal Supplement Series},
  vol. 199, no.~1, p.~14, 2012.

\bibitem[Gnedin and Abel(2001)]{GNEDIN2001437}
N.~Y. Gnedin and T.~Abel, ``Multi-dimensional cosmological radiative transfer
  with a variable {E}ddington tensor formalism,'' \emph{New Astronomy}, vol.~6,
  no.~7, pp. 437--455, 2001.

\bibitem[Gehmeyr and Mihalas(1994)]{GEHMEYR1994320}
M.~Gehmeyr and D.~Mihalas, ``Adaptive grid radiation hydrodynamics with
  {TITAN},'' \emph{Physica D: Nonlinear Phenomena}, vol.~77, no.~1, pp.
  320--341, 1994.

\bibitem[Davis et~al.(2012)Davis, Stone, and Jiang]{Davis_2012}
S.~W. Davis, J.~M. Stone, and Y.-F. Jiang, ``A radiation transfer solver for
  {Athena} using short characteristics,'' \emph{The Astrophysical Journal
  Supplement Series}, vol. 199, no.~1, p.~9, feb 2012.

\bibitem[Olivier and Morel(2017)]{me}
S.~S. Olivier and J.~E. Morel, ``Variable {E}ddington factor method for the
  {$S_N$} equations with lumped discontinuous {G}alerkin spatial discretization
  coupled to a drift-diffusion acceleration equation with mixed finite-element
  discretization,'' \emph{Journal of Computational and Theoretical Transport},
  vol.~46, no. 6-7, pp. 480--496, 2017.

\bibitem[Lou et~al.(2019)Lou, Morel, and Gentile]{LOU2019258}
J.~Lou, J.~E. Morel, and N.~Gentile, ``A variable {E}ddington factor method for
  the {1-D} grey radiative transfer equations with discontinuous {G}alerkin and
  mixed finite-element spatial differencing,'' \emph{Journal of Computational
  Physics}, vol. 393, pp. 258--277, 2019.

\bibitem[Yee et~al.(2020)Yee, Olivier, Haut, Holec, Tomov, and
  Maginot]{YEE2020109696}
B.~C. Yee, S.~S. Olivier, T.~S. Haut, M.~Holec, V.~Z. Tomov, and P.~G. Maginot,
  ``A quadratic programming flux correction method for high-order {DG}
  discretizations of {$S_N$} transport,'' \emph{Journal of Computational
  Physics}, vol. 419, p. 109696, 2020.

\bibitem[Anistratov(2019)]{ANISTRATOV2019186}
D.~Y. Anistratov, ``Stability analysis of a multilevel quasidiffusion method
  for thermal radiative transfer problems,'' \emph{Journal of Computational
  Physics}, vol. 376, pp. 186--209, 2019.

\bibitem[Adams and Larsen(2002)]{AL}
M.~Adams and E.~Larsen, ``Fast iterative methods for discrete-ordinates
  particle transport calculations,'' \emph{Progress in Nuclear Energy},
  vol.~40, no.~1, pp. 3--159, 2002.

\bibitem[Miften and Larsen(1993)]{QDBC}
M.~Miften and E.~Larsen, ``The quasi-diffusion method for solving transport
  problems in planar and spherical geometries,'' \emph{Journal of Transport
  Theory and Statistical Physics}, vol. 22(2-3), pp. 165--186, 1993.

\bibitem[Jones(2019)]{Jones2019TheQM}
J.~P. Jones, ``The quasidiffusion method for solving radiation transport
  problems on arbitrary quadrilateral meshes in {2D} r-z geometry.'' Ph.D.
  dissertation, North Carolina State University, 2019.

\bibitem[Wieselquist et~al.(2014)Wieselquist, Anistratov, and
  Morel]{WIESELQUIST2014343}
W.~A. Wieselquist, D.~Y. Anistratov, and J.~E. Morel, ``A cell-local finite
  difference discretization of the low-order quasidiffusion equations for
  neutral particle transport on unstructured quadrilateral meshes,''
  \emph{Journal of Computational Physics}, vol. 273, pp. 343--357, 2014.

\bibitem[Vallette(2002)]{vallette}
N.~D. Vallette, ``Discretisation and solution of quasi-diffusion equations,''
  Master's thesis, Texas A\&M University, 2002.

\bibitem[Olivier et~al.(2019)Olivier, Maginot, and Haut]{olivier_mandc}
S.~Olivier, P.~Maginot, and T.~Haut, ``High order mixed finite element
  discretization for the variable {E}ddington factor equations,'' in
  \emph{Proceedings of the International Conference on Mathematics and
  Computational Methods Applied to Nuclear Science and Engineering (M\&C
  2019)}, 2019.

\bibitem[{W.A. Wieselquist}(2010)]{wieselquist}
{W.A. Wieselquist}, ``A low-order quasidiffusion discretization via
  linear-continuous finite-elements on unstructured triangular meshes,'' in
  \emph{Proceedings of PHYSOR 2010: Advances in Reactor Physics to Power the
  Nuclear Renaissance}.\hskip 1em plus 0.5em minus 0.4em\relax The American
  Nuclear Society, 2010.

\bibitem[Anistratov and Warsa(2018)]{dima_dfem}
D.~Y. Anistratov and J.~S. Warsa, ``Discontinuous finite element
  quasi-diffusion methods,'' \emph{Nuclear Science and Engineering}, vol. 191,
  no.~2, pp. 105--120, 2018.

\bibitem[Benzi et~al.(2005)Benzi, Golub, and Liesen]{benzi_golub_liesen_2005}
M.~Benzi, G.~H. Golub, and J.~Liesen, ``Numerical solution of saddle point
  problems,'' \emph{Acta Numerica}, vol.~14, p. 1–137, 2005.

\bibitem[Warsa et~al.(2004{\natexlab{a}})Warsa, Benzi, Wareing, and
  Morel]{warsa_mfem}
J.~S. Warsa, M.~Benzi, T.~A. Wareing, and J.~E. Morel, ``Preconditioning a
  mixed discontinuous finite element method for radiation diffusion,''
  \emph{Numerical Linear Algebra with Applications}, vol.~11, pp. 795--811,
  2004.

\bibitem[Dobrev et~al.(2012)Dobrev, Kolev, and Rieben]{blast}
V.~Dobrev, T.~Kolev, and R.~Rieben, ``High-order curvilinear finite element
  methods for {L}agrangian hydrodynamics,'' \emph{SIAM Journal on Scientific
  Computing}, vol.~34, pp. B606--B641, 2012.

\bibitem[Dobrev et~al.(2013)Dobrev, Ellis, Kolev, and Rieben]{blast2}
V.~Dobrev, T.~Ellis, T.~Z. Kolev, and R.~Rieben, ``High-order curvilinear
  finite elements for axisymmetric {L}agrangian hydrodynamics,'' \emph{Comput.
  Fluids}, pp. 58--69, 2013.

\bibitem[Anderson et~al.(2018)Anderson, Dobrev, Kolev, Rieben, and
  Tomov]{blast3}
R.~W. Anderson, V.~A. Dobrev, T.~V. Kolev, R.~N. Rieben, and V.~Z. Tomov,
  ``High-order multi-material {ALE} hydrodynamics,'' \emph{SIAM J. Sci. Comp.},
  vol.~40, no.~1, pp. B32--B58, 2018.

\bibitem[Woods(2018)]{woods_thesis}
D.~Woods, ``Discrete ordinates radiation transport using high-order finite
  element spatial discretizations on meshes with curved surfaces,'' Ph.D.
  dissertation, Oregon State University, 2018.

\bibitem[Haut et~al.(2020)Haut, Southworth, Maginot, and Tomov]{ldrd_dsa}
T.~S. Haut, B.~S. Southworth, P.~G. Maginot, and V.~Z. Tomov, ``Diffusion
  synthetic acceleration preconditioning for discontinuous {G}alerkin
  discretizations of {$S_N$} transport on high-order curved meshes,''
  \emph{SIAM Journal on Scientific Computing}, vol.~42, no.~5, pp.
  B1271--B1301, 2020.

\bibitem[Southworth et~al.(2021)Southworth, Holec, and
  Haut]{doi:10.1080/00295639.2020.1799603}
B.~S. Southworth, M.~Holec, and T.~S. Haut, ``Diffusion synthetic acceleration
  for heterogeneous domains, compatible with voids,'' \emph{Nuclear Science and
  Engineering}, vol. 195, no.~2, pp. 119--136, 2021.

\bibitem[Warsa et~al.(2004{\natexlab{b}})Warsa, Wareing, and
  Morel]{doi:10.13182/NSE02-14}
J.~S. Warsa, T.~A. Wareing, and J.~E. Morel, ``Krylov iterative methods and the
  degraded effectiveness of diffusion synthetic acceleration for
  multidimensional {$S_N$} calculations in problems with material
  discontinuities,'' \emph{Nuclear Science and Engineering}, vol. 147, no.~3,
  pp. 218--248, 2004.

\bibitem[Ciarlet(2002)]{ciarlet2002finite}
P.~Ciarlet, \emph{The Finite Element Method for Elliptic Problems}, ser.
  Classics in Applied Mathematics.\hskip 1em plus 0.5em minus 0.4em\relax
  Society for Industrial and Applied Mathematics, 2002.

\bibitem[Brezzi et~al.(2000{\natexlab{a}})Brezzi, Manzini, Marini, Pietra, and
  Russo]{brstab}
\BIBentryALTinterwordspacing
F.~Brezzi, G.~Manzini, D.~Marini, P.~Pietra, and A.~Russo, ``Discontinuous
  galerkin approximations for elliptic problems,'' \emph{Numerical Methods for
  Partial Differential Equations}, vol.~16, no.~4, pp. 365--378, 2000.
  [Online]. Available:
  \url{https://onlinelibrary.wiley.com/doi/abs/10.1002/1098-2426%28200007%2916%3A4%3C365%3A%3AAID-NUM2%3E3.0.CO%3B2-Y}
\BIBentrySTDinterwordspacing

\bibitem[Cockburn and Dong(2007)]{10.1007/s10915-007-9130-3}
B.~Cockburn and B.~Dong, ``An analysis of the minimal dissipation local
  discontinuous {G}alerkin method for convection--diffusion problems,''
  \emph{J. Sci. Comput.}, vol.~32, no.~2, p. 233–262, Aug. 2007.

\bibitem[Quarteroni and Valli(1994)]{quateroni}
A.~Quarteroni and A.~Valli, \emph{Numerical Approximation of Partial
  Differential Equations}.\hskip 1em plus 0.5em minus 0.4em\relax Springer,
  Berlin, Heidelberg, 1994.

\bibitem[Xu(1992)]{Xu1992}
J.~Xu, ``Iterative methods by space decomposition and subspace correction,''
  \emph{SIAM Review}, vol.~34, no.~4, pp. 581--613, 1992.

\bibitem[Xu and Zikatanov(2002)]{Xu2002}
J.~Xu and L.~Zikatanov, ``The method of alternating projections and the method
  of subspace corrections in {H}ilbert space,'' \emph{Journal of the American
  Mathematical Society}, vol.~15, no.~03, pp. 573--598, Jul 2002.

\bibitem[Toselli and Widlund(2005)]{Toselli2005}
A.~Toselli and O.~B. Widlund, \emph{Domain Decomposition Methods --- Algorithms
  and Theory}.\hskip 1em plus 0.5em minus 0.4em\relax Springer Berlin
  Heidelberg, 2005.

\bibitem[Antonietti et~al.(2016)Antonietti, Sarti, Verani, and
  Zikatanov]{Antonietti2016}
P.~F. Antonietti, M.~Sarti, M.~Verani, and L.~T. Zikatanov, ``A uniform
  additive {S}chwarz preconditioner for high-order discontinuous {G}alerkin
  approximations of elliptic problems,'' \emph{Journal of Scientific
  Computing}, vol.~70, no.~2, pp. 608--630, Aug 2016.

\bibitem[Dobrev et~al.(2006)Dobrev, Lazarov, Vassilevski, and
  Zikatanov]{Dobrev2006}
V.~A. Dobrev, R.~D. Lazarov, P.~S. Vassilevski, and L.~T. Zikatanov,
  ``Two-level preconditioning of discontinuous {G}alerkin approximations of
  second-order elliptic equations,'' \emph{Numerical Linear Algebra with
  Applications}, vol.~13, no.~9, pp. 753--770, 2006.

\bibitem[O'Malley et~al.(2017)O'Malley, K{\'{o}}ph{\'{a}}zi, Smedley-Stevenson,
  and Eaton]{OMalley2017}
\BIBentryALTinterwordspacing
B.~O'Malley, J.~K{\'{o}}ph{\'{a}}zi, R.~Smedley-Stevenson, and M.~Eaton,
  ``Hybrid multi-level solvers for discontinuous {G}alerkin finite element
  discrete ordinate diffusion synthetic acceleration of radiation transport
  algorithms,'' \emph{Annals of Nuclear Energy}, vol. 102, pp. 134--147, Apr.
  2017. [Online]. Available:
  \url{https://doi.org/10.1016%2Fj.anucene.2016.11.048}
\BIBentrySTDinterwordspacing

\bibitem[Warsa et~al.(2003)Warsa, Benzi, Wareing, and Morel]{Warsa2003}
\BIBentryALTinterwordspacing
J.~S. Warsa, M.~Benzi, T.~A. Wareing, and J.~E. Morel, ``Two-level
  preconditioning of a discontinuous {G}alerkin method for radiation
  diffusion,'' in \emph{Numerical Mathematics and Advanced Applications}.\hskip
  1em plus 0.5em minus 0.4em\relax Springer Milan, 2003, pp. 967--977.
  [Online]. Available: \url{https://doi.org/10.1007%2F978-88-470-2089-4_88}
\BIBentrySTDinterwordspacing

\bibitem[Pazner(2020)]{Pazner2020}
W.~Pazner, ``Efficient low-order refined preconditioners for high-order
  matrix-free continuous and discontinuous {G}alerkin methods,'' \emph{{SIAM}
  Journal on Scientific Computing}, vol.~42, no.~5, pp. A3055--A3083, Jan.
  2020.

\bibitem[Lions(1988)]{Lions1988}
P.~L. Lions, ``On the {S}chwarz alternating method. {I}.'' in \emph{Domain
  Decomposition Methods for Partial Differential Equations}, R.~Glowinski,
  G.~Golub, G.~Meurant, and J.~P{\'e}riaux, Eds.\hskip 1em plus 0.5em minus
  0.4em\relax SIAM, 1988.

\bibitem[Falgout and Yang(2002)]{hypre}
R.~D. Falgout and U.~M. Yang, ``Hypre: A library of high performance
  preconditioners,'' in \emph{Proceedings of the International Conference on
  Computational Science-Part III}, ser. ICCS '02.\hskip 1em plus 0.5em minus
  0.4em\relax Berlin, Heidelberg: Springer-Verlag, 2002, p. 632–641.

\bibitem[Antonietti and Houston(2010)]{Antonietti2010}
P.~F. Antonietti and P.~Houston, ``A class of domain decomposition
  preconditioners for $hp$-discontinuous {G}alerkin finite element methods,''
  \emph{Journal of Scientific Computing}, vol.~46, no.~1, pp. 124--149, Jun.
  2010.

\bibitem[Brezzi et~al.(2000{\natexlab{b}})Brezzi, Manzini, Marini, Pietra, and
  Russo]{Brezzi2000}
F.~Brezzi, G.~Manzini, D.~Marini, P.~Pietra, and A.~Russo, ``Discontinuous
  {G}alerkin approximations for elliptic problems,'' \emph{Numerical Methods
  for Partial Differential Equations}, vol.~16, no.~4, pp. 365--378, 2000.

\bibitem[Eisenstat et~al.(1983)Eisenstat, Elman, and Schultz]{Eisenstat1983}
S.~C. Eisenstat, H.~C. Elman, and M.~H. Schultz, ``Variational iterative
  methods for nonsymmetric systems of linear equations,'' \emph{SIAM Journal on
  Numerical Analysis}, vol.~20, no.~2, pp. 345--357, 1983.

\bibitem[Cai(1989)]{Cai1989}
X.-C. Cai, ``Some domain decomposition algorithms for nonselfadjoint elliptic
  and parabolic partial differential equations,'' Ph.D. dissertation, New York
  University, 1989.

\bibitem[Cai(1990)]{Cai1990}
------, ``An additive schwarz algorithm for nonselfadjoint elliptic
  equations,'' in \emph{Third International Symposium on Domain Decomposition
  Methods for Partial Differential Equations}.\hskip 1em plus 0.5em minus
  0.4em\relax SIAM, Philadelphia, 1990, pp. 232--244.

\bibitem[Anderson et~al.(2020)Anderson, Andrej, Barker, Bramwell, Camier,
  Cerveny, Dobrev, Dudouit, Fisher, Kolev, Pazner, Stowell, Tomov, Dahm,
  Medina, and Zampini]{MFEM}
R.~Anderson, J.~Andrej, A.~Barker, J.~Bramwell, J.-S. Camier, J.~Cerveny,
  V.~Dobrev, Y.~Dudouit, A.~Fisher, T.~Kolev, W.~Pazner, M.~Stowell, V.~Tomov,
  J.~Dahm, D.~Medina, and S.~Zampini, ``{MFEM}: a modular finite element
  methods library,'' \emph{Computers {\&} Mathematics with Applications}, Jul.
  2020.

\bibitem[MFEM()]{mfem-web}
``{MFEM}: Modular finite element methods {[Software]},''
  \url{https://mfem.org}, 2010.

\bibitem[Hindmarsh et~al.(2005)Hindmarsh, Brown, Grant, Lee, Serban, Shumaker,
  and Woodward]{hindmarsh2005sundials}
A.~C. Hindmarsh, P.~N. Brown, K.~E. Grant, S.~L. Lee, R.~Serban, D.~E.
  Shumaker, and C.~S. Woodward, ``{SUNDIALS}: Suite of nonlinear and
  differential/algebraic equation solvers,'' \emph{ACM Transactions on
  Mathematical Software (TOMS)}, vol.~31, no.~3, pp. 363--396, 2005.

\bibitem[Li and Demmel(2003)]{lidemmel03}
X.~S. Li and J.~W. Demmel, ``{SuperLU\_DIST}: A scalable distributed-memory
  sparse direct solver for unsymmetric linear systems,'' \emph{ACM Trans.
  Mathematical Software}, vol.~29, no.~2, pp. 110--140, June 2003.

\bibitem[Bassi and Rebay(2000)]{Bassi2000}
F.~Bassi and S.~Rebay, ``A high order discontinuous {G}alerkin method for
  compressible turbulent flows,'' in \emph{Discontinuous Galerkin Methods},
  B.~Cockburn, G.~E. Karniadakis, and C.-W. Shu, Eds.\hskip 1em plus 0.5em
  minus 0.4em\relax Springer Berlin Heidelberg, 2000, pp. 77--88.

\bibitem[Arnold(1982)]{Arnold1982}
D.~N. Arnold, ``An interior penalty finite element method with discontinuous
  elements,'' \emph{SIAM Journal on Numerical Analysis}, vol.~19, no.~4, pp.
  742--760, Aug. 1982.

\bibitem[Ainsworth et~al.(2011)Ainsworth, Andriamaro, and
  Davydov]{doi:10.1137/11082539X}
M.~Ainsworth, G.~Andriamaro, and O.~Davydov, ``Bernstein–bézier finite
  elements of arbitrary order and optimal assembly procedures,'' \emph{SIAM
  Journal on Scientific Computing}, vol.~33, no.~6, pp. 3087--3109, 2011.

\bibitem[Hamilton et~al.(2009)Hamilton, Benzi, and Warsa]{hamilton2009negative}
S.~Hamilton, M.~Benzi, and J.~Warsa, ``Negative flux fixups in discontinuous
  finite element {$S_N$} transport,'' in \emph{International Conference on
  Mathematics, Computational Methods and Reactor Physics (M\&C 2009), American
  Nuclear Society, LaGrange Park, Illinois, USA}.\hskip 1em plus 0.5em minus
  0.4em\relax Citeseer, 2009.

\end{thebibliography}

\appendix 
\section{Implementation of Lifting Operators} \label{sec:lifting}
Consider the face-local lifting operator $\vec{\rho}_f(\omega)$ used in the BR2 stabilization term defined in Eq.~\ref{eq:rhof} with $\omega = \jump{u}$ which satisfies
	\begin{equation} \label{eq:lift_particularize}
		\int \vec{v}\cdot\vec{\rho}_f(\jump{u}) \ud \x = -\int_{f} \avg{\vec{v}\cdot\hat{n}} \jump{u} \ud s \,, \quad \forall \vec{v} \in W_p\,, \quad \text{on} \ f \in \Gamma_0 \,. 
	\end{equation}
Let $\underline{y}$ represent the vector of DOFs corresponding to a $Y_p$ or $W_p$ grid function $y$. Let $\vec{v},\vec{w} \in W_p$ and define 
	\begin{equation}
		\underline{v}^T \mat{M} \underline{w} = \int \vec{v}\cdot\vec{w} \ud \x 
	\end{equation}
as the $W_p$ mass matrix. Further, define 
	\begin{equation}
		\underline{v}^T \mat{A}_f \underline{u} = - \int_f \avg{\vec{v}\cdot\hat{n}} \jump{u} \ud s \,, \quad \text{on} \ f \in \Gamma_0 \,, 
	\end{equation}
for $u \in Y_p$. 
Equation \ref{eq:lift_particularize} is then equivalent to 
	\begin{equation}
		\mat{M} \underline{\rho}_f(\jump{u}) = \mat{A}_f \underline{u} \iff \underline{\rho}_f(\jump{u}) = \mat{M}^{-1} \mat{A}_f \underline{u} \,. 
	\end{equation}
Since the $W_p$ mass matrix is block diagonal by element, its inverse can be computed and stored without fill-in by simply inverting each block individually. The BR2 stabilization term can then be written as 
	\begin{equation}
	\begin{aligned}
		\sum_{f\in\Gamma_0}\int \vec{\rho}_f(\jump{u}) \cdot\vec{\rho}_f(\jump{\varphi}) \ud \x &= \sum_{f\in\Gamma_0}\underline{\rho}_f(\jump{u})^T \mat{M} \underline{\rho}_f(\jump{\varphi}) \\
		&= \sum_{f\in\Gamma_0}\underline{u}^T \mat{A}_f^T \mat{M}^{-T} \mat{M} \mat{M}^{-1} \mat{A}_f \underline{\varphi} \\
		&= \sum_{f\in\Gamma_0}\underline{u}^T \mat{A}_f^T \mat{M}^{-1} \mat{A}_f \underline{\varphi} 
	\end{aligned}
	\end{equation}
since $\mat{M}$ is symmetric. Again, since $\mat{M}^{-1}$ is block diagonal by element and the products $\mat{A}_f\underline{\varphi}$ and $\underline{u}^T \mat{A}_f^T$ are non-zero only on DOFs that share the face $f$, each argument of the sum only contributes to the DOFs that share the face $f$. Due to this, the matrix $\sum_{f\in\Gamma_0} \mat{A}_f^T \mat{M}^{-1} \mat{A}_f$ can be assembled face by face. 

Next, consider one part of the LDG stabilization term: 
	\begin{equation}
		\int \vec{\rho}_0(\jump{u}) \cdot \vec{r}_0(\jump{\E\varphi\hat{n}}) \ud \x \,. 
	\end{equation}
Let, 
	\begin{equation}
		\underline{v}^T \mat{B} \underline{\varphi} = -\int_{\Gamma_0} \avg{\vec{v}}\cdot \jump{\E\varphi\hat{n}} \ud s \,,
	\end{equation}
and further define the total interaction $W_p$ mass matrix as 
	\begin{equation}
		\underline{v}^T \mat{M}_t \underline{w} = \int \sigma_t\, \vec{v} \cdot \vec{w} \ud \x \,, 
	\end{equation}
so that $\underline{r}_0\!\paren{\jump{\E\varphi\hat{n}}} = \mat{M}_t^{-1} \mat{B} \underline{\varphi}$. In addition, define 
	\begin{equation}
		\underline{v}^T \mat{A} \underline{u} = - \int_{\Gamma_0} \avg{\vec{v}\cdot\hat{n}} \jump{u} \ud s \,, 
	\end{equation}
such that $\mat{A} = \sum_{f\in\Gamma_0} \mat{A}_f$. The LDG stabilization term under consideration is then 
	\begin{equation}
	\begin{aligned}
	 	\int \vec{\rho}_0(\jump{u})\cdot\vec{r}_0\!\paren{\jump{\E\varphi\hat{n}}} \ud \x &= \underline{\rho}_0(\jump{u})^T \mat{M} \underline{r}_0\!\paren{\jump{\E\varphi\hat{n}}} \\
		&= \underline{u}^T \mat{A}^T \mat{M}^{-T} \mat{M} \mat{M}_t^{-1} \mat{B} \underline{\varphi} \\
		&= \underline{u}^T \mat{A}^T\mat{M}_t^{-1}\mat{B} \underline{\varphi} \,.
	\end{aligned}
	\end{equation} 
Note that since the matrices $\mat{A}$ and $\mat{B}$ are not face-local, this term cannot be assembled locally. The LDG stabilization term is instead formed through matrix multiplication as $\mat{A}^T\mat{M}_t^{-1}\mat{B}$. 

\end{document}